\documentclass{amsart}
\usepackage{amssymb,epsfig}
\usepackage{color}
\usepackage{epsfig}
\usepackage{algorithm, algpseudocode}
\usepackage{subfigure}
\usepackage{multirow}
\usepackage{tikz}
\usepackage{pgfplots}
\usepackage{pgfplotstable}
\usepackage{hyperref}

\newcommand{\stab}{\textsc{ExpandStable}}
\newcommand{\solve}{\textsc{2DSOLVER}}
\newcommand{\adapt}{\textsc{Approx}}
\newcommand{\dataspace}{\textsc{UpdateDataSpace}}
\newcommand{\grid}{\mathcal{T}}

\newtheorem{theorem}{Theorem}[section]
\newtheorem{lemma}[theorem]{Lemma}

\newtheorem{assumption}[theorem]{Assumption}

\theoremstyle{definition}

\theoremstyle{remark}

\numberwithin{equation}{section}

\def\dis{\displaystyle}
\def\IR{{\mathbb R}}
\def\ov{\overline}

\def\t#1{\tilde{#1}}

\newcommand{\PP}{\mathbb{P}}

\newcommand{\R}{\mathbb R}

\newcommand{\N}{\mathbb N}
\newcommand{\Z}{\mathbb Z}

\newcommand{\bw}{{\bf w}}

\newcommand{\argmin}{\operatornamewithlimits{argmin}}
\newcommand{\argmax}{\operatornamewithlimits{argmax}}

\newcommand{\beqn}{\begin{equation}}
\newcommand{\eeqn}{\end{equation}}

\newtheorem{prop}{Proposition}[section]

\def\lll{\langle}
\def\rr{\rangle}

\newcommand\eref[1]{(\ref{#1})}
\def\e{\epsilon}

\def\cC{{\mathcal C}}

\def\cM{{\mathcal M}}

\def\cT{{\mathcal T}}

\def\cS{{\mathcal S}}
\def\cP{{\mathcal P}}

\def\cR{{\mathcal R}}

\newcommand{\Ao}{A_\circ}

\newcommand{\vb}{\vec{b}}
\newcommand{\vn}{\vec{n}}
\def\bs{\vec{s}}
\def\bn{\vec{n}}

\newcommand{\bitem}{\begin{itemize}}
\newcommand{\eitem}{\end{itemize}}
\newcommand{\benum}{\begin{enumerate}}
\newcommand{\eenum}{\end{enumerate}}
\newcommand{\beq}{\begin{equation}}
\newcommand{\eeq}{\end{equation}}

\newcommand{\norm}[1]{\|#1\|}

\def\ZZ{\mathbb{Z}}

\def\ZZ{\mathbb{Z}}

\def\cC{{\mathcal{C}}}

\def\cP{{\mathcal{P}}}
\def\cT{{\mathcal{T}}}
\def\cS{{\mathcal{S}}}

\def\cR{{\mathcal{R}}}
\def\cH{{\mathcal{H}}}

\def\cM{{\mathcal{M}}}

\def\cH{\mathcal{H}}

\newcommand{\cs}[1]{#1}

\newcommand{\wq}[1]{#1}

\newcommand{\css}[1]{#1}
\newcommand{\dww}[1]{#1}

\definecolor{dgreen}{RGB}{0,130,150}%[RGB={0,130,150}]{structure}

\title{%DRAFT \\
       Adaptive Anisotropic Petrov-Galerkin Methods
       for First Order Transport Equations
       }\thanks{This research was supported by the
        Priority Research Programme SPP1324 of the
        German Research Foundation (DFG) within the
        project
         {``Anisotropic Adaptive Discretization Concepts''}, and by the SFB-TR40, funded by DFG, and
         in part by the Excellence Initiative of the German federal and state governments.
         }
\date{\today}

\author{W. Dahmen, G. Kutyniok, W.-Q Lim, C. Schwab, and G. Welper}
\keywords{Transport equations,  shearlets, anisotropic meshes, adaptivity,
          computational complexity, best $N$-term approximation}

\begin{document}

\begin{abstract}
This paper builds on recent developments of adaptive methods for linear transport equations based on certain stable variational formulations
 of Petrov-Galerkin type. The key issues can be summarized as follows. The variational formulations allow us to employ meshes with
 cells of arbitrary aspect ratios. We develop a refinement scheme generating highly anisotropic partitions that is inspired by {\em shearlet systems}.
We establish approximation rates for $N$-term approximations from corresponding piecewise polynomials for certain compact {\em cartoon classes}
of functions.
In contrast to earlier results in a curvelet or shearlet context the
cartoon classes are concisely defined through
certain characteristic parameters
\cs{and the dependence of the approximation rates on these parameters
is made explicit here.
The approximation rate} results serve then
as a benchmark for subsequent applications to
\cs{adaptive Galerkin solvers for} transport equations.
\cs{We outline a new class of directionally adaptive, Petrov-Galerkin
discretizations for such equations. In numerical experiments, the new
algorithms track $C^2$-curved shear layers and discontinuities stably and accurately,
and realize} essentially optimal rates.
Finally, we treat {\em parameter dependent transport problems}, which
arise in kinetic models as well as in radiative transfer.
\cs{In heterogeneous media these problems feature propagation of
singularities along curved characteristics precluding, in particular,
fast marching methods based on ray-tracing.
Since now the solutions are functions of spatial variables and parameters
one has to address the curse of dimensionality.
We show computationally,
for a model parametric transport problem in heterogeneous
media in $2+1$ dimension,
that sparse tensorization of
the presently proposed spatial directionally adaptive scheme
with hierarchic collocation in ordinate space
based on a stable variational formulation high-dimensional phase space,
the curse of dimensionality can be removed when approximating
averaged bulk quantities.
}
\end{abstract}

\maketitle

%\tableofcontents

\noindent
{\bf AMS Subject Classification:}
Primary: 65N30,65J15, 65N12, 65N15\\[1mm]
%65J15, 65N12, 65N15, 35A15, 35A35, 35J60, 41A60, 46A45, 47H17
%
\noindent
{\bf Key Words:}
  Linear transport problems, $L_2$-stable Petrov-Galerkin formulations,
  $\delta$-proximality, adaptive refinements, anisotropic discretizations,
  best $N$-term approximation.

%%%%%%%%%%%%%%%%%%%%%%%%%%%%%%%%%
\setcounter{section}{0}
\section{Introduction}\label{sect1}
%%%%%%%%%%%%%%%%%%%%%%%%%%%%%%%%
%
%%%%%%%%%%%%%%%
%\subsection{Background and Motivation}
%%%%%%%%%%%%%%%%%%%%%%

Attempts to efficiently resolve anisotropic features in images and
large data sets have triggered %on the one hand
the development of {\em directional representation systems} like curvelets,
bandlets or, more recently shearlets, see
e.g. \cite{CD,KLTAMS09,KKL10a,KLL10b,KW10Compactly,Lim09b}.
These studies are motivated by
the observation that  ``{\em cartoon-like} functions'',
\cs{which are} roughly speaking piecewise smooth functions with possible jump
discontinuities along smooth curves (a precise definition will be given \cs{ahead}),
can be approximated by such systems at a better rate
than those achievable by classical {\em isotropic} wavelet systems.
Similarly, anisotropic meshes also allow one to obtain
better approximation rates than classical shape regular meshes when
the approximated functions exhibit strongly anisotropic features, see
e.g.  \cite{CM1,CoDyHeMi12,ChSuXu07,DLS}.
Both  concepts have been applied to images, or more generally,
to (I) data directly representing the  sought   object,
see e.g. \cite{BPC,DI}.
For recent advances in understanding the approximation properties of various types of
(possibly discontinuous) piecewise polynomials on such meshes,
ranging from hierarchies of nested meshes
with and without  irregular nodes to non-nested meshes obtained by local distortions,
we refer  e.g.  to \cite{CM2,Mi11}.
In addition, piecewise polynomials on anisotropic meshes have also been used
to recover (II)  an {\em implicitly given object} such as the solution of an
operator equation, see e.g. \cite{ChSuXu07,Do}.
However, in this latter context, the corresponding mesh generation is mostly
heuristic in the sense that distorting or stretching a mesh is based on
a given previously obtained approximation to the target function with
little chance to rigorously assessing the overall accuracy or to
design a refinement strategy with provable error reduction.
There seem to be only very few attempts to employ directional
representation systems in this context.
Here one should mention the work in \cite{demanet} where curvelets are used to
obtain sparse approximate representations of the wave propagation operator.
This is a somewhat different view of solving an operator equation
than the common approaches based on adaptive discretizations. It can perhaps be better compared
with computing an approximate inverse for the iterative solution of a system of equations.

The central objective of the present paper is to explore anisotropic recovery techniques in the context (II),
namely for approximating
solutions to operator equations which exhibit strongly anisotropic features, in such a way that the error assessment is more rigorously
founded than in the presently known approaches.

Specifically, we confine the discussion to a very simple model problem, namely a
linear transport equation where already discontinuities can arise, for instance, as shear layers.
The reasons for this choice are the following. On one hand, linear transport is an important
benchmark for stability issues.
It can also be viewed as the reduced model for more complex
convection-diffusion processes, highlighting the particular challenges
of \cs{vanishing} viscosity. Perhaps, more importantly,
linear transport is the core constituent of kinetic models and Boltzmann type
equations with a wide scope of applications.
Therefore, we address, in particular,  the treatment of {\em parametric} transport equations
looking for sparse representations of the solutions as functions \cs{on a high-dimensional
phase-space, i.e., the cartesian product} of the spatial variables as well as of the parameters.

The layout of the paper is as follows.
We briefly recall in Section \ref{sect:var-form} an adaptive refinement scheme from \cite{DHSW}
for unsymmetric problems, specialized here to linear transport equations.
It is based on a particular {\em well conditioned} variational formulation of the transport
equation (in both the parametric and non-parametric case) and a
strategy to realize uniformly stable {\em Petrov-Galerkin} schemes
for such problems.
The underlying functional analytic principles are very closely related to
\cite{BaMo84} and the so called DPG concepts, see e.g. \cite{DG-DPG1,DJ2}
as well as earlier least squares approaches \cite{CMMR,MMRS,MRS}.
The key feature of the present adaptive strategy is that the refinement decision
hinges on the approximation of an explicitly given {\em lifted residual} in $L_2(D)$.
This leaves considerable flexibility as to how such approximations are realized.
In particular, it accommodates, in principle,
classical finite element techniques with essentially no
constraints on the shape of the elements,
as well as directional representation systems whose frame properties are
typically known in $L_2$.

The solutions to transport equations very much resemble the cartoon model.
Depending on the right hand side,  the boundary data,
and the flow field, one typically encounters \cs{piecewise smooth solutions,
with pieces separated by a $C^2$ curve.}
%along a curve, even jump discontinuities.
This scenario will serve as our main orientation.
In Section \ref{sect:wang-Q} we develop an anisotropic approximation scheme which, on the one hand, plays well with the adaptation concept, mentioned above, and performs well in recovering cartoon like functions, viewed as a benchmark for this kind of applications.
The scheme is very much inspired by recent developments centering on shearlet systems,
see e.g. \cite{KLTAMS09,KW10Compactly}.
Moreover, we derive some implications on the performance of schemes based on anisotropic bisections of partitions comprised of
triangle and quadrilaterals.
%As it is technically somewhat involved we discuss a second alternative based on anisotropic bisections of triangles or more generally of
%simplices. This is technically somewhat simpler but exhibits slightly worse approximation rates (by a factor $\sqrt{\log N}$, $N$ the number
%of triangles).

In Section \ref{sect:numerical} we apply the shearlet inspired scheme to the solution of transport equations. The results quantify the
predicted near-optimal performance. They also demonstrate the beneficial effect of directional adaptation to the Gibbs phenomena
along jump discontinuities usually encountered.

In Section \ref{sect:parametric} we describe a numerical setup for
treating parametric problems \cs{analogous to sparse tensorization
of}  {\em discrete ordinate methods} as in \cite{KGCS2010a,KGCS2010b,Gre2014},
\cs{but with the directionally adaptive discretizations in physical space}.
%%%%%%%%%%%%%%%%%%%%%%%%%%%%
\section{Model Problem - First Order Transport Equations}\label{sect:var-form}
%%%%%%%%%%%%%%%%%%%%%%
%%%%%%%%%%%%%%%%%%%%%%%%%%%
\subsection{First Order Linear Transport Equations}\label{sect2.2}
%%%%%%%%%%%%%%%%%%
We consider the domain $D \subset \R^d$,  $d= 2,3$, denoting as usual by $\vn=\vn(x)$ the
unit outward normal at $x\in \partial D$.
Moreover, we consider velocity fields $\vb(x)$, $x\in D$, which for
simplicity will always be assumed to be differentiable, i.e. $\vec{b}(x) \in
C^1(\ov{D})^d$.
%although some of the statements remain valid under weaker assumptions.
Likewise $c(x) \in C^0(\ov{D})$ will serve as the reaction term in the {\em
first order transport equation}
\begin{align} \Ao u : = \vb \cdot \nabla u + cu & = \;\mbox{$f_\circ$ in
$D$}\,, \label{2.1} \quad u  = \; \mbox{$g$ on $\partial D_-$}\,, \end{align}
where ${\partial D}_\pm :=  \{x \in \partial D: \; \pm \vec{b} (x) \cdot \vec{n}(x) >
0\}$ denotes the {\em inflow, outflow boundary}, respectively, and $f_\circ \in L_2(D)$.
%in the following partition of ${\partial D}=\partial D$
%
Furthermore, to simplify the exposition we shall always assume that
\begin{equation} \label{coerccond} c -  \frac 12 \nabla\cdot  \vb \geq c_0 >
0 \quad \mbox{in}\,\, D \eeqn
holds.   In order to make use of the approximation results for shearlet inspired
approximation methods in Section \ref{sect:wang-Q}, we introduce a variational
formulation for which the solution $u$ is a function in $L_2(D)$.  To this end,
we first multiply \eqref{2.1} by a test function and use integration by parts
so that there are no derivatives applied to $u$.
Following \cite{DHSW}, one can show that the resulting bilinear form
\beqn \label{auv} a(w,v):=   \dis\int_D \,w(- \vec{b} \cdot {\nabla} v +  (c-
{\nabla} \cdot \vec{b})v)  \;dx, \eeqn
%
%which in view of \eqref{Aest}
is bounded on $L_2(D)\times W_0(-\vb,D)$, where
\beqn
\label{W0}
W_0(\mp \vb,D)
:=
{\rm clos}_{\|\cdot\|_{W(\vb,D)}}\{ v\in C^1(D)\cap C(\ov{D}),\,
v\mid_{{\partial D}_{\pm}}\equiv 0\}.
\eeqn
and
\beqn
\label{graph}
\|v\|_{W(\vb,D)}:=
\left( \|v\|_{L_2(D)}^2 + \int_D |\vb\cdot\nabla v|^2\,dx\right)^{1/2}.
\eeqn
Moreover, the trace $\gamma_-: W_0(-\vb,D) \to L_2({\partial D}_-,\omega)$
on the inflow boundary ${\partial D}_-$ exists and is bounded,
where the latter space is endowed with the weighted
$L_2$-norm $\|g\|^2_{L_2({\partial D}_\pm,\omega)} = \int_{{\partial D}_\pm} |g|^2\omega ds$
with weight $\omega = |\vb \cdot n|$
see, eg. \cite{Bardos}.
Thus,  the linear functional $f(\cdot)$, given by
\beqn
\label{ell}
f(v)
:= \int_D f_\circ vdx
%{}_{(W_0(-\vb,D))'}\lll f_\circ,v\rr_{W_0(-\vb,D)}
+
\int_{{\partial D}_-}g\gamma_-(v)|\vb\cdot\vn|ds,
\eeqn
belongs to  $(W_0(\vb,D))'$ and the variational problem
\beqn
\label{varprob}
a(u,v) = f(v),\quad v\in W_0(-\vb,D),
\eeqn
admits
a unique solution in $L_2(D)$ which, when regular enough,
coincides with the classical solution of \eqref{2.1}, see \cite[Theorem 2.2]{DHSW}.
Taking the boundary conditions on the test functions into account, the additional boundary integral in the right hand side  results from   integration by parts, in analogy  to
the weak imposition of Neumann boundary conditions for elliptic problems.

Defining for $v\in W_0(-\vb,D)$,
the adjoint $A^*$ of $A$ by $a(w,v)=\lll w,A^*v\rr$, $(w,v)\in L_2(D) \times W_0(-\vb,D)$
the quantity $\|v\|_Y:= \|A^* v\|_{L_2(D)}$ is an equivalent norm on $W_0(-\vb,D)$.
Moreover, $A: L_2(D) \to (W_0(\vb,D))'$ and $A^*: W_0(\vb,D)\to L_2(D)$
are {\em isometries}  (see \cite[Proposition 4.1]{DHSW})
\beqn
\label{isos}
\begin{aligned}
  \|A\|_{L(L_2(D),W_0(\vb,D)')} & = \|A^{-1}\|_{L(W_0(\vb,D)', L_2(D))} = 1,
\\
  \|A^*\|_{L(W_0(\vb,D),L_2(D))} & = \|(A^*)^{-1}\|_{L(L_2(D), W_0(\vb,D))} = 1\;.
\end{aligned}
\eeqn
%
%%%%%%%%%%%%%%%%%%%%%%%%
\subsection{Parametric Transport Problems}\label{sect2.3}
%%%%%%%%%%%%%%%%%%%%%%%%%%%%
%
 {In \eqref{2.1} $\vb(x)$ is a fixed single given flow-field.
Here we are also interested in the case when the convection field
$\vb(x,\bs)$ depends on a {\em parameter} that may range over some parameter set $\cS$.
%We consider the following model problem,
%as core constituent of radiative transfer models,
%%%%%%%%%%%%%%%%%%%%%%%
%\subsubsection{Variational Formulation}\label{ssec:vf}
%%%%%%%%%%%%%%%%%%%%%%%%%%
Specifically, we now consider  the parametric family of  transport problems
%\eqref{rad1.1}, ie.,
%
\beq\label{eq:rads}
\begin{array}{rcl}
A_\circ u(x,\bs) &=&  \vb(x,\bs) \cdot \nabla u(x,\bs) + c(x) u(x,\bs) =  f(x),\;\; x\in D\subset \IR^2\;,
\\
u(x,\bs)  &=& g(x,\bs)\;,\;\; x\in {\partial D}_-(\bs) ,
\end{array}
\eeq
where  the parameter dependent inflow-boundary  $ {\partial D}_-(\bs) $ % in \eqref{eq:Gammapm}
is now given by
\beq\label{eq:Gammapms}
{\partial D}_\pm(\bs)
:=
\{ x\in \partial D : \pm \vb(x,\bs) \cdot \vec{n}(x) > 0 \},\;\; \bs\in \cS
\;.
\eeq
For the special case $\vb(x,\bs) = \bs$ the parametric transport problems \eqref{eq:rads} form
a core constituent of {\em radiative transfer} models. While in the latter special case singularities of the parametric solution
propagate along {\em straight lines} so that efficient directional approximation methods
such as {\rm ray tracing} are applicable, the present approach covers parameter dependent families of  {\em curved} characteristics.
%
%\beqn
%\label{rad1.1}
%\begin{array}{rcl}
%\Ao u(x,\bs) = \bs\cdot \nabla u(x,\bs)+\kappa (x)u(x,\bs)&=&%\kappa(x)
%f_\circ(x),\quad x\in D\subset \R^d,\,\, d=2,3,
%\\
%u(x,\bs)&=& g(x,\bs),\, x\in {\partial D}_-(\bs),
%\end{array}
%\eeqn
%
Now the solution $u$ depends also on the prameter $\bs$ in the convection field $\vb$.
%which, however, varies over a set  $\cS$ of directions $\bs$.
Thus, for instance, when $\cS=S^2$, the unit $2-$sphere, $u$
 could be considered as a function of five variables,
namely  of $d=3$ spatial variables  $x$ and  of parameters  $\bs$ from the sphere $\cS$.
We always assume that $\vb$ depends smoothly on the parameter $\bs$ while
for each $\bs\in \cS$, as a function of the spatial variable $x$, we continue to impose the same conditions as
in the preceding section.
%Clearly, the  in- and outflow  boundary now depends on $\bs$, i.e.,
%
%\beqn\label{eq:Gammapm}
%\gk{\partial D}_\pm(\bs):= \{x\in \partial D: \mp \bs\cdot {\bf n}(x)<0\},\qquad \bs\in \cS\;.
%\eeqn
%
The absorption coefficient $c\in L_\infty(D)$ will always be assumed to be nonnegative in
 the physical domain $D$ and to satisfy \eref{coerccond}.
 }

Following \cite{DHSW} we can resort to similar concepts as above to
obtain a variational formulation for
\eqref{eq:rads} over $2d-1$-dimensional phase domain
$D\times \cS$.
To that end, let
\beqn
\label{gammatotal}
\partial D
:=
\partial D \times \cS,\quad {\partial D}_\pm
:=
\{(x,\bs)\in \partial D : \mp \vb(x,\bs)\cdot \bn(x)<0\},\quad
{\partial D}_0:= {\partial D}\setminus ({\partial D}_-\cup {\partial D}_+),
\eeqn
and  denote as before
$
(v,w):=(v,w)_{D\times \cS}=\int_{D\times \cS}v(x,\bs)w(x,\bs)dxd\bs
$
where, however, $d\bs$ is for simplicity the
normalized Haar measure on $\cS$, i.e.
$
\int_{\cS}d\bs =1.
$
Defining  the Hilbert space
\beqn\label{eq:defWDS}
W(D\times \cS)
:=
\{ v\in L_2(D\times \cS):
   \int_{\cS\times D}|\vb(\bs)\cdot \nabla v|^2 dx d\bs <\infty\}
\eeqn
(in the sense of distributions where the gradient $\nabla$
 always refers to the $x$-variable in $D$),
 endowed with the norm  $\|v\|_{W(D\times \cS)}$ given by
 \beqn\label{eq:normWDS}
\|v\|_{W(D\times \cS)}^2
:=
\|v\|_{L_2(D\times \cS)}^{2} + \int_{\cS\times D}|\vb(\bs)\cdot \nabla v|^2 dx d\bs,
\eeqn
the counterpart to \eqref{W0} is given by
\beqn\label{eq:Wpm0}
W^\pm_0(D\times \cS)
:=
{\rm clos}_{\|\cdot\|_{W(D\times \cS)}} \{v\in C(\cS,C^1(D)):
             v|_{{\partial D}_\pm}\equiv 0\}
\eeqn
which is again a Hilbert space under the norm $\|v\|_{W(D\times \cS)}$.
It is shown in \cite{DHSW} that the problem
\beqn
\label{rad1}
a(u,v)= f(v)
\eeqn
with
\begin{align*}
  a(u,v) & := \int_{D \times \cS} u(- \vb(\bs) \cdot {\nabla} v + (c-\nabla\cdot \vb) v )  \,dx \, d \bs,
  \\
  f(v) & := %{}_{W^+_0(D\times \cS)'}\lll f_\circ , v\rr_{W^+_0(D\times \cS)}
  \lll f_\circ , v\rr
   +\int_{{\partial D}_-}g\gamma_-(v)|\bs\cdot\bn|d\cs{\bs},
                        \quad v\in W^+_0(D\times \cS),
\end{align*}
has a unique solution in  $L_2(D \times \cS)$  %\Omega used later for cartoon
which, when sufficiently regular, agrees
for  almost every $\bs\in \cS$ with the corresponding classical solution of \eqref{eq:rads}. Here $\lll f_\circ , v\rr$ is the dual pairing
on $W^+_0(D\times \cS)' \times W^+_0(D\times \cS)$.
%\subsection{Some Numerical Observations}
%%%%%%%%%%%%%%%%%%%%%%%%%%%%%%
%The nature of the solutions to the above transport problems is best demonstrated for problem \eqref{2.1}.
%Specifically, we consider the following cases for $d=2$.\\

%\noindent
%{\bf (I)}

%%%%%%%%%%%%%%%%%%%%%%%%%%%%
\section{Adaptive Scheme}\label{sect:ad}
%%%%%%%%%%%%%%%%%%%%%%%%%%%%%%%
Both variational formulations of Sections \ref{sect2.2} and \ref{sect2.3}
are instances of the generic variational problem of
finding $u \in X$ such that
\beqn
\label{A}
a(u,v) = \lll f,v\rr ,\quad v\in Y,
\eeqn
which is either \eqref{varprob} or  \eqref{rad1}.
The respective spaces are $X = L_2(D)$ or
$X = L_2(D\times \cS)$ and $Y = (W_0(-\vb,D))$ or $Y = (W^+_0(D\times \cS))$.

The relevance of the
isometry properties \eqref{isos} %and \eqref{iso2}
lies in the fact that errors in $X$ are equal to the corresponding residuals in $Y'$,
i.e.,
\beqn
\label{best}
\|u - v\|_X = \| f -A v\|_{Y'}\;.
\eeqn
We use this equality of errors and residuals to deal with the following
two tasks:

\begin{itemize}
  \item[(I)] Devise a numerical scheme which,
  for a given trial space $X_h\subset X$, yields an approximation $u_h\in X_h$ that is near best,
  ie.,
  \beqn
  \label{nearbest}
  \|u-u_h\|_X \leq C \inf_{w\in X_h}\| u-w\|_X;
  \eeqn

  \item[(II)] Estimate errors via the residual $\| f -A u_h\|_{Y'}$.
\end{itemize}

As for (I), the main idea, which \cs{also} plays a central role in the recent
developments of DPG schemes \cite{DG-DPG1,DJ2},
is to construct uniformly stable Petrov-Galerkin schemes.
\cs{D}ue to \eqref{isos}, $Y_h:= A^{-*}X_h$ would be the {\em ideal test space}
for a given trial space $X_h$, in the sense that
\beqn
\label{infsup1}
\inf_{w\in X_h}\sup_{v\in Y_h}\frac{a(w,v)}{\|w\|_X \|v\|_Y}=1,
\eeqn
but \cs{$Y_h$} is practically inaccessible.
%\gw{Note that we can choose the constant $1$ on the right hand side because of the isometry properties \eqref{isos} and \eqref{iso2}.}
 One therefore settles for an approximate test space $Y_h$
which is sufficiently close to ideal to perserve uniform (with respect to $h$) stability.
This has been quantified in \cite{DHSW} through the notion of
{\em $\delta$-proximality}:
a finite dimensional space $Z_h\subset Y$ is called
$\delta$-proximal for $X_h\subset X$ for some $\delta \in (0,1)$ if the
$Y$-orthogonal projection $P_{Y,Z_h}: Y \to Z_h$ satisfies
\beqn
\label{delta}
\|y - P_{Y,Z_h}y\|_Y \leq \delta \|y\|_Y,\quad   y\in Y_h = A^{-*}X_h
\;.
\eeqn
Note that for any $y=A^{-*}w$, $w\in X_h$,
the projection $P_{Y,Z_h} y= \t y$ is given by the Galerkin projection
\beqn
\label{deltaGal}
(A^*\t y,A^* z)_{X'}= (w, A^*z)_{X'},\quad v\in Z_h,
\eeqn
and hence computable at a cost depending on ${\rm dim}\, Z_h$. When $Z_h$ is $\delta$-proximal for $X_h$
the {\em test space} $\t Y_h := P_{Y,Z_h}(A^{-*}X_h)$
turns out to give rise to uniformly stable Petrov-Galerkin schemes
in a sense to be made precise in a moment.
Before, we wish to point out two different ways of using this fact:
in the context of discontinuous Galerkin formulations with suitably adjusted mesh-dependent norms,
it is actually possible to compute  a basis
%CS2all: removed ``explicit'', as that's implied by the word ``compute''
for $\t Y_h$ at a computational cost that under suitable circumstances
stays essentially proportional to the dimension of $X_h$  \cite{DG-DPG1,DJ2}.
For conforming discretizations determining each test basis function would require
solving a problem of essentially the original size.
However, there is a way of realizing the Petrov-Galerkin scheme without
computing the test basis functions explicitly while keeping the computational
expense again essentially proportional to ${\rm dim}\, X_h$, \cite{DHSW}.
To this end, using the optimal test functions in
the Petrov-Galerkin scheme $a(u_h, A^{-*}v_h) = f(A^{-*}v_h)$
and, using that $A^{-*} = (AA^*)^{-1}A$,
we obtain
\begin{align*}
  \langle f-A u_h, A^{-*}v_h \rangle = \langle (AA^*)^{-1}(f-A u_h), A v_h \rangle & = 0 & v_h & \in X_h,
\end{align*}
where $\langle \cdot, \cdot \rangle$ is the $Y', X$ dual pairing.
Singling out $r_h = (AA^*)^{-1}(f-A u_h)$, or rather an approximation thereof,
as a new variable, we obtain the saddle point problem
\beqn
\label{system-disc}
\begin{array}{lcll}
 (A^* r_{h} ,A^* z_h)_{X'}  + a( u_{h} ,z_h)
  & = & \lll f,z_h\rr, & z_h \in Z_h,\\
a(v_h,r_h) %\lll \bo_\mu^* u_f, q_f\rr
&=& 0,& v_h \in X_h.
\end{array}
\eeqn
The error estimates from \cite{DHSW}
can be summarized as follows.
\begin{prop}
\label{prop:stab2}
Assume that $Z_h\subset Y$ is $\delta$-proximal for $X_h\subset X$ with some $\delta\in (0,1)$.
Then
the solution component $u_h$ of the saddle point problem \eqref{system-disc}
solves the Petrov-Galerkin problem
\beqn
\label{PG}
a(u_h,v)=\lll f,v\rr, \quad v\in P_{Y,Z_h}(A^{-*}X_h)\;.
\eeqn
Moreover, one has
\beqn
\label{BAP1}
\|u- u_h\|_{X}\leq \frac{1}{1-\delta}\inf_{w\in X_h}\|u-w\|_{X},
\eeqn
and
\beqn
\label{BAP2}
\|u- u_h\|_{X}+ \| r_h \|_{Y}\leq \frac{2 }{1-\delta}
\inf_{w\in X_h}\|u-w\|_{X}.
\eeqn
Finally, one has
\beqn
\label{inf-sup-0}
\inf_{v_h \in X_h} \sup_{z_h \in Z_h}\frac{a(v_h , z_h)}{\|z_h\|_{Y}\|v_h \|_{X}} \geq  {\sqrt{1-\delta^2}} .
\eeqn
\end{prop}

The benefit of the above saddle point formulation is not only
\cs{that it allows to bypass the computation} of the test basis functions
but \cs{also that provides} an {\em error estimator}
based on the {\em lifted residual} $r_h=r_h(u_h,f)$ defined by \eqref{system-disc}.
In fact,
by construction $r_h$ is an approximation of $(AA^*)^{-1}(f-A u_h)$ where the additional $(AA^*)^{-1}$ makes sure that $r_h$ is measured in the localizable $Y$-norm
instead of the computationally  generally inaccessible $Y'$-norm of the residual.
It is shown in \cite{DHSW} that when $Z_h\subset Y$ is even $\delta$-proximal
for $X_h + A^{-1}F_h$  with $F_h\subset Y'$ being
a computationally accessible subspace, one has
\beqn
\label{rh}
(1-\delta)\|f_h - A v_h\|_{Y'} \leq \|r_h(v_h,f_h)\|_Y \leq \|f_h - A v_h\|_{Y'},
\eeqn
and hence a tight bound on the residual $\|f_h - A v_h\|_{Y'}$ obtained for the approximation $f_h$ of $f$.
Generally, a vector $r_h(v_h, f)$, which is calculated  only to finite precision  %byfinite means
on a computer,
cannot accurately represent every possible error for {\em any} right hand side $f$
in an infinite dimensional space but only for the prortion of $f$ that  can be captured through the current discretization.
Therefore the approximation $f_h \in F_h$ represents the components of $f$ visible to the error estimator.
Analogously to finite element error estimators for elliptic problems,
$f-f_h$ can be regarded  as %a
data oscillation error.
Since $ \|r_h\|_Y = \|A^* r_h\|_X$,
the current error of the Petrov-Galerkin approximation $u_h$
is tightly estimated from below and above by the quantity
$\|A^*r_h(u_h,f_h)\|_X$.

As shown in \cite{DHSW},
the saddle point problem \eqref{system-disc} is efficiently solved approximately
by the Uzawa iteration
\begin{equation}\label{uzawa}
\begin{array}{rcll}
  (A^*  r^k_h, A^* z_h) & = & (f-A u^k_h,z_h), & \quad z_h \in Z_h,
\\
  (u^{k+1}_h,v_h) & = & (u_h^k,v_h) + (A^*  r^k_h,v_h), & \quad v_h \in X_h\;.
\end{array}
\end{equation}
The approximations $u_h^k$ to $u_h$ are obtained in \eqref{uzawa}
through simple updates while the lifted residuals $r_h^k$ result
from solving a {\em symmetric positive definite system} in $Z_h$.
See \cite[Section~4.2]{DHSW} for more details and a proof of
the following proposition.
\begin{prop}
\label{prop:reduction}
Assume that $\delta$-proximality \eqref{delta} is satisfied.
Then the iterates $u_h^{k}$ of the Uzawa iteration converge
to the solution $u_h$ of the saddle point problem \eqref{system-disc}
and
\begin{equation*}
  \|u_h-u_h^{k+1}\|_X \le \delta \|u_h - u_h^k\|_X.
  \end{equation*}
\end{prop}
\cs{
An adaptive algorithm based on the Uzawa iteration
is based on the following subroutines.
}
\begin{description}
\item[\rm $\stab$]
For a given space $X_h$ return a space $Z_h$ such that
$\delta$-proximality \eqref{delta} is satisfied.
\item[\rm $\adapt$]
For a function $w \in X$, error tolerance $\epsilon$ and a space $X_h$ return
an approximation $w_a \in \hat X_h$ and a new space $\hat X_h \supset X_h$ such that
\begin{equation*}
\|w - w_a\|_X \le \epsilon\;.
\end{equation*}
\item[\rm $\dataspace$]
For a function $f \in Y'$, error tolerance $\epsilon$ and spaces $F_h$ and $X_h$,
return an approximation $f_a \in \hat F_h + A X_h$ and
\cs{enlarged} space $\hat F_h \supset F_h$ such that
    \begin{equation*}
      \|f - f_a\|_{Y'} \le \epsilon\;.
    \end{equation*}
\end{description}

According to the second row of the Uzawa iteration
$u^{k+1}_h = u_h^k + A^*r_h^k$ the solution $u_h^{k+1}$ is %the old one
updated by a residual term.
Since $r_h^k \in Z_h$ for a stable space $Z_h \supset X_h$ and $A_h$ is a differential operator,
we generally do not have $A^* r_h^k \in X_h$.
In the standard Uzawa iteration this function is projected onto $X_h$.
Instead,  here we use the adaptive approximation $\adapt$,
giving rise to more accurate updates and better \cs{space adaptation of} $X_h$
to approximate the solution $u$.
The precise formulation of the adaptive refinement scheme is given in \cite{DHSW}
and summarized in Algorithm \ref{alg:uzawa}.

\begin{algorithm}[htb]
  \caption{Adaptive Uzawa iteration}
  \label{alg:uzawa}
  \begin{algorithmic}[1]
    \State Initialization: Choose a target accuracy $\epsilon$ and an initial trial space $X_h$. Set the initial guess, initial error bound and test space
    \begin{align*}
      u_a & = 0, & errbound & =\|f\|_{Y'}, & Z_h = \stab (X_h),
    \end{align*}
    respectively. Choose parameters $\rho, \eta, \alpha, \in (0,1)$ and $K \in \N$.
    \While{$errbound> \epsilon$}
      \State Compute $\hat u_a$, $\hat r_a$ as the result of $K$ Uzawa iterations with initial value $u_a, r_a$ and right hand side $f_a$.
      \State $(update, \hat X_h) = \adapt(A^* \hat r_a, \eta\|A^* \hat r_a\|_X, X_h)$
      \State $X_h := X_h + \hat X_h$
      \State $(f_a, F_h) := \dataspace(f, \alpha \rho \, errbound, F_h, X_h)$
      \State $Y_h := \stab(X_h+A^{-1}F_h)$
      \State $u_a := u_a + update$, $errbound := \rho \, errbound$
     \EndWhile
  \end{algorithmic}
\end{algorithm}

It is shown in \cite[Proposition 4.7]{DHSW} that Algorithm \ref{alg:uzawa}
terminates after finitely many steps and outputs an approximate
solution $u_a$ satisfying $\|u-u_a\|_X \leq \e$.
Due to the fixed error reduction per step in Proposition \ref{prop:reduction},
a fixed number of Uzawa iterations $K$ in each loop of the algorithm is sufficient.
%However, using information \cs{obtained from}
%the error estimators, one can also stop the iteration earlier.

In the envisaged applications we either have
$X=L_2(D)$ or $X= L_2(D\times \cS)$.
Thus the actual adaptive refinement $\adapt$
is based on an $L_2$-approximation of an explicitly given function $A^*\hat r_a \in X$.
This leaves some flexibility in the choice of the trial spaces $X_h$
and their adaptive enlargement in $\adapt$.
The approximation in $L_2$ is by no means constrained by any
{\em shape regularity} conditions \cs{in the subdivisions of $D$
used for constructing the trial space}.
In particular, we are free to employ highly anisotropic approximation
systems, as e.g. the shearlet inspired implementation of $\adapt$ in
Section \ref{sect:aniso} below.
For the actual enlargement of $X_h$ in $\adapt$,
one could think, for instance, of an $L_2$-projection
$P_{X,\bar X_h}(A^* \hat r_a)$ onto a  fixed, but sufficiently large number of levels of
uniform refinement $\bar X_h$ of $X_h$,
followed by a coarsening step which, in principle, would even
allow one to formulate optimality statements.
As an alternative one could consider extremal problems
$$
\inf_{\| A^* \hat r_a - w\|_X\leq \eta \e; w\in \bar X_h}\|\bw\|_{\ell_1},
$$
where $\bw$ is the coefficient vector of $w$.
Moreover, in many settings an orthonormal basis of $\bar X_h$ is available so that
a greedy enlargement of $X_h$ is efficient and yields a best $N$-term approximation
of $A^* \hat r_a$ in $\bar X_h$.
%%%%%%%%%%%%%%%%%%%%%
\section{Adaptive Anisotropic Approximation}\label{sect:aniso}
%%%%%%%%%%%%%%%%%%%%%%%%%%%%%
The nature of solutions to  linear transport equations such as \eqref{2.1}
is illustrated by considering the following instances  of problem \eqref{2.1}.
\\
 (I): When the right hand side $f_\circ$ has a jump discontinuity along a smooth curve,  the boundary conditions
 $g$ are zero, and the convection field $\vb$ is $C^1$, the solution $u$ is only Lipschitz continuous with a
 jump of first order derivatives in directions not parallel to $\vb$.
 If the discontinuity of $f_\circ$ is parallel to $\vb$, the solution even exhibits a jump
  discontinuity.

 (II): When $f_\circ$ is smooth, say constant, $\vb$ is $C^1$, and the boundary data
  is piecewise smooth with a jump discontinuity,
 the solution  will exhibit a jump discontinuity along a $C^2$ curve with tangents $\vb$.

 (III): When both, right hand side and boundary data are discontinuous both above effects superimpose.

 At any rate, the solutions $u$ in the above scenarios are piecewise smooth on regions separated by a $C^2$ curve.
 In the context of image processing such functions are usually referred to as {\em cartoon-like} and serve as
 benchmarks for image compression schemes that are able to economically encode edges and curve-like discontinuities.
  In what follows we do the same in the context of solving transport equations. The methods employed
  in image processing cannot be directly applied in this context.
One reason is the more delicate role of boundary conditions
  in the present context. The second one is more subtle. Since we approximate solutions in $L_2$ one could, in principle,
  employ directional representation systems such as shearlets and curvelets because they form frames in $L_2(\R^d)$.
But even setting aside the issue of boundary conditions (which would certainly destroy the tightness of such frames),
  one would have to determine suitable test frames whose elements are obtained by solving a global problem which renders
  resulting schemes inefficient.
%\gw{(As far as I remember that was our first idea before we had the saddle point formulation. Thus it is not clear that we really have to compute the test frame explicitly, just as we don't compute test bases explicitly.)}
Moreover, quadrature becomes a major issue since elements of directional
representation systems do not stem from nested hierarchies of multiresolution sequences.
We shall therefore discuss modifications that work in the present context as well.
%%%%%%%%%%%%%%%%%%%%%%%%%%%%
\subsection{Isotropic Schemes}\label{sect:iso}
%%%%%%%%%%%%%%%%%%%%%%%%%%%%
We begin recalling briefly some typical results for numerical approximations to solutions of transport equations
using uniform and adaptive but {\em isotropic} meshes in both cases (I) and (II) described
at the beginning of Section \ref{sect:aniso}.

To this end  the procedures $\stab$ and $\adapt$ need to be specified.
Since $X_h$ is an isotropic finite element space, we can define $\stab(X_h)$
simply as the finite element space obtained by refining each cell of the
\cs{partition} of $X_h$ once or by increasing the polynomial degree of the finite elements.
%by one or two orders.

To define $\adapt(w, \epsilon, X_h)$,
the algorithm has to find an adaptive approximation of a known function
$w$ up to  {prescribed accuracy $\epsilon>0$.
 A straightforward strategy to this end is as follows:
first compute a best $L_2$ approximation $w_h \in X_h$ of $w$. The local errors are
\begin{equation*}
 e(T) = \|w-w_h\|_{L_2(T)}^2 = \int_T |w-w_h|^2 \text{d} x
\end{equation*}
for each cell $T$ in the \cs{partition} $\grid_h$ of $X_h$.
Next, refine cells in a subset $\mathcal{S} \subset \grid_h$
that are selected by a bulk criterion
\begin{equation*}
 \sum_{T \in \mathcal{S}} e(T) \ge \theta \sum_{T \in \grid_h} e(T)
\end{equation*}
for some constant $0 < \theta < 1$.
This is a first, rough realization for
the overall algorithm and numerical experiments which are provided
below.
To guarantee the error bound $\epsilon>0$, we can explicitly compute the best
$L_2(D)$ approximation $w_{ref}$ on the refined finite element space and the error
$\|w - w_{ref}\|_{L_2}$.
If the error bound $\epsilon$ is met, we stop.
If it is not met, we recursively refine the finite element space again.

%%%%%%%%%%%%%%%%%%%%%%%%%%%%%%%%%%%%%%%%%%%%%%%55
\subsection{Benchmark Class}\label{sect2}
%%%%%%%%%%%%%%%%%%%%%%%%%%%%%%%%%%%%%%%%%%%%%%%%%%

Let $D = (0,1)^2$ and assume that  {$\mathcal{O}$ is the class of  nonempty simply connected domains $\Omega \subset \R^2$
  with a  $C^2$-boundary}.
%\todo{[Further conditions, such as connectedness?]}
Let ${\rm curv}(\partial\Omega )$
denote the curvature of $\partial\Omega\cap D$ and $|\partial \Omega\cap D|$ its length in $D$. Finally,
let $C_a(x):= \{x'\in \R^2: |x-x'|<a\}$ denote the disc of radius $a$ around $x$ and define the
{\em separation width} of $\Omega$
\beqn
\label{sepwidth}
\rho=\rho(\Omega):= \sup\,\{\zeta >0 : \forall\, x\in \partial\Omega,\,\,\partial\Omega \cap B_\zeta(x) \mbox{ contains a single
connected arc of } \partial\Omega\}.
\eeqn
%\todo{[provide definition of {\rm curv}]}.
Then, we consider the following model   class of {\em cartoon-like functions} on $D = (0,1)^2$,  defined by
\begin{eqnarray}\label{cartoon}
{\cC}( \kappa, L,M,\omega) &:=& \{ f_1\chi_{\Omega } + f_2\chi_{D\setminus\Omega}: \,\,
\Omega\subset D,\,  \Omega\in \mathcal{O},\,  |\partial\Omega \cap D|\leq L, \partial\Omega \cap D \in C^2,
\nonumber \\
&& \rho(\Omega)\ge \omega,\,  {\rm curv}(\partial\Omega )\leq \kappa,  \|f_i^{(l)}\|_{L_\infty(D)}\leq M, \, l \leq 2,\, \, i=1,2\},
\end{eqnarray}
(where the parameters $\kappa,L$ are not mutually independent).
We derive approximation error bounds for anisotropic refinements  for the
simple scenario \eqref{cartoon} making an effort to exhibit the dependence of involved constants on the
class parameters $\kappa,L,M,\omega$.
We shall indicate subsequently how the results extend to
natural extensions of this class, e.g. admitting several components or curves with selfintersections
and arbitrarily small separation width.

In the area of imaging sciences, this
model has become a well accepted benchmark for sparse approximation, see, e.g., \cite{D2001}.
One prominent system, which does provide (near-) optimal sparse approximations
for such classes are compactly supported
shearlet systems for $L^2(\mathbb{R}^2)$ \cite{KKL10a,KW10Compactly}.
Our main point here is that, as the discussion at the
beginning of this section shows,  such cartoons
also exhibit similar features as solutions to transport problems.

Unfortunately, even compactly supported shearlets do not comply well with quadrature  and boundary adaptation tasks
faced in variational methods for PDEs.
%\gk{The approach recently presented in \cite{GKMP15} is capable of coping with
%boundary adaption, is, however, not adapted to quadrature tasks.}
We are therefore interested in generating locally refinable anisotropic partitions
for which corresponding piecewise polynomial approximations realize the favorable near-optimal approximation rates for
cartoon functions achieved by shearlet systems. Unfortunately, as shown in \cite[Chapter 9.3]{Mi11}, simple triangular
bisections connecting the midpoint of an edge to the opposite vertex is not sufficient for warranting such rates, see
\cite{ChSuXu07,CoDyHeMi12} for related work. In fact, a key feature would be to realize a ``parabolic scaling law''
similar to the shearlet setting, see e.g. \cite{KLTAMS09}.
%\todo{[GKWQL: References to ``parabolic scaling here?]}
By this we mean a  sufficiently efficient
directional resolution by anisotropic cells whose
width scales like the square of the diameter.
To achieve this we consider partitions comprising triangles {\em and}
\cs{parallelograms} pointed out to us in \cite{priv-com}.
%%%%%%%%%%%%%%%%%%%%%%%%%%%%%%%%%%%%%%
\newcommand{\dI}{\Delta}
\newcommand{\TT}{\mathfrak{T}}

%In the following subsections, we will describe two refinement schemes, the
%first consisting of simplicial bisections
%and the second being inspired by the shearlet framework.
%It will turn out that the second approach indeed improves the
%approximation error by a $\log$-factor.
%Our numerical experiments will be based upon this approach.

\subsection{Approximation Rates}\label{sect:wang-Q}
%%%%%%%%%%%%%%%%%%%%%%%%%%%%%%%%
The central objective of this section is to describe and analyze a
nonlinear piecewise polynomial approximation scheme on {\em anisotropic} meshes
\cs{in physical space $D$}
for elements from \cs{the} cartoon class \eqref{cartoon}.
Specifically, we establish approximation rates
matching those obtained for  {\em shearlet systems} (see \cite{KKL10a,KW10Compactly}
or the survey \cite{KL12}).
In fact, the scheme is directly inspired by shearlet concepts in that a key element is the ability to realize the parabolic scaling
law for a sufficiently rapid directional adaptation.
The latter point is the main distinction \cs{of the present work from}
technically simpler schemes based on bisecting triangles.
The envisaged error bounds, \cs{which are proved} using {\em full knowledge}
about the approximated function, will serve later as
a benchmark for the performance of our {\em adaptive}
scheme that is to approximate the {\em unknown} solution of a transport problem.

As hinted at before, we  wish to analyze the approximation power of piecewise polynomial approximations
based on successively constructed anisotropic meshes,
utilizing the operations of parabolic scaling and shearing.
The matrices associated to those two operations are defined by
\begin{equation}\label{eq:shear1}
D_{j} = \begin{pmatrix}
2^j & 0 \\ 0 & 2^{\lfloor j/2 \rfloor}
\end{pmatrix}
\quad \text{and} \quad
S_k = \begin{pmatrix}
1 & k \\ 0 & 1
\end{pmatrix},
\qquad j \ge 0, k \in \mathbb{Z}.
\end{equation}
The relation to shearlet systems can be seen by observing that shearlets
are a function system which is -- except for a low frequency component -- generated
by application of parabolic scaling, shearing and translation operators to
 at least two ``mother functions''.
Regarding parallelograms as ``prototype supports'' of compactly supported shearlets  we consider next partitions comprised of
parallelograms and triangles only. Specifically, we set
\beqn
\label{Sigman}
\Sigma_N := \bigcup\{ \PP_1(\cT):  \cT \, \mbox{a partition of $D$ into parallelograms
and triangles with}\, \#(\cT)\leq N\},
\eeqn
where for any partition $\cT$ of $D$ we denote by $\PP_1(\cT)$
the space of \cs{discontinuous,} piecewise affine functions \cs{on} the partition $\cT$.
%In this sense  the parallelograms in such anisotropic meshes
%can be regarded as potential supports of elements of a compactly supported shearlet system.

%We analyze the approximation properties of $\sigma_n$ first for the following restricted {\em Horizon model} of cartoons:
%\beqn
%\label{Horizon}
%\cC_H(\kappa,L,M,\omega):= \{f_1\chi_\Omega + f_2\chi_{D\setminus\Omega}\in \cC(\kappa,L,M,\omega): \Omega \mbox{ has the form
%%\eref{eq:boundary}} \}
%\eeqn
%where
%
%\begin{equation}\label{eq:boundary}
%\Omega = \{(x_1,x_2) \in D : x_1 \leq E(x_2)\},\quad \partial\Omega\cap D = \{(E(x_2),x_2): x_2\in (0,1)\},
%\end{equation}
%
%and %$D = (0,1)^2$ and
% {$E \in C^2([0,1])$ satisfies
%\begin{equation}
%\label{Eprop}
% \|E'\|_{L_\infty(D)} \leq 2,\quad \|E''\|_{L_\infty(D)} \leq \kappa.
% \end{equation}
%Thus,  $\gk{\partial D} = \partial \Omega \cap D $
%is the discontinuity curve referred to in \eqref{cartoon}.

The main result shows that approximations by elements of $\Sigma_N$
%based on this refinement scheme
 realize (and even slightly improve on)
the known rates obtained for shearlet systems for the class of cartoon-like functions \cite{KW10Compactly}.
\begin{theorem}
\label{ShearCart}
Let $f \in \cC(\kappa,M,L)$ with given cartoon $\Omega \subset D = (0,1)^2$,
defined by \eqref{cartoon}.
Let $\Sigma_N$ denote the class of all piecewise affine
 functions subordinate to any partition of $D$ that consists
of at most $N$ cells which are either parallelograms or triangles.
%Further, let \wq{$\partial \Omega$ be a $C^2$ curve containing} \cs{at most}
%\wq{$\tilde M$ intersection points.}
Then
\begin{equation}
\label{eq:main_thm}
\inf_{\varphi \in \Sigma_N}\|f-\varphi\|_{L_2(D)}
\leq
{C(\kappa)L} \cdot M %\cdot \sqrt{\tilde M}
\cdot \sqrt{\log N} \cdot N^{-1},
\end{equation}
where  {$C( \kappa )$} is an absolute constant
depending only on $\kappa$. % and on the cartoon $\Omega$.
%\dw{the original version makes no sense to me. What is the mysterious dependence on $\Omega$?. All dependencies
%should be expressed by the parameters $\kappa,M,L$!}
\end{theorem}

The proof of Theorem \ref{ShearCart}, given in the following subsections, is based on constructing a specific sequence
of admissible partitions, where the
refinement decisions use full knowledge of the approximated function $f$. The above rate will then later be used as a
{\em benchmark} for the performance of a fully adaptive scheme for transport equations where such knowledge about
the unknown solution is, of course, not available. When repeatedly referring below to the construction in the proof of Theorem \ref{ShearCart}
as a ``scheme'' this is therefore not to be confused with the adaptive scheme presented later in Section \ref{sec:shearlet-implementation}.

%**************************************************************************************************
\subsubsection{Anisotropic Refinement Scheme}\label{subsect:refinement}
The main idea is to generate anisotropic meshes for which the
cells have an approximate length $2^{-j/2}$ and width $2^{-j}$, i.e.,
following the {\em parabolic scaling law ``width $\approx$ length$^2$''}.  This will later be seen to  guarantee  an improved
approximation rate for the class of cartoons defined in \eqref{cartoon}.
A very general framework for systems based on a parabolic scaling law
with an analysis of such approximation rates can be found in
\cite{GK14},\cs{ where however only unbounded domains are treated.
For bounded domains $D \subset \IR^2$ as we consider here,
$L^2(D)$-stable systems are available in \cite{GKMP15}.
}
As in the previous section, given an element $f$ from
the class \eqref{cartoon}, we shall construct now a hierarchy of
specific partitions of the domain, using at this
stage complete knowledge of the discontinuities of $f$.
%We refer in this section  to this process, underlying the
%constructive proof,  as a ``scheme''
%which is, however, not to be confused with the adaptive method used later for the solution of \eqref{varprob}
%where the searched function is, of course, not known.

We now describe the scheme which %uses the above rules but
strongly prioritizes the generation of
parallelograms.
% scales $j = 0,\dots,2J$,
For any $J\in \N$ (sufficiently large depending on the parameters $\kappa, L,M$ of the cartoon class), it successively generates anisotropic meshes.
For each scale $j = 0,\dots,2J$, it anisotropically refines only
cells intersecting the discontinuity of $f$, which will lead to about $2^J$ anisotropic
cells whose union contains all points of discontinuity of $f$.
Summing all areas of those cells gives about $2^{J} ( 2^{-2J} \times 2^{-J} ) = 2^{-2J}$, which provides
an $L_2(D)$ approximation error $\sim 2^{-J}$ for the non-smooth part of $f$.
On those cells that do not intersect any discontinuity of $f$ piecewise linear approximation
gives rise to an $L_2$ approximation error of the order of $\sim 2^{-J}$.
The number of those cells turns out to be of the order of $2^J$ times a $\log$-factor.
Combining these two error estimates, will lead  to the final bound.

The refinement scheme will  be based on an anisotropic refinement operator and an isotropic one, being applied in an
alternating fashion. The reason is to ensure, at each scale,
that the generated parallelograms follow the parabolic scaling law.
Applying only the anisotropic refinement operator would, of course,
lead to cells of approximate area $2^{-j} \times 1$,
in violation of the \cs{parabolic scaling} law.

We now proceed with studying the approximation power of our
anisotropic scheme by analyzing its performance on cartoon-like functions
$f \in \cC( \kappa,L,M,\omega)$.
%As in the previous section, given $f \in \cC( \kappa,L,M)$
%we shall construct a partition, using full knowledge about $f$,
%which will be shown to realize a certain error bound
%expressed in terms of the number of used cells.
%This is not to be confused with the adaptive scheme presented in the
%next section but merely serves to provide a benchmark.
Throughout this section, we denote by $P$, $Q$, $\cC$ and $\cP(\cC)$,
a parallelogram, a polygon (either \cs{parallelogram} or triangle),
a collection of cells and the power set of $\cC$, respectively.

We describe the refinement scheme first for a ``horizon model'', i.e., we  assume
\begin{equation}\label{eq:boundary}
\Omega = \{(x_1,x_2) \in D : x_1 \leq E(x_2)\},\quad \partial\Omega\cap D = \{(E(x_2),x_2): x_2\in (0,1)\},
\end{equation}
where %$D = (0,1)^2$ and
$E \in C^2([0,1])$ satisfies
\begin{equation}
\label{Eprop}
 \|E'\|_{L_\infty(D)} \leq 2,\quad \|E''\|_{L_\infty(D)} \leq \kappa.
 \end{equation}
Thus,  $\Gamma := \partial \Omega \cap D $
is the discontinuity curve referred to in \eqref{cartoon}.
Accordingly, in what follows we consider  for $\Omega$ from \eqref{eq:boundary} any
$f=f_1\chi_\Omega + f_2\chi_{D\setminus\Omega}$
with
$\|f_i^{(l)}\|_{L_\infty(D)}\leq M, \, l \leq 2,\,  i=1,2$,
which therefore belongs to a cartoon class $\cC(\kappa,M,L)$
with $L$ constrained by \eqref{Eprop}.

The main building blocks will be parallelograms
 obtained by parabolic scaling and shearing,
 which are defined as follows:
first, choose some $J_0 > 0$ so that
\begin{equation}\label{eq:initialchoice}
4 \cdot 5^{3/2} \cdot \kappa \leq 2^{J_0}.
\end{equation}
Moreover, we can choose $J_0$ large enough, depending on the separation parameter $\omega$, such that for every dyadic subsquare $Q\subset D$ of side-length $2^{-J_0}$
the curve $\Gamma$ intersects $\partial Q$ in at most two points.
%the bounbdary portion $\partial\Omega \cap Q$ is a graph of a function satisfying \eref{Eprop} or the analogous condition
%with the roles of $x_1,x_2$ exchanged.

The generation of an  adapted partition starts with a uniform partition
of $D= [0,1]^2$ into translates of the square $2^{-J_0}\cdot D$ of side length $2^{-J_0}$.
These squares will be further subdivided as follows:
for scale $j \in \Z$, shear $k \in \ZZ$, with corresponding operators $D_j,S_k$, defined in \eqref{eq:shear1}, and position $m \in 2^{-J_0}D^{-1}_j\ZZ^2$,
we define the parallelogram $P_{j,k,m}$ by
\begin{equation}\label{eq:shear_refine}
P_{j,k,m} := D_j^{-1}S_k(2^{-J_0}\cdot D) + m
\;.
\end{equation}
In the sequel, we refer for each fixed %$j \in \Z \backslash \{0\}$
$j \ge 0$ to the set of all such parallelograms as ${\bf P}_j$, that is,
\[
\mathcal{P}_j = \{P_{j,k,m} : k \in \Z, m \in 2^{-J_0}D^{-1}_j\Z^2\}\;.
\]
%For the initial stage, we also define
%\[
%{\mathcal P}_0 = \{P_{0,\wq{k},m} : \wq{k \in \Z,} m \in 2^{-J_0} \Z^2\}\;.
%\]
In addition to the ${\mathcal P}_j$ we consider the ``trimmed'' version
\[
\tilde{\mathcal{P}}_j = \{P_{j,k,m} \cap D \cap P : P \in {\mathcal P}_{j-1}, k \in \Z, m \in 2^{-J_0}D^{-1}_j\Z^2\},
\]
comprised of cells that are contained in $D$ as well as in parent cells from $\mathcal{P}_{j-1}$.
For $\iota \in \{0,\pm1\}$, indicating an updated shear direction,
and nonnegative {\em odd} integer $j$,
we introduce the following basic refinement operators:
for $P = P_{j-1,k,m}$,
$T_{-1} = \{(-2,0),(-1,0),(0,0)\}$,
$T_0 = \{(0,0), (1,0)\}$,
and $T_1 = \{(-1,0),(0,0),(1,0)\}$,
let
\beqn
\label{Ran}
R^{an}_{\iota}  : {\mathcal P}_{j-1} \rightarrow \cP({ \tilde{\mathcal P}}_{j}), \quad R^{an}_{\iota}(P)
= \{P_{j,2k+\iota,m+2^{-J_0-j}t} \cap D \cap P : t \in T_{\iota}\}
\eeqn
which act only on parallelograms as illustrated in Figure~\ref{fig:refinement}.
%Regarding   the basic refinement
%rules, note that $R^{an}_{\iota}$ is generated by applying the refinement operator
%$R_{\iota_Q}$ of type (R3)  when $\iota = 0$ and by
%$R_{\iota_Q}$ of type (R1) twice when $\iota = \pm 1$.
%
\begin{figure}[htb]
\begin{center}
\hspace*{\fill}
\subfigure[\label{fig:refinement-a}]{\includegraphics[width=0.4\textwidth]{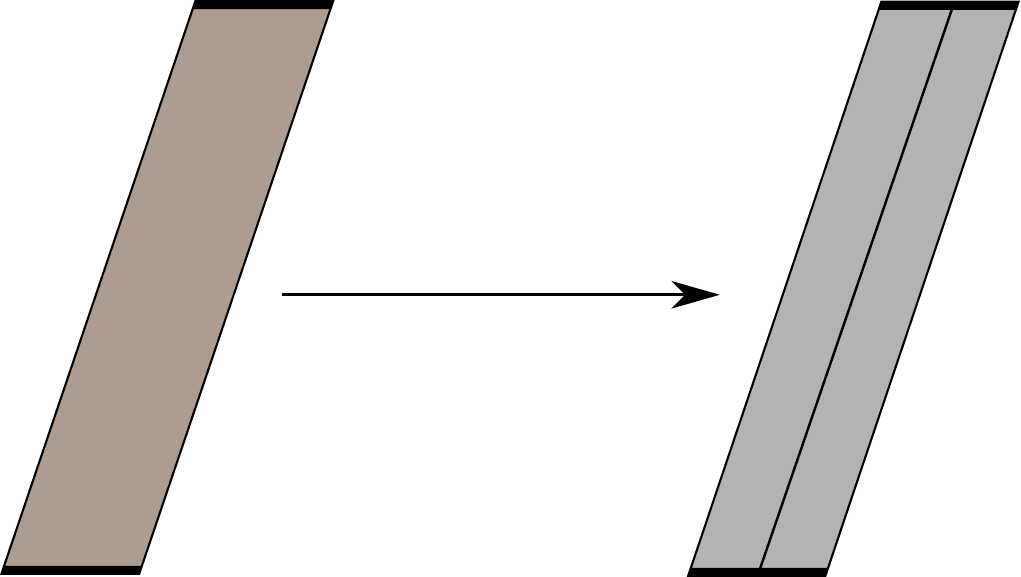}}
\hfill
\subfigure[\label{fig:refinement-b}]{\includegraphics[width=0.4\textwidth]{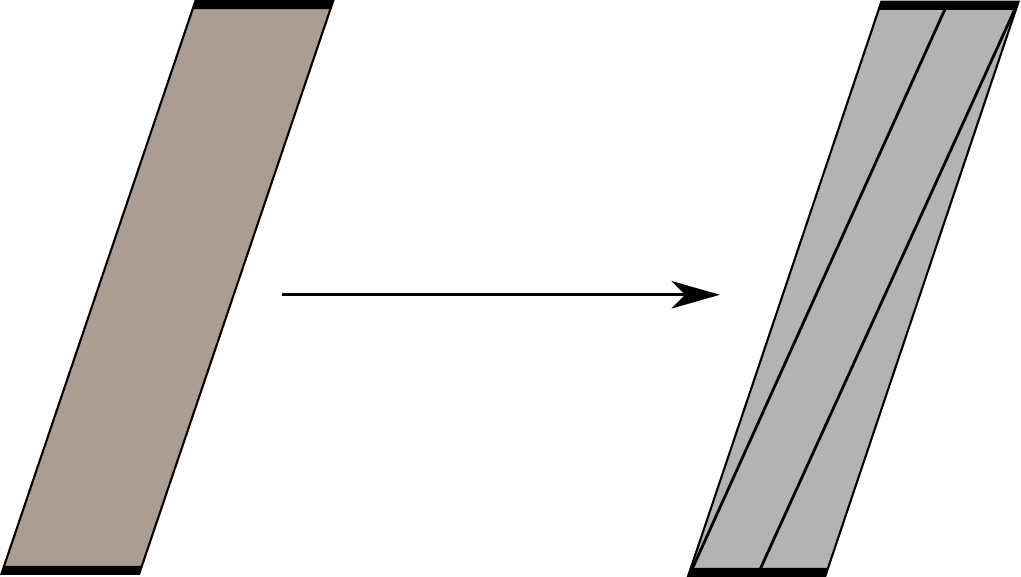}}
\hspace*{\fill}
\end{center}
\caption{Anisotropic refinements $R^{an}_{\iota}$ for $\iota = 0,1$ : \subref{fig:refinement-a}
$R^{an}_{0}(P_{j,k,m})$. \subref{fig:refinement-b} $R^{an}_{1}(P_{j,k,m})$.}
\label{fig:refinement}
\end{figure}

For $j\ge 1$ we will employ next the above basic refinement operators to construct refinement operators $R_j$, which
depend on the parity of $j$.
Specifically, for positive odd integers $j$ and any parallelogram $P=P_{j-1,k,m}\in \mathcal{P}_{j-1}$, $m=(m_1,m_2)$,
we define the orientation $\iota_P$ for splitting $P$ by
\begin{equation}\label{eq:direction01}
\iota_P= \iota_P(E): = \argmin_{\iota \in \{0,\pm 1\}} \Bigl| E'(m_2) - \frac{2k+\iota}{2^{\lceil j/2 \rceil}}\Bigr|.
\end{equation}
Thus, the orientation is determined by the slope of $E$ at the vertical position given by $m_2$. Therefore, a little care is needed
when $P$ is close to the left boundary $\{(0,x_2): x_2\in[0,1]\}$ because $\Gamma$ may enter $D$  above $m_2$, so that $E'(m_2)$ is not defined.
 In such a case,, i.e.,
$|\Gamma \cap P|>0$ but $E'(m_2)$ is not defined, we replace $m_2$ by $m_2'$ which is the second coordinate of the intersection of $\Gamma$
with $D$. We shall alwas assume this default choice in what follows without any formal distinction.

In addition, let
$R^{iso}(Q)$ be a dyadic refinement of $Q$ into four (congruent) children  of $Q$.
%Thus, $R^{iso}$ is generated by applying $R_{\iota_Q}$ of type (R3) twice
%when $Q$ is a parallelogram and by (R3) (R1) (R2) when $Q$ is a triangle.
Furthermore, define $(R^{iso})^{\ell}$ for $\ell \ge 0$ by the following recursive formula:
\beqn
\label{isol}
(R^{iso})^{\ell+1}(Q) = \{Q'' \in R^{iso}(Q') : Q' \in (R^{iso})^{\ell}(Q)\},
\eeqn
where $(R^{iso})^{0}(Q) = \{ Q \}$.

Our goal is to construct next partitions of cardinality $O(N)$ with $N=2^J$ such that the $L_2$-error is of the order at most $\sqrt{J}2^{-J}$.
To that end,  we  consider any such  $J \in \N$ which can  be thought of as significantly larger than $J_0$.
Moreover, we define for $0 \leq j \leq 2J-J_0$ and any cell $Q$  %, choose
\begin{align}
\label{lj}
\ell_j^0(Q) &:= \argmin_{\ell\in \N_0}\{\exists \mbox{  rectangle $R$ of area width $\times$ height $\sim 2^{-J/2-j/4} \times 2^{-J/2+j/4}$ such that }\nonumber\\
&\qquad (R^{iso})^{\ell^0_j(Q)}(Q) \subset R\}
 \end{align}
 as the smallest nonnegative integer such that each  cell in $(R^{iso})^{\ell_j(Q)}(Q)$ is
contained in a rectangle of area width $\times$ height $\sim 2^{-J/2-j/4} \times 2^{-J/2+j/4}$.

For $j=0$ the initial refinement $R_0(P)$ is defined as in Figure \ref{fig:init_refine}
and $R_0(P)=P$ when $\iota_P=0$.
%\dw{[WD: why can't we simply take $R_0(P)=P$, $P\in \mathcal{P}_0$?]}
For $j > 0$, we then define the operator
$R_{j}:  {\mathcal P}_{j-1}\to \cP({\tilde P}_{j})$,
by
\begin{equation}
\label{eq:refinement_op}
R_{j}(P) := \left\{
\begin{array}{rl}
R^{an}_{\iota_{P}}(P), & P\in {\mathcal P}_{j-1},  j \mbox{ odd},
\\
R^{iso}(P), & P \in {\mathcal P}_{j-1},  j \mbox{ even},
\end{array}
\right.
\end{equation}

We  describe next the refinement procedure, which
generates a sequence $\mathcal{C}_j$ of admissible partitions.
To this end, as indicated above, we first define $\cC_0$
by the partition of $D$ into dyadic squares of side length $2^{-J_0}$ as follows.
Set
\begin{equation}\label{eq:initial_cell}
\cC_0 := \{P_{0,0,m}: m_1,m_2 = 0,2^{-J_0},\dots,1-2^{-J_0}\},% = \{P\in {\bf P}_0 \cap D\},
\quad
\cC_0(\Gamma):= \{P\in\cC_0: |P\cap\Gamma|>0\}.
\end{equation}
The subsequent refinements of $\cC_0(\Gamma)$ require one further ingredient
which is the operator $\mbox{MERGE}$ that
reduces the number of triangles in favor of parallelograms.
To explain this it is useful to introduce what we call a {\em horizontal $\Gamma$-chain} of length $l\ge 1$
$$
\cH_{j,m_2}(\Gamma) := \left\{P_{j,k,(m_{1,i},m_2)}: \mbox{for some k, }
\begin{array}{l}
 |P_{j,k,(m_{1}+i-1,m_2)}\cap \Gamma|>0,
i=1,\ldots,l,\\
|P_{j,k,(m_{1}-1,m_2)}\cap \Gamma|= |P_{j,k,(m_{1}+l,m_2)}\cap \Gamma|=0
\end{array}
\right\}
$$
%of length $l\ge 1$,
comprised of a maximal collection of parallelograms at a fixed vertical
level $m_2$ which have  the same orientation and intersect $\Gamma$ substantially.
In particular, we can decompose $\cC_0(\Gamma)= \bigcup_{m_2} \cH_{0,m_2}(\Gamma)$.

Defining for any collection of parallelograms $\cC$
$$
R_j(\cC):= \{R_j(P): P\in \cC\},
$$
in particular, when $j$ is odd the collection
$R_j\big(\cH_{j,m_2}(\Gamma)\big)$ contains in general parallelograms
and triangles. We make essential use of the following observation.

\begin{lemma}
\label{lem:merge}
Assume that $j$ is even and $\cH_{j,m_2}(\Gamma)$ is a horizontal
$\Gamma$-chain of \cs{length} $l$ none of whose elements
intersects the left or right vertical part of $\partial D$.
Then, replacing all pairs of contiguous triangles in $R_{j+1}\big(\cH_{j,m_2}(\Gamma)\big)$ whose union is a parallelogram
by this parallelogram while keeping all other cells in  $R_{j+1}\big(\cH_{j,m_2}(\Gamma)\big)$ the same, yields
a collection
%\beqn
$$%\label{merge}
\mbox{\cs{\rm MERGE}}\big(\cH_{j,m_2}(\Gamma)\big) %=: \cM_{j+1}(m_2,\Gamma)
$$
with the following properties: all cells $Q$ in
$\mbox{\cs{\rm MERGE}}\big(\cH_{j,m_2}(\Gamma)\big)$ %$\cM_{j+1}(m_2,\Gamma)$
that substantially intersect $\Gamma$, i.e., $|Q\cap\Gamma|>0$,
form a horizontal $\Gamma$-chain of length not larger than $2l-1$.
\end{lemma}
\cs{\begin{proof}
The validity of the above claim is easily confirmed by noting that by
\eref{eq:direction01} $\iota_P=\iota_{P'}$ for any two $P,P'\in \cH_{j,m_2}$
and using the maximality of a horizontal $\Gamma$-chain.
\end{proof}}

\cs{For arbitrary, finite union $\cH$} of horizontal $\Gamma$-chains we denote
$$
\mbox{MERGE}\big(\cH\big) := \bigcup\,\big\{\mbox{MERGE}\big(\cH_{j,m_2}(\Gamma)\big): \cH_{j,m_2}(\Gamma) \in \cH\big\}.
$$

The   statement in Lemma \ref{lem:merge} regarding the fact that $\Gamma$ is only intersected by
parallelograms after merging,  is not necessarily true when the first or last parallelogram $P$ in a horizontal $\Gamma$-chain
intersects the left or right vertical boundary part of $D$, respectively. In fact, when $\iota_P \neq 0$ the left and rightmost triangles cannot be
merged but may still be substantially intersected by $\Gamma$. To account for this fact,
note first that due to the constraints on $\Gamma$ there exist at most a finite number of intersections of $\Gamma$ with
vertical boundary portions of $\partial D$.
We may assume without loss of generality that
$J_0$ is large enough (depending only on the parameters $\kappa,L,\omega$)
such that a cell in $\cC_0$ contains at most one intersection of $\Gamma$ with the
vertical boundary
$$
(\partial D)^v:= \{(\bar x,x_2): \bar x \in \{0,1\},\, x_2\in [0,1]\}.
$$
 Let us denote
\beqn
\label{CB}
\cC_0^B := \{P\in \cC_0: P\cap \Gamma \cap (\partial D)^v \neq \emptyset\}, \quad \cC_0'(\Gamma):= \cC_0(\Gamma)\setminus
\cC_0^B.
\eeqn
As said before $\#(\cC_0^B)$ is finite depending only on $\kappa,L,\omega$.

We are now prepared to describe the refinement process first
only for the cells in
$\cC_0'(\Gamma)$, $\cC_0':= \cC_0\setminus \cC_0^B$.
As mentioned before $\cC'_0(\Gamma)$ can be decomposed
into horizontal $\Gamma$-chains so that the collection
\begin{equation}
\label{eq:init01}
\cC^{an}_{0,0} = \mbox{MERGE}(\{ Q \in R_0(P): P \in \cC'_0(\Gamma) \}) = \mbox{MERGE}\big(R_1(\cC'_0(\Gamma))\big)
\end{equation}
%\dw{[WD: I think the $R_1(P)$ above comes from my understanding of the initialization of the refinement process, see my eralier
%comment. Feel free to change if you think my simpler initialization does not work.]}
is well defined and has, by Lemma \ref{lem:merge}, the following property: all cells $Q$
in $\cC^{an}_{0,0}$ that substantially intersect $\Gamma$ are parallelograms with the same shear orientation for each vertical level.

 Given $\cC^{an}_{0,0}$, we can initialize the collection %of isotropically refined cells
\begin{equation}
\label{eq:init02}
\cC^{iso}_{0,0} := \{Q' \in (R^{iso})^{\ell^0_0(Q)}(Q) : Q \in \cC^{an}_{0,0} \cup \cC'_0, |Q \cap \Gamma| = 0\}
\end{equation}
of those cells that result from an isotropic refinement of the cells $Q$ with $|Q \cap \Gamma| = 0$,
see Figure \ref{fig:adaptive} and the definition \eref{lj} of the refinement depth $\ell^0_0(Q)$.

%\vspace*{3cm}

On account of Lemma \ref{lem:merge}, we can now recursively define $(\cC_{j}^{an})_{j,0}$ and $(\cC_{j}^{iso})_{j,0}$, $j>0$ by
\begin{equation}\label{eq:aniso_cell01}
{\cC}^{an}_{j,0} := \mbox{MERGE}(\{Q' \in {{R}_{j}(Q)} : Q \in \cC^{an}_{j-1,0}, |Q \cap \Gamma| \neq 0\}),\quad 0< j\le 2J-J_0,
\end{equation}
and
\begin{equation}\label{eq:iso_cell01}
\cC^{iso}_{j,0} := \{Q' \in (R^{iso})^{\ell^0_j(Q)}(Q) : Q \in \cC^{an}_{j,0}, |Q \cap \Gamma| = 0\},
\end{equation}
see \eref{lj} for the definition of $\ell^0_j(Q)$ as the smallest integer for which
  each  cell in $(R^{iso})^{\ell^0_j(Q)}(Q)$ is
contained in a rectangle of area width $\times$ height $\sim 2^{-J/2-j/4} \times 2^{-J/2+j/4}$.
We also set ${\cC}^{an}_{j,0}:= \emptyset$ for $j> 2J-J_0$.
%Using the recursion \eqref{eq:aniso_cell01} described above,
%we iteratively generate $\cC^{an}_j$ for
%$j = 0,1,\dots,2J$, and set $\cC^{an}_j := \emptyset$ for
%$j > 2J$, where $J$ is some positive integer we have chosen before.
%
\begin{figure}[htb]
\begin{center}
\hspace*{\fill}
\subfigure[\label{fig:refine01-a}]{\includegraphics[width=0.2\textwidth]{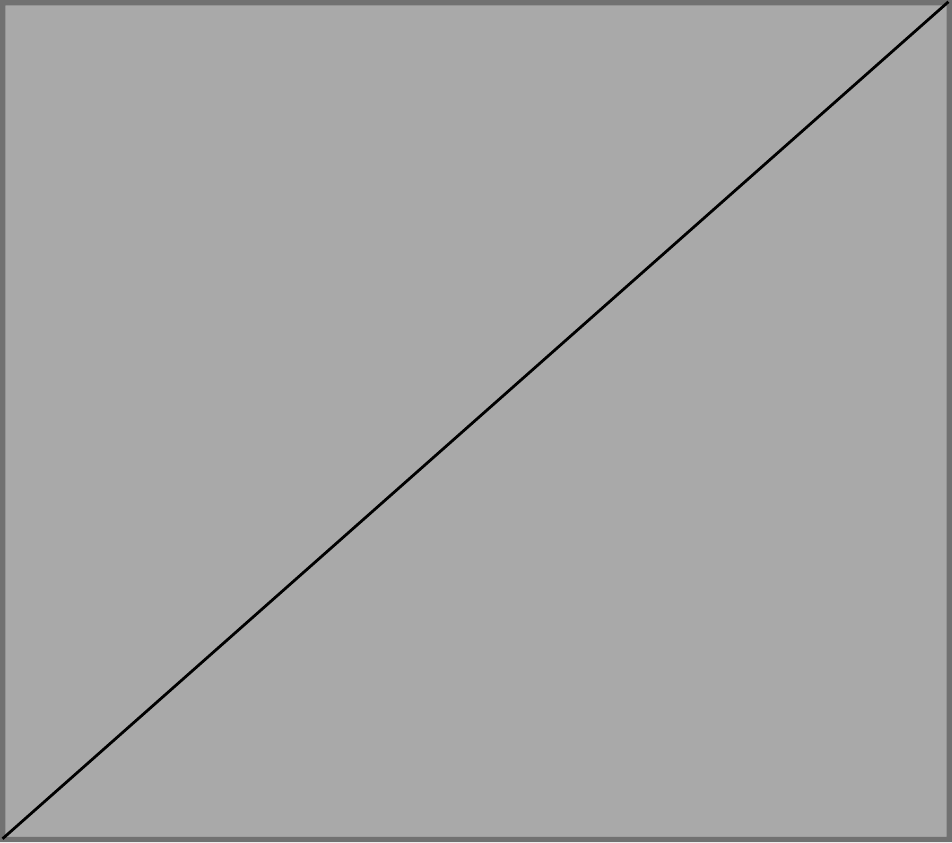}}
\hfill
\subfigure[\label{fig:refine02-b}]{\includegraphics[width=0.2\textwidth]{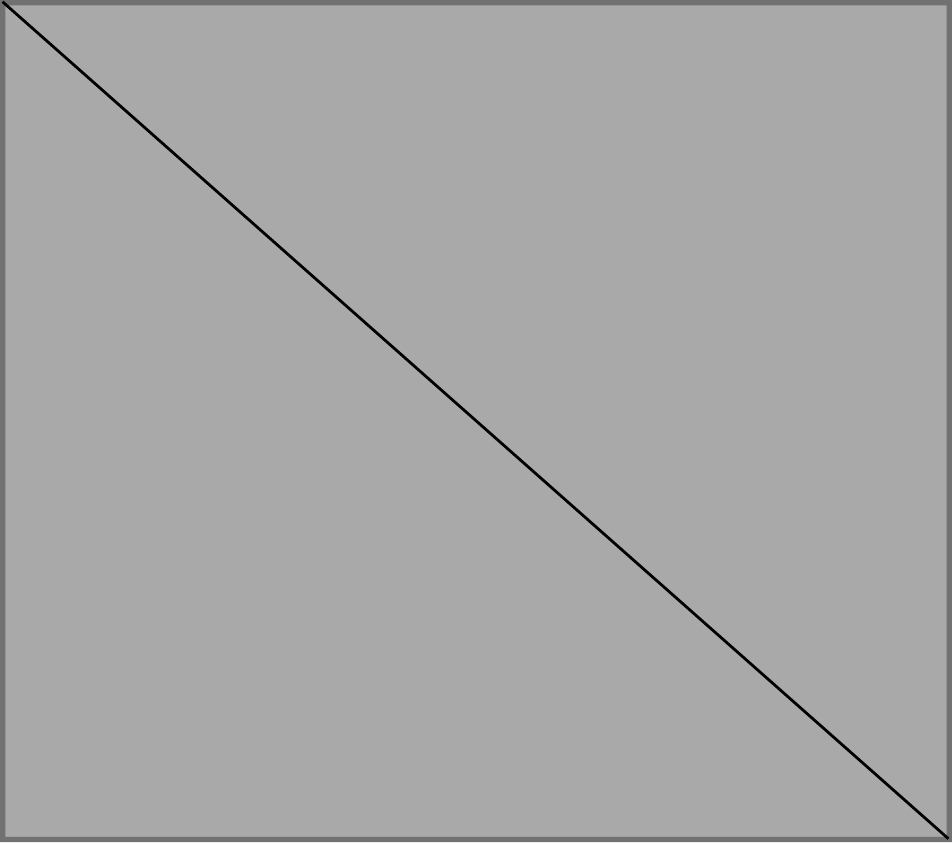}}
\hspace*{\fill}
\end{center}
\caption{Initial refinement $R_{0}(P)$ : \subref{fig:refine01-a}
$R_0(P)$ for $P \in {\mathcal P}_0$ when $\iota_P = 1$.
\subref{fig:refine02-b} $R_0(P)$ for $P \in {\mathcal P}_0$ when $\iota_P = -1$.
%\dw{[WD: here again Wang-Q suggests a change because of the initialization he proposes. I would prefer using the simpler one, but
%if I made a mistake, cahnge it.]}
}
\label{fig:init_refine}
\end{figure}
Note that the collections $\cC^{an}_{j,0}, \cC^{iso}_{j,0}$ give rise to partitions.
Since the initializing collections
contain only pairwise disjoint cells it is clear that the collections
$\cC_{j,0}^{an}$ and $\cC_{j,0}^{iso}$ are also only
comprised of cells with pairwise disjoint interiors.
More precisely, for any $m\in \N$, the collection
\beqn
\label{Cm}
\cC_{m,0}(f):= \cC^{an}_{m,0} \cup \bigcup_{j < m} \cC^{iso}_{j,0}
\eeqn
forms a partition of
$$
D_0 : =[0,1]^2\setminus \bigcup\,\{ P\in \cC_0^B\}
$$
comprised of parallelograms and triangles.

\begin{figure}[htb]
\begin{center}
\includegraphics[height=2in]{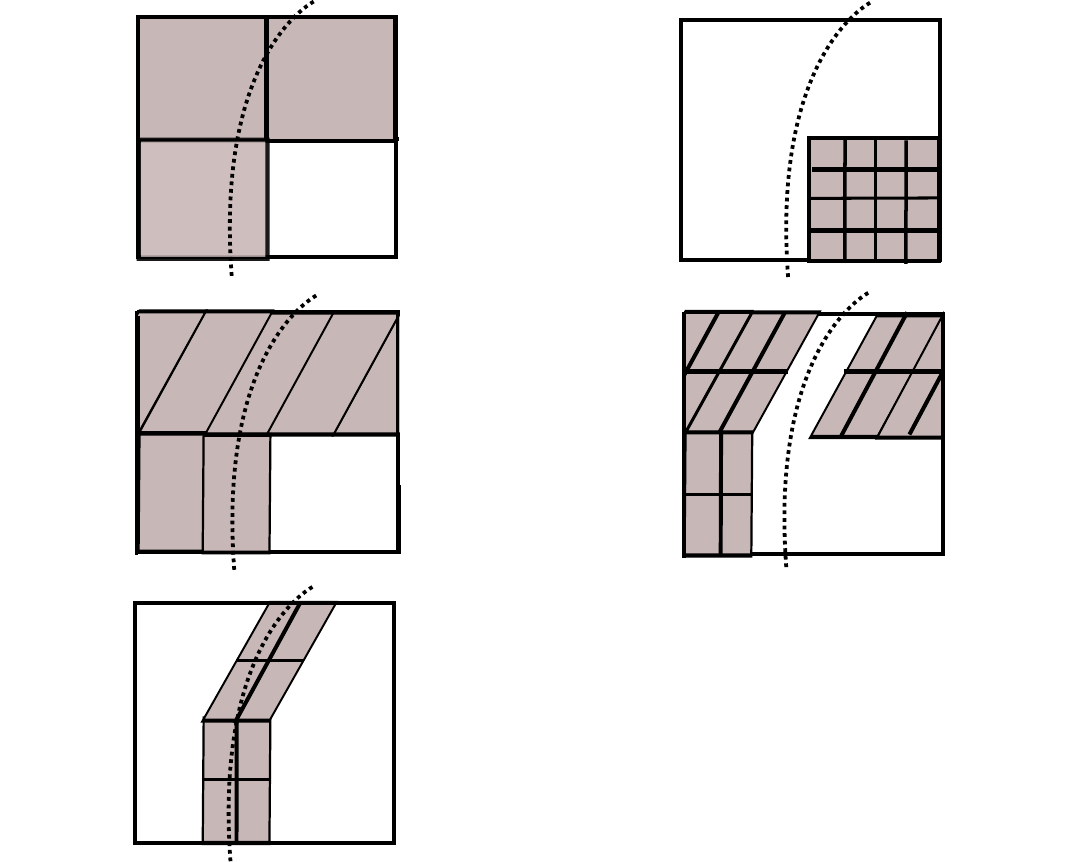}
\put(-110,100){$\cC^{an}_0$}
\put(-110,60){$\cC^{an}_1$}
\put(-110,20){$\cC^{an}_2$}
\put(-20,100){$\cC^{iso}_0$}
\put(-20,60){$\cC^{iso}_1$}
\end{center}
\caption{Recursive construction of $\cC^{an}_j$ and $\cC^{iso}_j$.}
\label{fig:adaptive}
\end{figure}

%Let us denote now by  $\cB^{an}_{j,0}$, $\cB^{iso}_{j,0}$, $j\geq 0$, piecewise polynomial bases
%for the partitions $\cC^{an}_{j,0}$, $\cC^{iso}_{j,0}$, respectively, so that
%
%\begin{equation}\label{eq:basis}
%\mathbb{P}_1(\cC^{an}_{j,0}) = \text{{\bf span}}(\cB^{an}_{j,0}) \quad \text{and} \quad  \mathbb{P}_1(\cC^{iso}_{j,0}) = \text{{\bf span}}(\cB^{iso}_{j,0})
%\quad \text{for} \,\, j \ge 0.
%\end{equation}
%
%Obviously, we have
%\beqn
%\label{Bcard} {
%(d+1)^{-1}\#(\cB^{an/iso}_j) \leq \#(\cC_j^{an/iso}) \leq (d+1)\#(\cB^{an/iso}_j).
%,\quad (d+1)^{-1}\#(\cB^{iso}_j) \leq \#(Q\cC_j^{iso}) \leq (d+1)\#(\cB^{iso}_j).
%3\#(\cC_{j,0}^{an/iso}) = \#(\cB^{an/iso}_{j,0}).}
%\eeqn
%

%We will employ linear combinations of elements from the  redundant system
%This construction yields a function system associated with adaptively refined cells by setting
%\begin{equation}\label{eq:total_basis}
%\cB_0(f) =\cB_0(f;J_0,J)
%:=
%\bigcup_{j \ge 0} \Bigl(\cB^{iso}_{j,0} \cup \cB^{an}_{j,0}\Bigr),
%\end{equation}
%
%to approximate $f$. Specifically, for any $m\in \N$, the collection
%
%\beqn
%\label{Bm}
%\cB_{m,0}(f) := \{\chi_Q g: g\in \PP_1,\, Q \in \cC_{m,0}(f)\}
%\eeqn
%
%forms a basis for $\cS_m({D_0}):=\PP_1(\cC_{m,0}(f))$ because $\cC_{m,0}(f)$ is a partition of $D\setminus \bigcup\,\{P\in \cC_0^B\}$.

It remains to construct similar partitions for the remaining squares $P_0\in \cC_0^B$
\cs{with nonempty intersection with $(\partial D)^v$.}
%having a nontrivial contact with $(\partial D)^v$.
There is a uniform upper bound $B$ of such cells over the given cartoon class,
depending only on the parameters $\kappa,L,\omega$.
Thus it suffices to discuss
one such square $P_0$ of side length $2^{-J_0}$ in $\cC_0^B$ whose vertical boundary has exactly one point $g$ of intersection with $\Gamma$.
Hence, there exists exactly one square $P_1\in R^{iso}(P_0)$ that contains $g$.
Let
$$
\cC'_{0,1}:= R^{iso}(P_0) \setminus \{P_1\},\quad \cC'_{0,1}(\Gamma):= \{P'\in \cC'_{0,1}: |P'\cap\Gamma|>0\}.
$$
More generally, for $r\in\N_0$, let $P_r(P_0)\subset P_0$ be the unique square of side length $2^{-J_0-r}$ containing the intersection point
$g$ and define %for $r\le J-J_0-1$
\beqn
\label{Cj}
\cC'_{0,r+1}:= R^{iso}(P_r) \setminus \{P_{r+1}\},\quad \cC'_{0,r+1}(\Gamma):= \{P'\in \cC'_{0,r+1}: |P'\cap\Gamma|>0\}.
\eeqn
For each collection $\cC'_{0,r}$, $0\le r\le J-J_0$, we can now apply the above construction to arrive at collections
$\cC^{an}_{j,r}, \cC^{iso}_{j,r}$, $j>0$, providing partitions
\beqn
\label{Cmr}
\cC'_{m,r}(f,P_0):=% \{P_{r+1}(P_0)\}\cup
\cC^{an}_{m,r} \cup \bigcup_{j < m} \cC^{iso}_{j,r}
\eeqn
of  the $L$-shaped domains
$$
D_r=D_r(P_0) := P_r\setminus P_{r+1} \subset P_0. %, \quad r\le J-J_0-1.
$$
Note that it suffices to stop this zoom as soon as $P_{r+1}$ has side length $2^{-J}$ because when approximating
$f$ on $P_{r+1}$ by zero, the squared $L_2$-error on this square  is bounded by $M^2$ times  its area which
is $2^{-2J}$ and there are at most $\#(\cC_0^B)$ such squares each contributing a squared error of the order $2^{-2J}$
which is within our target accuracy, i.e., $0\le r \le J-J_0-1$.
For each $r$ in this range we need to consider only anisotropic refinements of level up to $2J$ which yields the following ranges
% consider only the range
\beqn
\label{ranger}
0\le j \le 2J- J_r, \quad 0\le r \le J-J_0-1,  \quad \mbox{where}\quad J_r := J_0 + r.
\eeqn
 %, rooted in $P_0\in \cC_0^B$.
%along with the associated spaces of piecewise polynomials.
%Thus, defining
%\beqn
%\label{BmP0}
%\cS_{m,r}(P_0) := \PP_1\big(\cC'_{m,r}(f)\big),
%\eeqn
Thus, on each $P_0\in \cC_0^B$ we are led to consider the partition
 \beqn
\label{TP0}
\cT_J(P_0) := \{P_{J-J_0}\} \cup \bigcup_{r=0}^{J-J_0-1}\cC'_{J+1,r}(f,P_0)
\eeqn
of $P_0$ so that
\beqn
\label{TJ}
\cT_J(f) := \cC_{2J-J_0,0}(f) \bigcup_{P_0\in \cC_0^B} \cT_J(P_0)
\eeqn
is a partition of $D$ adapted to $f$, see \eref{Cm}.

%***************************************************************************************************************************
\subsubsection{{Proof of Theorem \ref{ShearCart}}}
We are now ready to provide the proof of Theorem \ref{ShearCart} by analyzing the approximation error provided by $\cT_J(f)$
first under the above restricted assumptions on $f\in \cC(\kappa,L,M,\omega)$. Once such a result is at hand, we point out how
to extend it to the full cartoon class.

We adhere to the notation of the preceding section and denote
for any partition $\cT$ by $\PP_p(\cT)$ the space of piecewise
polynomials of degree at most $p\in \N$ subordinate to the partition $\cT$.

% which
%shows that the above shearlet-type partitions
%realize approximation rates that are essentially optimal for directional representation systems in the sense of \cite{D2001}) as
%well as for anisotropic triangulations \cite{CoDyHeMi12,CM1,CM2,M,ChSuXu07}.
%\gw{(Reference \cs{\cite{ChSuXu07} Xu et.al. (Determinant of Hessian)?)}
Recalling \css{the definition \eqref{Sigman} of $\Sigma_N$,}
%that
%$$
%\Sigma_N := \bigcup\{ \PP_1(\cT):  \cT \, \mbox{a partition of $D$ into parallelograms
%and triangles with}\, \#(\cT)\leq N\},
%$$
we determine for given cartoon
$f$ in $\cC(\kappa,L,M,\omega)$ the cardinality
$N= N(J):= \sharp(\cT_J(f))$, for $\cT_J(f)$
defined by \eref{TJ}, and then estimate
$\inf_{g\in \PP_1(\cT_J(f))}\|f-g\|_{L_2(D)}$.

\cs{By} construction, we have $\PP_1(\cT_m(f)) \subset \Sigma_{\#(\cT_m(f))}$.
%\noindent {\bf Proof of Theorem \ref{ShearCart}:}
We assume first that the discontinuity curve
$\Gamma = \partial \Omega \cap D$ is given as in \eqref{eq:boundary}. Moreover, we analyze first the approximation from
the spaces $\PP_1(\cC_{m,0}(f))$ on the domain $D_0$ (excluding the squares with vertical boundary intersections of $\Gamma$),
see \eref{Cm}.
%The asserted bound will be established by approximations based on suitable subsets of
%$\cB(f;J_0,J)$ from \eqref{eq:total_basis}.
To simplify notation we abbreviate in this part of the proof
$$\cC^{an}_j:= \cC^{an}_{j,0},\quad    \cC^{iso}_j:= \cC^{iso}_{j,0},\quad  \ell_j(Q):= \ell_j^0(Q).
$$

We begin with recording the following properties of the collections  $\cC^{an}_j$ that will be used repeatedly:
%These are the following:
%
\begin{itemize}
\item[{\bf (C1)}] $\Bigl|\Gamma \cap \Bigl(\bigcup_{Q \in \mathcal{C}^{an}_j}Q\Bigl)^c\Bigr| = 0$.
\item[{\bf (C2)}] For each $Q \in \mathcal{C}^{an}_j$ with $|Q\cap \Gamma|>0$, we have $Q = P \subset D$
for some $P = P_{j,k,m}$ and
\begin{equation}\label{eq:first}
\sup_{ (x_1,x_2)\in Q } \Bigl| E'(x_2) - \frac{k}{2^{ \lceil j/2 \rceil }}\Bigr| \leq \frac{1}{2^{ \lceil j/2 \rceil }}.
\end{equation}
\item[{\bf (C3)}] In {\bf (C2)}, $P$ is a parallelogram of size width $\times$ height $\lesssim 2^{-j} \times 2^{-j/2}$.
\end{itemize}

We now verify first  the validity of these properties. {\bf (C1)} and {\bf (C3)} follow directly from our construction. We verify
next {\bf (C2)} by induction on $j$. For $j=0$ note first that for $m_2\le s \le m_2+ 2^{-J_0}$
\beqn
\label{simple}
|E'(s)-E'(m_2)|\le \int_{m_2}^s |E''(t)|dt \le (4\cdot 5^{3/2})^{-1},
\eeqn
where we have used \eref{eq:initialchoice}.
\cs{Let us} consider first $E'(m_2)\ge 0$.
Then we have   $\iota_P=0$, if $E'(m_2)< 1/2$,
and $\iota_P= 1$ when $E'(m_2)\ge 1/2$. When   $E'(s)$ remains positive we have
$|E'(s)-1|\le 1$, by \eref{Eprop}, which confirms the claim in this case. When $E'(s)<0$ for some $s$ \eref{simple}
says that $|E'(s)-E'(m_2)|= |E'(m_2)| + |E'(s)| \le  (4\cdot 5^{3/2})^{-1}$ which gives $\iota_P=0$ and again confirms the
claim for $j=0$. An analogous argument can be used when $E'(m_2)<0$ to conclude that \eref{eq:first} holds for $j=0$.

Let us assume now that $j \ge 0$ is even and
that \eqref{eq:first} is satisfied for this $j$. We then prove that
\begin{equation}\label{eq:next}
\sup_{(x_1,x_2)\in Q'} \Bigl| E'(x_2) - \frac{2k+\hat\iota}{2^{ \lceil (j+1)/2 \rceil }}\Bigr| \leq \frac{1}{2^{ \lceil (j+1)/2 \rceil }}
\end{equation}
for any $Q' \in R^{an}_{\hat \iota}(P_{j,k,m})$ with $|Q' \cap \Gamma| \neq 0$, where
\begin{equation}\label{eq:direction}
\hat \iota = \argmin_{\iota \in \{0,\pm 1\}} \Bigl| E'(m_2) - \frac{2k+\iota}{2^{j/2+1}}\Bigr|.
\end{equation}
%Let $Q' \in R^{an}_{\hat \iota}(P_{j,k,m})$ be arbitrarily chosen.
By Lemma \ref{lem:merge} and because we are working only on $D_0$, we can assume
without loss of generality that $Q'
= P_{j+1,2k+\hat{\iota},\tilde{m}} \in {\bf P}_{j+1}$ with $\tilde{m} = (\tilde m_1, m_2)$, for some $\tilde{m}_1 \in (0,1)$,
 is a parallelogram contained in $D$.
%Then observe that we have $Q' \subset P' \cap D$ and
%$P' = P_{j+1,2k+\hat{\iota},\tilde{m}} \in {\bf P}_{j+1}$ with $\tilde{m} = (\tilde m_1, m_2)$
%for some $\tilde{m}_1 \in (0,1)$.
If $(x_1,x_2) \in P_{j+1,2k+\hat{\iota},\tilde{m}}$,
then
\begin{equation}\label{eq:distance}
|x_2-m_2| \leq 2^{-J_0-j/2}
\end{equation}
Thus
\[
\Bigl| E'(x_2) - \frac{2k+\hat{\iota}}{2^{ j/2 + 1}}\Bigr| \leq
\Bigl| E'(m_2) - \frac{2k+\hat{\iota}}{2^{ j/2 + 1}}\Bigr| + |E'(m_2) - E'(x_2)|. %
\]
This now implies that
\[
\sup_{Q' \in (x_1,x_2)} \Bigl| E'(x_2) - \frac{2k+\hat\iota}{2^{ \lceil (j+1)/2 \rceil }}\Bigr|
\leq \Bigl| E'(m_2) - \frac{2k+\hat{\iota}}{2^{ j/2 + 1}}\Bigr| + |E'(m_2) - E'(x_2)| =  \text{(I)} + \text{(II)}.
\]
By \eqref{eq:first} and \eqref{eq:direction},
\[
\text{(I)} \leq (\frac{1}{2})2^{-j/2-1}.
\]
Again, by \eqref{eq:distance},
\[
\text{(II)} \leq \|E''\|_{\infty}2^{-J_0-j/2} \leq (\frac{1}{2})2^{-j/2-1},
\]
%with sufficiently large $J_0$ -- this is the role of the initial scale $J_0$ in \eqref{eq:initialchoice}.
where in the last step we have used the condition \eqref{eq:initialchoice} on the initial scale $J_0$.
%\dw{[WD: As far as I can see we don't use the factor $5^{3/2}$ from \eqref{eq:initialchoice} here nor in the case $j=0$, so perhaps this can be relaxed.
%I'll keep it for now in case another constraint comes up.]}
This
proves \eqref{eq:next}.

Finally, when $j$ is odd, we note  that \eqref{eq:first} obviously implies \eqref{eq:next}
for any $Q' \in R^{iso}(P_{j,k,m})$,
since $\lceil j/2 \rceil = \lceil (j+1)/2 \rceil $.
Thus, {\bf (C2)} follows   by induction.

We proceed bounding next $\#(\cC_{2J{-J_0}}(f))$, see \eref{Cm}, and
noting that ${\rm dim}\,(\PP_1(\cC_{2J{-J_0}}(f)))= 3 \#(\cC_{2J{-J_0}}(f))$.
To this end, {\bf (C2)} and {\bf (C3)} together with the fact that each $P \in \cC^{an}_j$ is for $j>0$  a refinement of
some $P' \in \cC^{an}_{j-1}$ intersecting $\Gamma$,   the length of $\Gamma\cap Q$ scales for $Q\in \cC_j^{an}$ asymptotically like ${\rm diam}\,(Q)\sim 2^{-J_0- j/2}$ so that
\beqn
\label{cardCanj}
\sharp(\mathcal{C}^{an}_j)
\lesssim L 2^{J_0 + j/2},\quad j\ge 0,
\eeqn
with a constant depending only on $\kappa$.
%imply
%$\sharp(\mathcal{C}^{an}_j) %\sim \sharp(\cB^{an}_{j})
%\lesssim L 2^{j/2}$,  $j \ge 0$,
%where the constant depends only on  $\kappa$.
 % for $f\in \cC(\cs{\Omega},\kappa, L,M)$ on \cs{$\Omega$, $\kappa$ and on $L$.}
 Hence, %Since $\sharp(\cB^{an}_{j}) \le 3 \sharp(\mathcal{C}^{an}_j)$
we obtain
\begin{equation}\label{eq:num1}
\sharp\Bigl(\bigcup_{j = 0}^{2J{-J_0}}\cC^{an}_{j}\Bigr) \lesssim {L2^{J_0}} \sum_{j = 0}^{2J{-J_0}} 2^{j/2} \lesssim L 2^{J{+J_0/2}},
\end{equation}
with a constant depending only on $\kappa$.
Note that our count is generous because of the overlap of cells on different levels $j$.
%since the union of bases is a
%linearly dependent set and we only need to count the basis
%functions on leaf-cells.
 By majorization with a summable geometric series,
such redundancies are always controlled by a constant uniform factor.

We now estimate $\sharp(\bigcup_{j = 0}^{2J}\cC^{iso}_{j})$
%First, note that $\sharp(\cB^{iso}_{j})
%=3 \sharp(\cC^{iso}_j)$ for each $j \ge 0.$
 observing first that for each $Q \in \cC^{an}_j$,
\[
\sharp((R^{iso})^{\ell_j(Q)}(Q)) = 2^{2\ell_j(Q)},
\]
 where $\ell_j(Q)=\ell_j^0(Q)$ (as defined in \eref{lj}) is chosen so that $2^{-\ell_j(Q)}Q$
is contained in a rectangle of size width $\times$ height
$\lesssim 2^{-\frac{J}{2}-\frac{j}{4}} \times 2^{-\frac{J}{2}+\frac{j}{4}}$.
 This is accomplished by taking %one may assume that
 \beqn
 \label{lj0}
\ell_j(Q)=\ell_j^0(Q) = \max\,\Big\{\Big\lceil \frac J2- \frac{3j}{4} - J_0\Big\rceil, 0\Big\},
\eeqn
i.e., $\ell_j(Q)=0$ when $j\ge 2(J-2J_0)/3$.
%if $0 \leq j \leq \lceil \frac{2}{3}J \rceil$
%and $\ell_j(Q) = 0$ otherwise.
Therefore, by \eqref{eq:iso_cell01}, we obtain, upon recalling that $\sharp(\cC^{an}_j) \lesssim L2^{{J_0+} j/2}$
\begin{eqnarray}
\label{eq:num2}
\sharp(\bigcup_{j = 0}^{2J{-J_0}}\cC^{iso}_{j}) &{\le}&  \sum_{j = 0}^{\lceil \frac{2J}{3}{- \frac{4J_0}{3}}\rceil}2^{2(J/2-3/4j  {-J_0})}\sharp(\cC^{an}_j)
+ \sum_{j > \lceil \frac{2J}{3}{- \frac{4J_0}{3}}\rceil} \sharp(\cC^{an}_j)\nonumber\\
& \lesssim &
%L \sum_{j = 0}^{\lceil \frac{2J}{3}\rceil}2^{2(J/2-3/4j)+j/2}
%+ L \sum_{j = \lceil \frac{2J}{3}\rceil+1}^{2J} 2^{j/2} \nonumber\\
%&=&
L \Big(\sum_{j = 0}^{\lceil \frac{2J}{3}{- \frac{4J_0}{3}}\rceil} %2^{J-j/2}
{2^{J-j-J_0}}
+ {2^{J_0}}  \sum_{j = \lceil \frac{2J}{3}{- \frac{4J_0}{3}}\rceil+1}^{2J{-J_0}} 2^{j/2}\Big)\lesssim L 2^{J{+J_0/2}},
\end{eqnarray}
 where we have used \eqref{eq:num1} in the last step, so that the constant in the above
estimate depends only on $\kappa$.

In summary, by \eqref{eq:num1} and \eqref{eq:num2},  we conclude that
\begin{equation}
\label{Nest}
N = {\rm dim}\,\PP_1(\cC_{2J{-J_0},0}(f))
\lesssim
{ L {2^{J_0/2}} 2^J =: {2^{J_0/2}}L \, N_0, \quad N_0\sim 2^J},
\end{equation}
 see \eqref{Cm}, with a constant depending only on $\kappa$ and $\omega$.
%The system $\cB_{2J}(f)$ contains a basis for $\mathbb{P}_1(\cC_{2J}(f))$,
%where $\cC_{2J}(f)$ is a partition of $D$ defined in \eqref{Cm}.

After bounding the number of degrees of freedom associated with the constructed partition
we now turn to estimating the approximation error.
To this end, recall that,  by \eqref{Cm}, $\cC_{2J{-J_0}}(f) = \cC^{an}_{2J{-J_0}} \cup \cC^{iso}$, where $\cC^{iso} = \cup_{j = 0}^{2J-1{-J_0}} \cC^{iso}_j$
and notice further that on account of
 {\bf (C1)}, \eqref{eq:aniso_cell01}, \eqref{eq:iso_cell01},
one has
 $|Q \cap \Gamma| = 0$ for any $Q \in \cC^{iso}$.
Thus, for $Q \in \cC^{iso}$, denoting by $P_{Q}(f)$ the best $L_\infty$-approximation
to $f$ from $\PP_1$ on $Q$,
a classical Whitney estimate states
$\|f-P_Q(f)\|_{L_\infty(Q)}\leq C ({\rm diam}\,(Q))^2\|f\|_{{W^{2,\infty}(Q)}}$
 where  $C>0$ is independent of the aspect ratio of $Q$.
Thus, since $f \in  \cC(\kappa,L,M,\omega)$,
we conclude that
\begin{equation}
\label{eq:error}
\|f-P_{Q}(f)\|^2_{L_2(Q)} \leq |Q|\|f - P_Q(f)\|_{L_\infty(Q)}^2   {\le  C ({\rm diam}\,(Q))^4|Q| M^2} \;.
\end{equation}
If $Q \in \cC^{iso}_{j}$ with $j > \lfloor \frac{2J}{3}{-\frac{4J_0}{3}} \rfloor$,
then $|Q| \lesssim 2^{-j {-J_0}} \times 2^{-j/2 {-J_0}}$
and
${\rm diam}\,(Q)\lesssim 2^{-j/2 {-J_0}}$, so that,
by \eqref{eq:error},
\beqn
\label{PQ1}
\|f-P_{Q}(f)\|^2_{L_2(Q)} \lesssim 2^{-\frac{7}{2}j}{2^{-6J_0}}M^2\;.
\eeqn
On the other hand, if
$Q \in \cC^{iso}_{j}$ with $j \leq \lfloor \frac{2J}{3} {-\frac{4J_0}{3}}  \rfloor$,
then
$|Q| \lesssim 2^{-j-\ell_j(Q) {-J_0}} \times 2^{-j/2-\ell_j(Q){-J_0}}$, {${\rm diam}\,(Q)\sim 2^{\frac j4-\frac J2}$}
with
$\ell_j(Q)  {\sim}  \lceil J/2-3j/4 {-J_0} \rceil$
and the same reasoning yields
\beqn
\label{PQ2}
\|f-P_{Q}(f)\|^2_{L_2(Q)} \lesssim %2^{-\frac{7}{2}j-6(J/2-3/4j)}M^2 =
2^{j-3J}M^2
\;.
\eeqn
Therefore, the approximation error away from $\Gamma$ can be bounded by
\begin{eqnarray}\label{eq:estimate01}
\sum_{Q \in \cC^{iso}}\|f-P_{Q}(f)\|^2_{L_2(Q)} &\lesssim&
\sum_{j = \lceil \frac{2J}{3}{- \frac{4J_0}{3}} \rceil+1}^{2J-1 {-J_0}} \sum_{Q \in {\cC}^{iso}_j}\|f-P_{Q}(f)\|^2_{L_2(Q)}\nonumber
\\ &&\quad
+
\sum_{j = 0}^{\lceil \frac{2J}{3}{- \frac{4J_0}{3}}\rceil}\sum_{Q \in \mathcal{C}^{iso}_{j}}\|f-P_{Q}(f)\|^2_{L_2(Q)}\nonumber
\\
&\lesssim& L
M^2\Big({2^{-7J_0}2^J} \sum_{j = \lceil \frac{2J}{3} {- \frac{4J_0}{3}}  \rceil}^{2J-1 {-J_0}}
%(2^{j/2})(2^{-\frac{7}{2}j})
{2^{-\frac{9j}{2}}}+\sum_{j = 0}^{\lceil \frac{2J}{3}{- \frac{4J_0}{3}}\rceil}
{2^{J-j-J_0} 2^{j-3J}}
%(2^{j/2+2(J/2-3/4j)})(2^{j-3J})
%2^{-\frac{7}{2}j-6(J/2-3/4j)})
\Big)
\nonumber
\\
&\lesssim& {2^{-J_0}} L M^2 2^{-2J} {\Big(1+ \left\lceil \frac{2J}{3}{- \frac{4J_0}{3}}\right\rceil_+\Big)}\nonumber\\
&  \sim &  L \,M^2 {2^{-J_0}} (\log N_0) N_0^{-2},
\end{eqnarray}
where the constant depends only on $\kappa$ and
$\left\lceil \frac{2J}{3}{- \frac{4J_0}{3}}\right\rceil_+:=   \max\left\{0,\left\lceil \frac{2J}{3}{- \frac{4J_0}{3}}\right\rceil\right\}$.
% and the Whitney constant from the local affine approximation
%{in the current special setting of the horizon model, $L$ is already controlled by $\kappa$
%because of the constraint on $E'$ and the Whitney constant is absolute because the spatial dimension and the degree
%is fixed}.

It remains to bound the error on the cells in $\cC^{an}_{2J{-J_0}}$.
Due to the jump of $f$ across $\Gamma$ we simply take $0$ as
an approximation to  $f$ on the cells in $\cC^{an}_{2J{-J_0}}$.
Thus, we obtain
\begin{equation}\label{eq:estimate02}
\sum_{Q \in \cC^{an}_{2J{-J_0}}}\|f\|^2_{L_2(Q)}
\lesssim \sharp(\cC^{an}_{2J{-J_0}})2^{-3J {-J_0/2}} {M^2}
\lesssim (2^{-2J})  {L\,M^2},
\end{equation}
 with a constant depending only on $\kappa$.
Therefore, \cs{adding} \eqref{eq:estimate01} and \eqref{eq:estimate02}, we obtain
\begin{equation}
\label{goodbound}
\inf_{\varphi \in \Sigma_N}\|f-\varphi\|_{L_2(D_0)}^2 \lesssim L\,M^2 (\log N_0)N_0^{-2},
\end{equation}
with a constant depending only on $\kappa$.
This proves our claim for discontinuity curves exhibiting the special form \eqref{eq:boundary}
under the constraints in \eqref{Eprop} on the restricted domain $D_0$.

We treat next the remaining cells in $\cC_0^B$.
Since there is a uniform finite bound on the number of these squares, depending only
on $\kappa, L$, it suffices to analyze the complexitity of the partitions $\cC'_{j,r}(f,P_0)$
where $P_0\in \cC_0^B$.
Recall that $\cC'_{0,r}$ is the initial partition of an $L$-shaped domain consisting
of $3$ squares of side length $2^{-J_r}$ where we abbreviate for convenience $J_r:= J_0 +r$.
We can estimate the cardinalities of the corresponding partitions $\cC_{j,r}^{an}, \cC_{j,r}^{iso}$
along the same lines as before, taking the respective scalings into account.
First we note that
\beqn
\label{cardCanjr}
\sharp(\cC^{an}_{j,r}) \lesssim L 2^{j/2},\quad j\ge 0.
\eeqn
Since the highest refinement level should not exceed $2J$ we obtain
\beqn
\label{eq:num1r}
\sharp\Big(\bigcup_{j=0}^{2J-J_r}\cC^{an}_{j,r}\Big) \lesssim L \sum_{j=0}^{2J-J_r} 2^{j/2} \lesssim
L 2^{J- J_r/2}.
\eeqn
The analogous isotropic refinement depth $\ell^r_j(Q)$ for $Q\in \cC^{an}_{j,r}$ not intersecting $\Gamma$,
that ensures that $(R^{iso})^{\ell^r_j(Q)}(Q)$ fits into a rectangle of width
$\sim 2^{-\frac J2 - \frac j4}$
and  length
$\sim 2^{-\frac J2 + \frac j4}$,
takes now the form
\beqn
\label{ljr}
\ell_j^r(Q) = \max\,\Big\{\Big\lceil \frac J2- \frac{3j}{4} - J_r\Big\rceil, 0\Big\},
\eeqn
i.e., $\ell^r_j(Q)=0$ when $j\ge 2(J-2J_r)/3$. Accordingly, we obtain
\begin{eqnarray}
\label{eq:num2r}
\sharp(\bigcup_{j = 0}^{2J{-J_r}}\cC^{iso}_{j,r}) &{\le}&  \sum_{j = 0}^{\lceil \frac{2J}{3}{- \frac{4J_r}{3}}\rceil}2^{2(J/2-3/4j  {-J_r})}\sharp(\cC^{an}_{j,r})
+ \sum_{j > \lceil \frac{2J}{3}{- \frac{4J_r}{3}}\rceil} \sharp(\cC^{an}_j)\nonumber\\
& \lesssim &
%L \sum_{j = 0}^{\lceil \frac{2J}{3}\rceil}2^{2(J/2-3/4j)+j/2}
%+ L \sum_{j = \lceil \frac{2J}{3}\rceil+1}^{2J} 2^{j/2} \nonumber\\
%&=&
L \Big(\sum_{j = 0}^{\lceil \frac{2J}{3}{- \frac{4J_r}{3}}\rceil} %2^{J-j/2}
{2^{J-j-J_r}}
+  \sum_{j = \lceil \frac{2J}{3}{- \frac{4J_r}{3}}\rceil+1}^{2J{-J_r}} 2^{j/2}\Big)\lesssim L 2^{J{-J_r/2}}.
\end{eqnarray}

In summary, recalling from \eref{TP0} that
$
\sharp(\cT_J(P_0)) = 1+ \sum_{r=0}^{J-J_0-1} \sharp(\cC'_{J+1,r}(f,P_0))$, and noting that $2J-J_r\ge J+1$
for $0\le r\le J-J_0-1$ we infer from \eref{cardCanjr} and \eref{eq:num2r} that
$\sharp(\cC'_{J+1,r}(f,P_0)) \lesssim L(2^{J/2}+ 2^{J-J_r/2})$ which yields
 \beqn
\label{rcount}
\sharp(\cT_J(P_0))   \lesssim L \Big((J-J_0-1)2^{J/2} +2^{-J_0} 2^J\Big)\lesssim L 2^J,
\eeqn
where the constant depends only on $\kappa$.
Thus, keeping \eref{Nest} in mind,  the cardinality of the partition $\cT_J(f)$,
defined in \eref{TJ}, is bounded by
\beqn
\label{countTJ}
\sharp(\cT_J(f)) \lesssim L \Big(2^{J_0/2} + \sharp(\cC_0^B) 2^{J}\Big)\lesssim L C 2^{J},
\eeqn
where $C$ depends only on $\kappa$.

We next have to adapt the local approximation error estimates \eref{PQ1}, \eref{PQ2}
which amounts to replacing $J_0$ by $J_r$,
\cs{ providing at least one of the following error bounds}
\beqn
\label{PQr}
\|f-P_{Q}(f)\|^2_{L_2(Q)} \lesssim 2^{-\frac{7}{2}j}{2^{-6J_r}}M^2,
\quad
\|f-P_{Q}(f)\|^2_{L_2(Q)} \lesssim %2^{-\frac{7}{2}j-6(J/2-3/4j)}M^2 =
2^{j-3J}M^2,
\eeqn
depending on whether
$j > \lfloor \frac{2J}{3}{-\frac{4J_r}{3}} \rfloor$,
or
$j \le  \lfloor \frac{2J}{3}{-\frac{4J_r}{3}} \rfloor$,
respectively.
Combining this with the above counts, one obtains
\begin{eqnarray}
\label{error-r}
\sum_{Q\in \bigcup_{j=0}^{2J-J_r-1}\cC^{iso}_{j,r}}\|f-P_{Q}(f)\|^2_{L_2(Q)}
&\lesssim &
LM^2 2^{-2J} 2^{-J_r} \Big( 2^{-J_r}
+
\left\lceil \frac{2J}{3}{- \frac{4J_r}{3}}\right\rceil_+\Big).
\end{eqnarray}
By \eref{cardCanjr}, one derives
\beq
\label{Canr}
\sum_{Q \in \cC^{an}_{2J{-J_r}}}\|f\|^2_{L_2(Q)}
\lesssim \sharp(\cC^{an}_{2J{-J_r}})2^{-3J {-J_r/2}} {M^2}
\lesssim 2^{-2J-J_r}  {L\,M^2},
\end{equation}
so that  we obtain the following bound for the total error on $D_r$
\beqn
\label{erroronDr}
\inf_{g\in \PP_1(\cC'_{2J-J_r}(f,P_0))}\|f - g\|_{L_2(D_r)}^2 \lesssim
LM^2 2^{-2J} 2^{-J_r} \Big( 1+ \left\lceil \frac{2J}{3}{- \frac{4J_r}{3}}\right\rceil_+\Big).
\eeqn

In conclusion, the total (squared) error \cs{contributions from}
$\PP_1(\cT_J(f))$ can be bounded by summing over the local (squared) errors on the domains
$D_r=D_r(P_0)$, $P_0\in \cC_0^B$, $r\le J-J_0-1$.
Specifically, we infer from \eref{goodbound} and \eref{erroronDr} that
$$
\inf_{g\in \PP_1(\cT_J(f))}\|f- g\|_{L_2(D)}^2 \lesssim
LM^2 \Big(2^{-2J} + \sharp(\cC_0^B) J 2^{-2J} \sum_{r=0}^{J-J_0-1}2^{-J_r}\Big),
$$
where the constant depends only on $\kappa$.
This proves the claim under assumptions \eqref{eq:boundary}, \eqref{Eprop}  on $f$.

%Note that $\cC^{an}_{2J} \cup \cC^{iso}$ forms a partition of $D$ and all of the refinement operators used to construct $\cC^{an}_j$ and $\cC^{iso}_j$ can be generated by the general refinement operators
%illustrated in Figure \ref{refinement}.

\begin{figure}[htb]
\begin{center}
\includegraphics[width=0.3\textwidth]{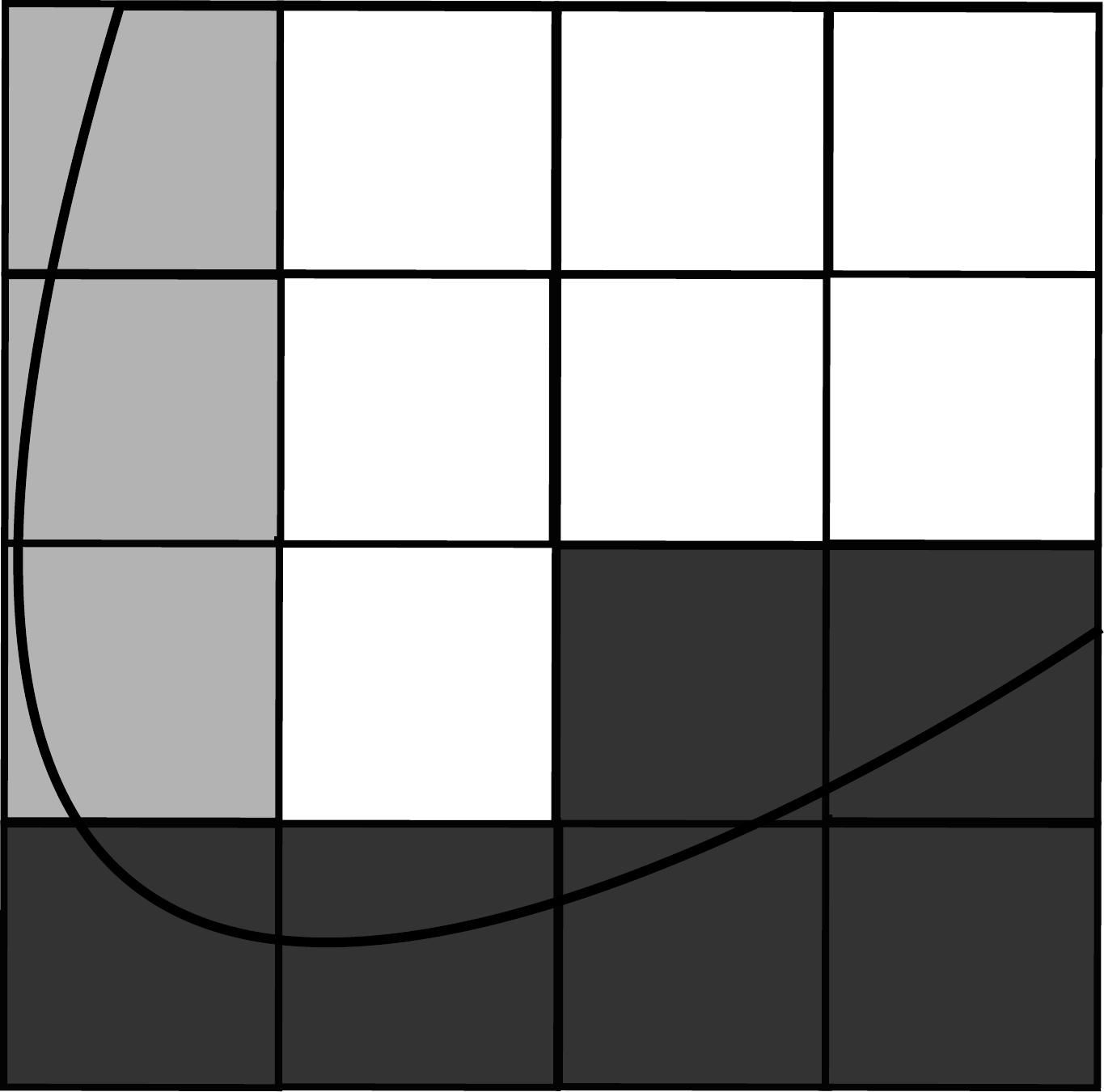}
\end{center}
\caption{Dark gray cells correspond to the region $D^v$
         while light gray cells to $D^h$ \cs{are} defined below.}
\label{fig:hor_ver_curve}
\end{figure}

Next, we prove the claim for a general $f\in \cC(\kappa,M,L,\omega)$.
First note that the above arguments apply in the same way
to a ``vertical horizon model'', i.e., when the roles of the coordinates
$x_1,x_2$ are interchanged.
Depending on $\kappa$ and on $\omega$,
we find next a fixed $I_0\in \N$ such that the partition $\cP_0$
comprising the squares $2^{-I_0}(D + m)$, $0\le m_i<2^{I_0}$, $i=1,2$,
can be be decomposed as
\[
\cP_0 = \cP^{v}_0 \cup \cP^{h}_0
\]
in such a way that the corresponding regions
$$
D^v:= \bigcup\,\{P\in \cP^{v}_0\},\quad D^h := \bigcup\,\{P\in \cP^{h}_0\}
$$
have the following properties:
each connected segment of $\Gamma \cap  D^{v}$, $\Gamma \cap  D^{h}$
can be written as the graph of a function in $x_2$, $x_1$, respectively, with slope bounded
by two, and hence satisfies \eref{Eprop} on each respective
region $\tilde D$ (see Figure \ref{fig:hor_ver_curve}).
Since the number of these regions is bounded with a constant depending
only on $\kappa$ and on  $\omega$,
it suffices to verify the claim for each such subdomain.
Hence,  the claim follows,  from the result for the preceding
special case with $D$ replaced by $\tilde D$,
\cs{with all constants in the above estimates appropriately rescaled.}
\hfill $\Box$\\

We conclude this section with a few remarks on possible extensions. The above cartoon class is certainly not the most general one for which
an analogous result would hold. For instance, one could permit several components where the constants would also
depend on a separation between these components. One could consider a finite number of non-tangential self-intersections
of the singularity curve, treating these self-intersections in a similar way as the intersection of $\Gamma$ with the vertical boundary portions in the above proof.
\cs{
As in \cite{CD}, one could relax the present
global}
$C^2$ continuity requirement on $\Gamma$ \cs{to piecewise continuity with finitely many pieces}
or even \cs{to} $B^2_\infty(L_\infty)$-Besov regularity.
%\dw{[WD: the notation is not really good, too many $\cC's$. Perhaps someone has a better idea.]}

\subsection{Refinement without merging}

The scheme used to prove Theorem \ref{ShearCart} strongly  prioritizes parallelograms in combination with a simple index structure. Moreover, the merging operation
supports a sufficiently rapid directional adaptation. The prize to be paid is to dispense with {\em nestedness} of the resulting partitions which causes
well known drawbacks. In particular, simple multi-resolution concepts, based for instance on Alpert wavelets, are no longer available for non-nested partitions.

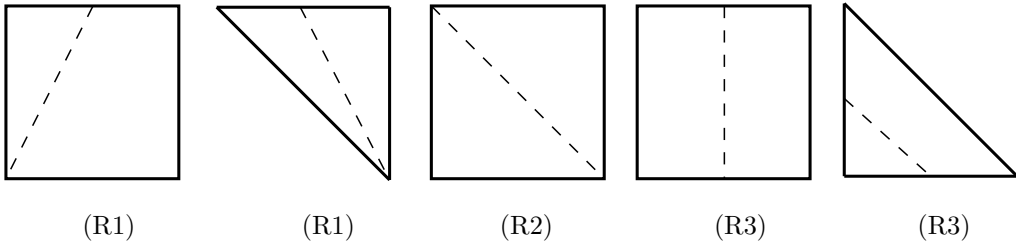
\begin{figure}[ht]
\begin{center}
\includegraphics[width=0.9\textwidth]{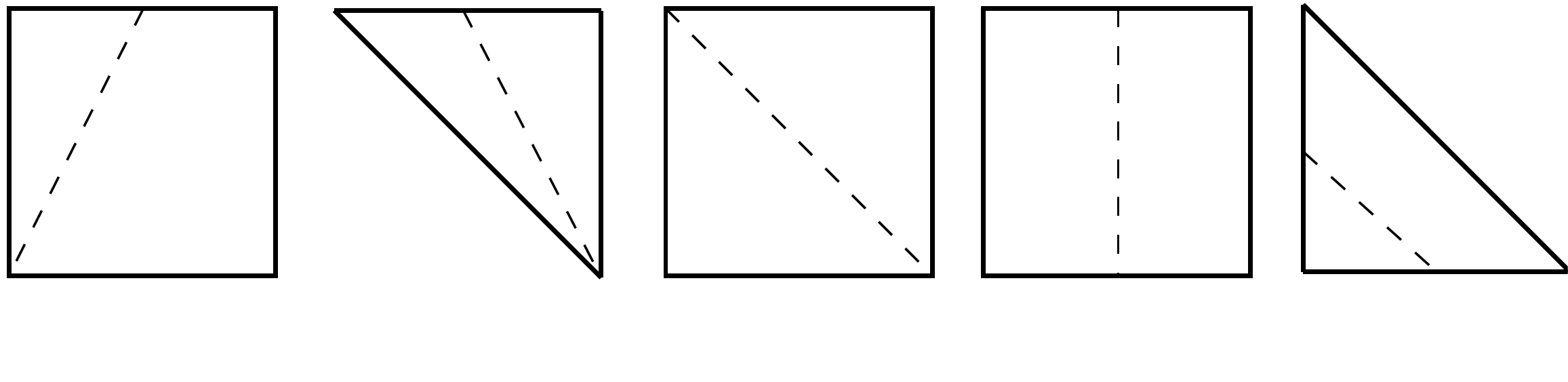}
%\put(-280,5){(1)}
%\put(-215,5){(2)}
%\put(-152,5){(3)}
%\put(-90,5){(4)}
%\put(-30,5){(5)}
\put(-353,5){(R1)}
\put(-270,5){(R1)}
\put(-195,5){(R2)}
\put(-114,5){(R3)}
\put(-37,5){(R3)}
\end{center}
\caption{Illustration of the partition rules.}
\label{refinement}
\end{figure}

As an alternative we discuss next a scheme for generating {\em nested} hierarchies of anisotropic partitions comprised of
triangles and quadrilaterals only, no longer insisting on the quadrilaterals to be parallelograms.
 Compared to the construction in Theorem \ref{ShearCart}, this scheme is is easier to code, especially for more complex domains. It is a bisection scheme using the following types of splits (see also \cite{priv-com,Mi11}).
%***************************************************************************************************************************
%\subsubsection{Basic Refinement Rules}
%\label{sec:BasRefRul}
%
%We start by describing the basic refinement rules,
%which are applied by our subsequent shearlet-inspired refinement scheme.
Starting
from some initial partition  of $D$  consisting of triangles and quadrilaterals,
{\em refined partitions} are obtained by splitting
a given cell $Q$ of a current partition according to one of the following rules, depicted in Figure \ref{refinement}:
\begin{enumerate}
\item[(R1)] Connect a vertex with the midpoint of an edge not containing the vertex.
\item[(R2)] Connect two vertices.
\item[(R3)] Connect the midpoints of two  edges which, when $Q$ is a quadrilateral, do not share any vertex.
\end{enumerate}

We emphasize that that these rules apply to general quadrilaterals as they involve only connecting vertices or midpoints of edges.
It is easily checked that these three refinement rules produce only
triangles and quadrilaterals. Moreover, the refinement operators used  in the proof of Theorem \ref{ShearCart} are covered by these rules as special cases. In fact, $R^{an}_{\iota}$ is generated by applying  one refinement of type (R3) when $\iota = 0$ and two refinements type (R1) when $\iota = \pm 1$. Furthermore, the isotropic refinement operator $R^{iso}$ is generated by  applying three type (R3) refinements when $Q$ is a parallelogram and by (R3) (R1) (R3) when $Q$ is a triangle.

 We briefly indicate now  that the bisection split rules (R1) -- (R3) offer, in principle, a sufficiently rapid direction adaptation
allowing one to create a {\em nested hierarchy} of partitions comprised of triangles and quadrilaterals with the following property: Let $\cT$ denote a partition obtained through $N$ \dww{successive} splits from (R1) -- (R3), \dww{starting from a given initial partition which for simplicity consists just of the unit square.} \dww{As before, let
$\PP_1(\cT)$ denote} the  space of discontinuous, piecewise linear  polynomials
\dww{on $\cT$} and $\Sigma_N$ the \dww{union} of  all such spaces \dww{on such partitions}
generated by any  combination of at most $N$ \dww{successive refinements} (R1) -- (R3).
%\css{[from which starting partition?]}
%
\begin{theorem}
\label{thm:split}
For any $f\in \cC(\kappa,L,M,\omega)$ one has
\beqn
\label{Rsplits}
\inf_{\varphi \in \Sigma_N} \|f- \varphi\|_{L_2(D)} \le C (\log N) N^{-1},\quad N\in \N,
\eeqn
where $C$ depends only on $\kappa,L,M,\omega$.
\end{theorem}
We do not insist in \eref{Rsplits} on the smallest possible power of $\log N$. The purpose of Theorem \ref{thm:split} is again to provide a benchmark for which rates could be achieved when adaptively using the above split rules. As in Theorem \ref{ShearCart} the subsequent proof is constructive but not meant to suggest a practical algorithm because the construction is based on an ideal oracle directing the choice of splits. However,
replacing the ideal oracle used below by some a posteriori criteria, as done later below,   leads to a practical algorithm whose performance
can then be compared with the rates given in Theorem \ref{thm:split}, see e.g. \cite{CoDyHeMi12,CM1} for algorithms with similar refinement rules.
\dww{Of course, the number $N$ of splits is a lower bound for the computational work required by any such numerical realization.}\\

\noindent
{\bf Proof:} The proof proceeds in three steps. First, we refine isotropically until we achieve the desired target error on all cells that do not not intersect the singularity curve. However, the diameter of these cells is still too large to achieve the target accuracies on the cells intersecting the singularity. Therefore, in a second step, we apply further isotropic refinements to these cells. Finally, in a third step, we ensure the correct parabolic scaling by applying purely anisotropic refinements to the remaining cells intersecting the singularity curve $\Gamma$.\\

\paragraph{\it Step 1}
We begin with the isotropic partition $\cP_{J/2}:= (R^{iso})^{J/2}(D)$ (where we assume for simplicity that $J$ is even). Here we assume  that $J$ is large enough
to ensure that  every cell boundary intersects the discontinuity curve just twice, i.e.} for every $Q\in \cP_{J/2}(D)$ with $|Q\cap\Gamma| >0$ we have  $\#(\Gamma \cap \partial Q)=2$.

Let $\cP_{J/2}(\Gamma):= \{Q\in \cP_{J/2}: |Q\cap \Gamma|>0\}$  be the collection of  the cells that intersect the singularity curve on level $J/2$. Then we have
\beqn
\label{err-case-2}
\inf_{g\in \PP_1}\|f-g\|_{L_2(Q)}^2 \le M^2 2^{-2J} |Q|
\eeqn
for all $\,Q\in \cP_{J/2}\setminus \cP_{J/2}(\Gamma)$.
The total projection error on the union of those squares outside $\cP_{J/2}(\Gamma)$ is bounded by $M 2^{-J}$.
\\

\paragraph{\it Step 2}
 In order to obtain the desired diameter for the cells intersecting the singularity, we successively subdivide isotropically only the cells intersecting $\Gamma$ until reaching level $J$. For each refinement level $j$ with $J/2 \le j \le J$ we have
\beqn
\label{n-case-2}
\#(\cP_{j}(\Gamma))\le C 2^{j}
\eeqn
cells intersecting $\Gamma$, where $C$ depends only on $L,\kappa,\omega$. Each isotropic refinement generates $4$ sub-cells so that we obtain at most
\[
  \sum_{j=J/2}^{j<J} 4 C 2^j \le 4 C 2^J
\]
new sub-cells. Each cell that does not intersect the singularity satisfies the error bound \eqref{err-case-2}.\\

%%For simplicity we do not insist here on  In fact, t
%To verify the above claim we begin with the isotropic
%partition $\cP_{J/2}:= (R^{iso})^{J/2}(D)$ (where we assume for simplicity that $J$ is even). Here we assume  that $J$ is large enough
%to ensure that for every $Q\in \cP_{J/2}(D)$ with $|Q\cap\Gamma| >0$ we have \gw{$\#(\Gamma \cap \partial Q)=2$} and for
%$\cP_{J/2}(\Gamma):= \{Q\in \cP_{J/2}: |Q\cap \Gamma|>0\}$ we have
%\beqn
%\label{G1}
%\#(\cP_{J/2}(\Gamma))\le C 2^{J/2},
%\eeqn
%where $C$ depends only on $L,\kappa,\omega$ while
%\beqn
%\label{G2}
%\inf_{g\in \PP_1}\|f-g\|_{L_2(Q)}^2 \le M^2 2^{-2J} |Q|,\quad Q\in \cP_{J/2}\setminus \cP_{J/2}(\Gamma)),
%\eeqn
%that is the total projection error on the union of those squares outside $\cP_{J/2}(\Gamma))$ is bounded by $M 2^{-J}$.
%
%
%So it remains to generate suitable further subdivisions of the elements in $\cP_{J/2}(\Gamma)$. The first step is to further
%subdivide only the elements in $\cP_{J/2}(\Gamma)$ isotropically $J/2$ times to arrive at
%$$
%\cP_J(\Gamma) = (R^{iso})^{J/2}(\cP_{J/2}(\Gamma)).
%$$
%By \eref{G1} we have
%\beqn
%\label{G3}
%\#(\cP_J(\Gamma)) \le C 2^J,
%\eeqn
%with the above dependencies of $C$.

\paragraph{\it Step 3} Finally, we apply some anisotropic refinements to achieve the parabolic scaling. To simplify the construction note that on each $Q\in \cP_J(\Gamma)$ we can replace $\Gamma \cap Q$ by the straight line connecting the two intersection points $\{p,q\}= \partial Q\cap\Gamma$. In fact the squared error on $Q$ caused this way is easily seen to be bounded by \dww{$\bar CM^2 2^{-3J}$ where $\bar C$ depends on the curvature bound
$\kappa$. Hence, by \eqref{n-case-2}, for $J$ larger than $J_0$, depending on $L,\kappa,\omega$, the total squared error on the union of all $Q\in \cP_J(\Gamma)$ is bounded by $\bar CM^2 2^{-2J}$ where $\bar C$ depends on $\kappa$}.

From now on we will therefore assume that $\Gamma$ is a polygon with vertices $p=p(Q), q=q(Q)$ on the boundaries $\partial Q$, i.e., $Q\cap \Gamma$ is a straight line segement contained in $Q$ with end points $p=p(Q), q=q(Q)$
located on $\partial Q$. The goal of the subsequent anisotropic refinements is to create a cell $Q(\Gamma)$ which contains the whole
line segment $\Gamma \cap Q$ and has area \dww{$\le  2^{-3J}$}, i.e.,
%\css{[implied constant in $\lesssim$ depends on $K,L,M,\omega$]}\dww{[Comment WD:  I don't think so, , the constant in the above inequality is simply one and does not
%depend on the parameters $\kappa, \omega,L$. A terminal quadrilateral resulting from the split strategy has longest edge $\le \sqrt{2}2^{-J}$ and
%both other edges $2^{-2J}$. Its area is easily seen to be bounded by $2^{-3J}$.]}
\beqn
\label{Qgamma}
\Gamma \cap Q \subset Q(\Gamma),\quad |Q(\Gamma)| \dww{\le 2^{-3J}}.
\eeqn
Since there are \dww{$\le  4C2^J$} elements in $\cP_J(\Gamma)$,
the total squared error on those squares is then
\[
\sum_{Q\in \cP_J(\Gamma)}\inf_{g\in {\PP}_1}\|f-g\|^2_{L_2(Q(\Gamma))}
\le \dww{4C }
M^2 2^{-2J} \sum_{Q\in \cP_J(\Gamma)} \|f \|^2_{L_2(Q)},
\]
and hence of the desired  same order as the total squared error on all cells that do not intersect $\Gamma$. The proof of Theorem \ref{thm:split}
follows as soon as we have shown that the number of refinements needed to achieve \eref{Qgamma} is
bounded by a fixed absolute multiple
%\css{(depending on $L$, $\kappa$ and on $\omega$)}
of $J$.
%\dww{[Comment WD: I don't think so, this multiple is independent of $L$, $\kappa$ and on $\omega$, see the subsequent argument, so this should be removed.}

We will show next that \eref{Qgamma} can indeed be achieved at the expense of the order of $J$ further refinements. These refinements
will create a sequence of cells $Q_k$ with $Q=Q_0$, always containing $\Gamma \cap Q$, which are either quadrilaterals or triangles.
We refer to the edges $E^k_p, E^k_q$ of $Q_k$  containing $p,q$, respectively, as split-edges. One then encounters one of the  following   two cases:\\

\begin{figure}[ht]
\centering
\includegraphics[width=0.4\textwidth]{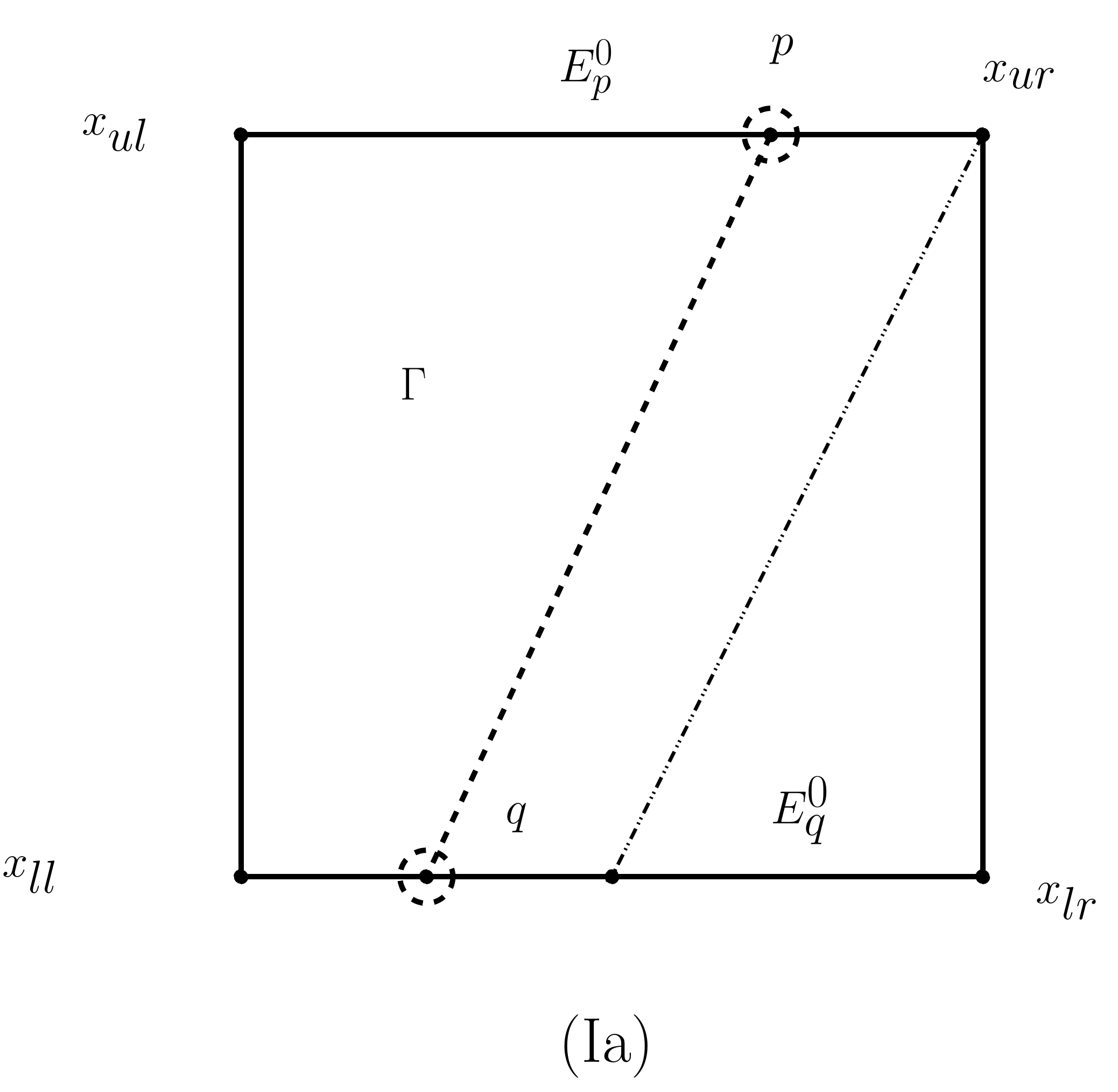}~
\includegraphics[width=0.38\textwidth]{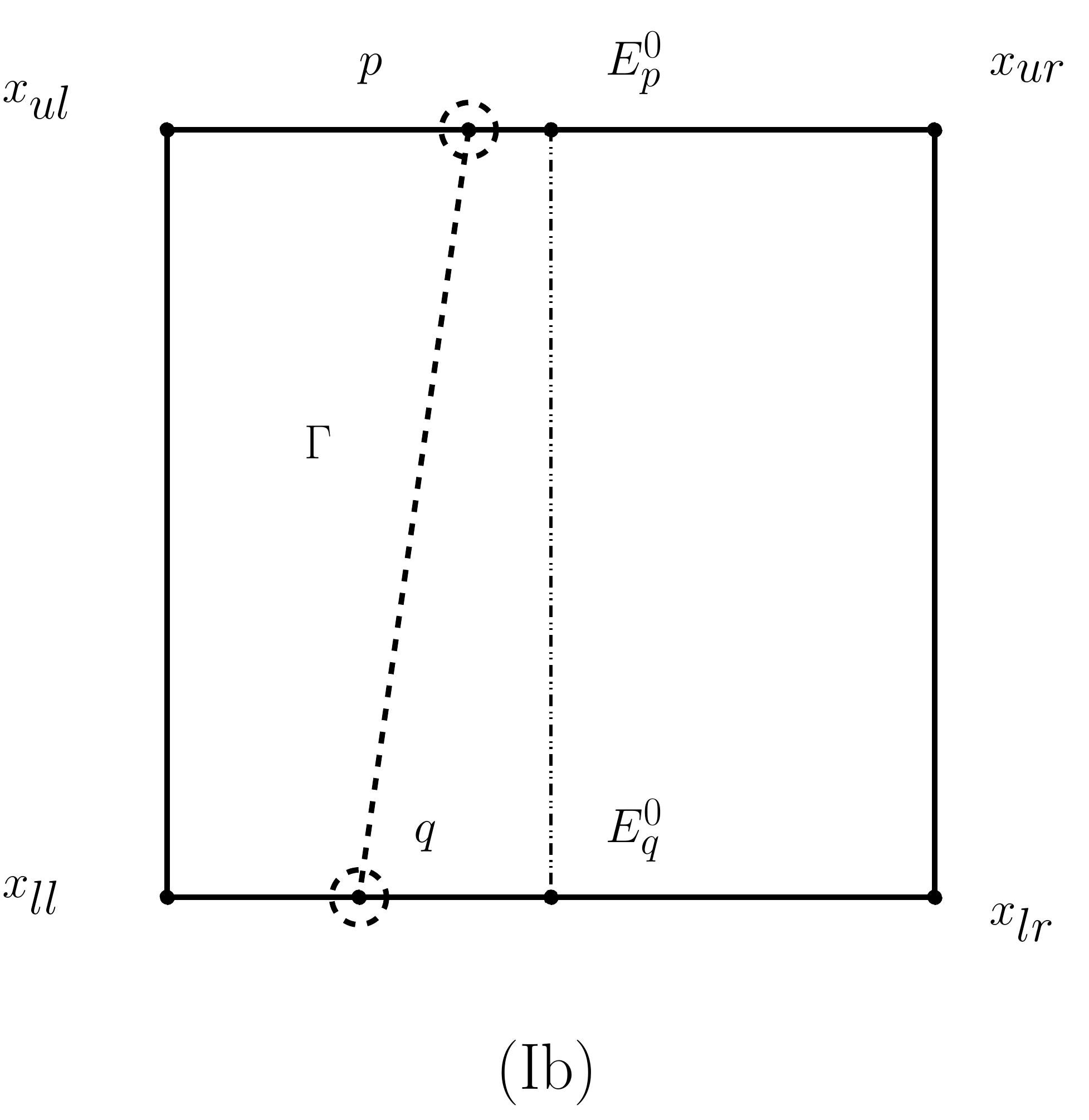}
\caption{(Ia): split type (R1);\hspace*{3mm} (Ib): split type (R3)}
\label{Fig:I1}
\end{figure}

(I): The endpoints $p,q$ of $\Gamma \cap Q_k$ are located on a pair of {\em opposite} split-edges $E^k_p, E^k_q$ of   the quadrilateral $Q_k$, i.e., these edges do not have any
  vertex of $Q_k$ in common, see the two cases in Figure \ref{Fig:I1}.\\

\begin{figure}[ht]
\centering
\includegraphics[width=0.32\textwidth]{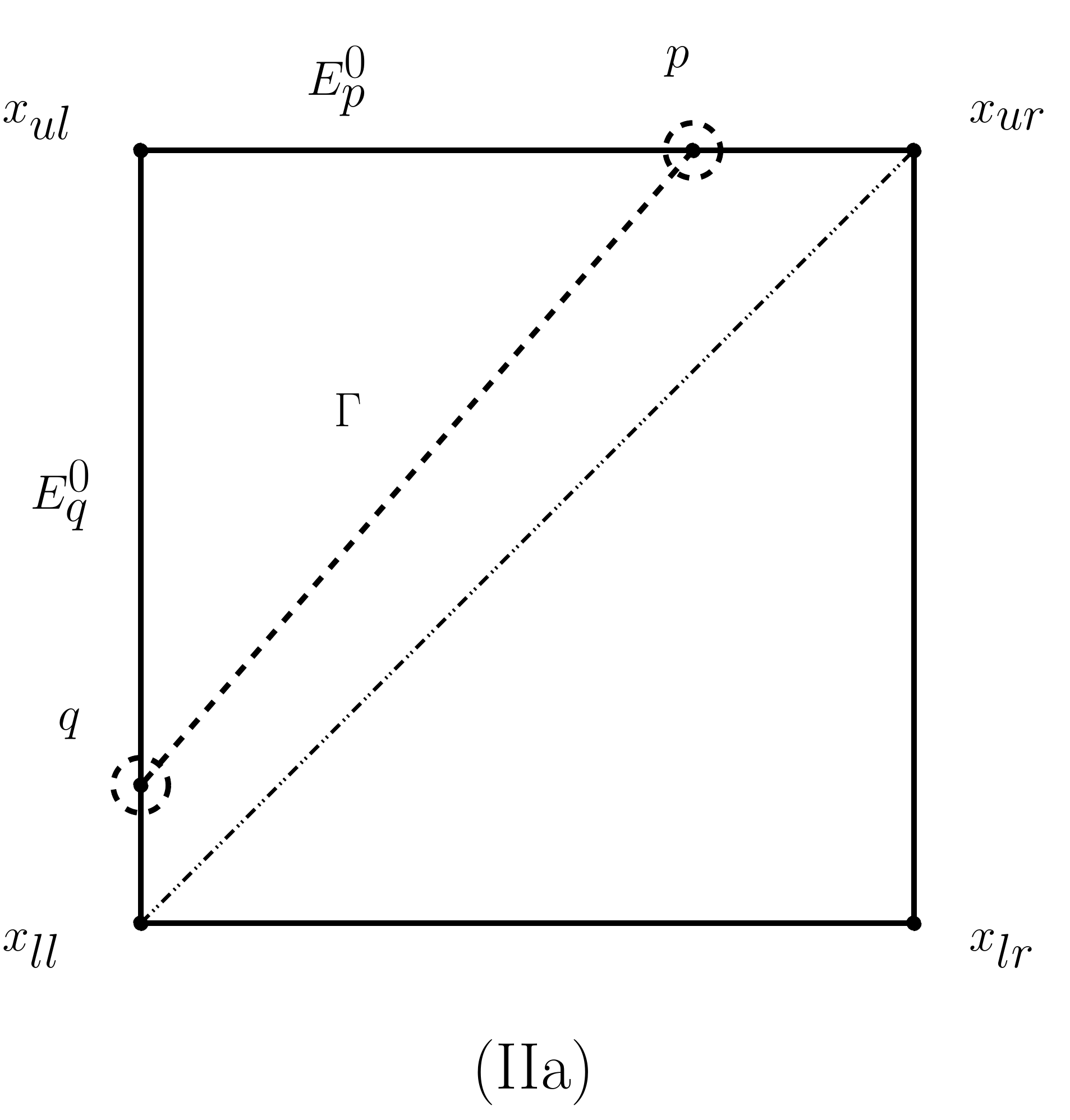} ~
\includegraphics[width=0.32\textwidth]{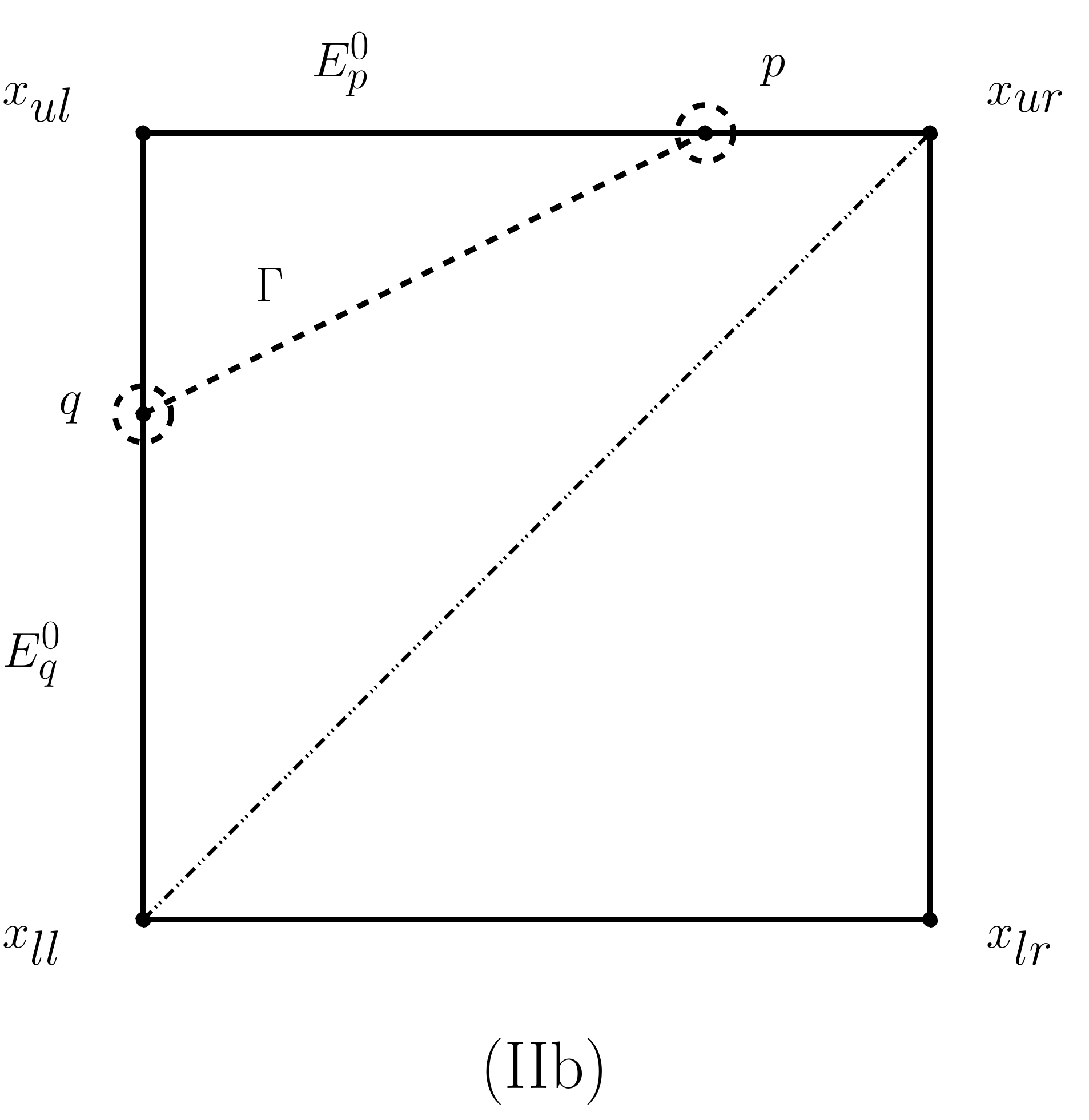} ~
\includegraphics[width=0.32\textwidth]{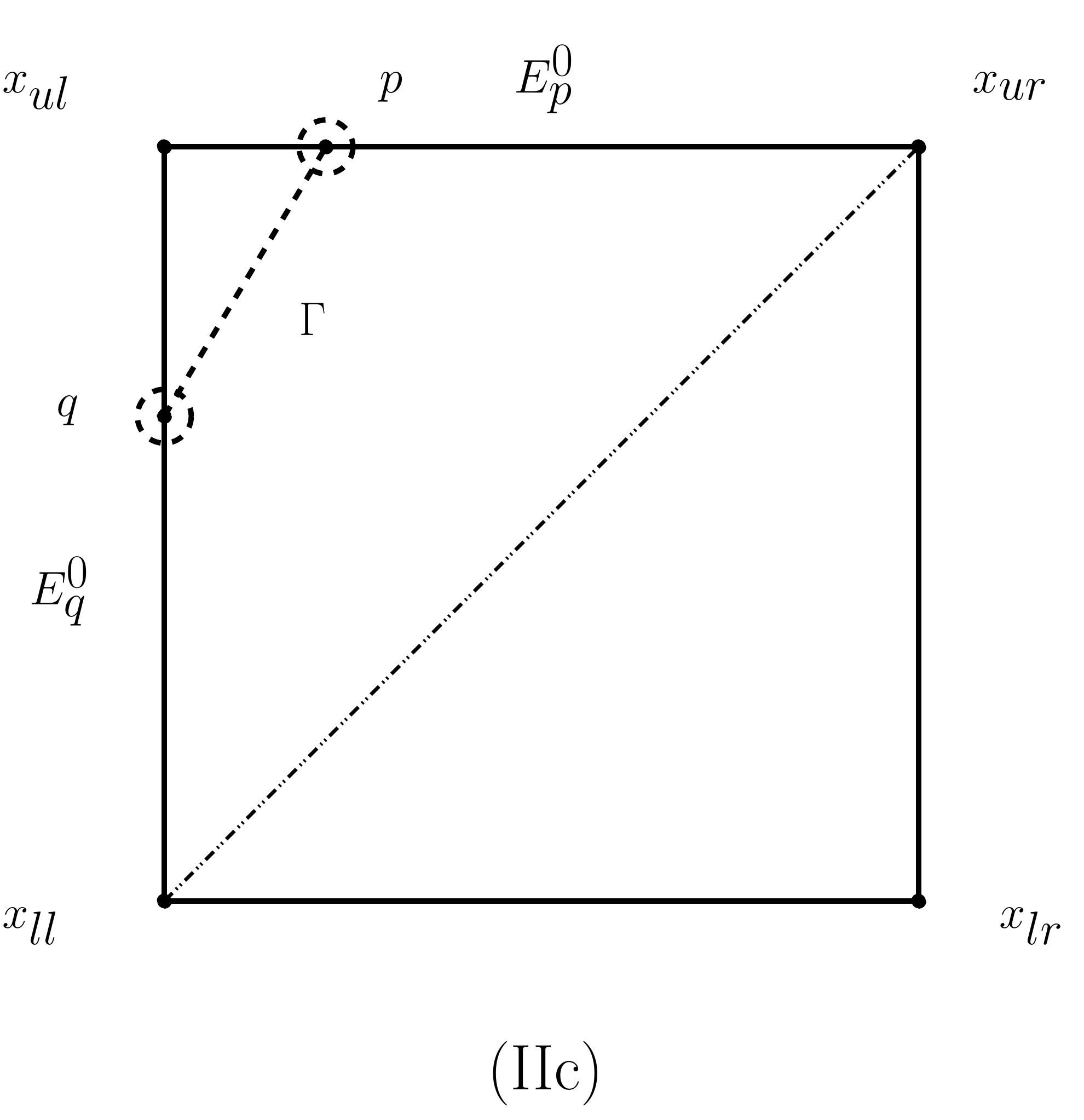}
\caption{$Q_1$ is a triangle after (R2) in the case (II)}
\label{Fig:II}
\end{figure}

(II): The split-edges $E^k_p, E^k_q$ of $Q_k$, share a common vertex and  $Q_k$ is either a triangle or a quadrilateral, see Figure \ref{Fig:II}.\\

We describe  first the refinements for case (I).   Figure \ref{Fig:I1} shows for    the initial cell $Q=Q_0$ all possible situations
modulo a rotation or swapping coordinates.
  To simplify the explanation we call the  endpoint of $E^k_p, E^k_q$ which has
the larger distance from $p,q$, respectively, its {\em far-vertex}. We call this distance the {\em far-distance}.
In these terms we apply the following  rule: \\

\noindent
{\bf IR:} Apply either (R1) or (R3) according to the following conditions:
\begin{itemize}
\item[(IR1)] If the far-vertex of $E_p$ and the far-vertex of $E_q$ are {\em not} connected by a common edge of $Q_k$ apply (R1)
by splitting  the longest split-edge as indicated in Figure \ref{Fig:I1}, (Ia)
(this leaves $\Gamma\cap Q_k$ within the new quadrilateral $Q_{k+1}$).
%  (If both split-edges have equal length we split the one with
%the larger far-distance.)
\item[(IR3)] If the far-vertices of $E^k_p, E^k_q$ are the endpoints of a {\em common} edge of $Q_k$, apply (R3) splitting both $E^k_p$ and
$E^k_q$, see Figure \ref{Fig:I1}, (Ib).
\end{itemize}
The rule says that one chooses the split that reduces as many far-distances as possible. If one can split only one split-edge one chooses the larger one.
The following properties are then an immediate consequence:\\

\noindent
{\bf IPr:} Given $Q_k$ of type (I), $Q_{k+1}$ still contains $\Gamma \cap Q$ and is also of type (I). Moreover, either both split-edges are halved, i.e., $|E_p^{k+1}|=
|E_p^{k}|/2$ and $|E_q^{k+1}|=
|E_q^{k}|/2$, or the longer one of the split-edges   is halved, i.e.,
$|E_p^{k+1}|=
|E_p^{k}|/2$ and $|E_q^{k+1}|=
|E_q^{k}|$.\\

 Hence, since the splits aim at reducing far-distances, after at most $2J$ such refinements one has $|E^J_p|\le 2^{-2J}$,  $|E^J_q|\le 2^{-2J}$,
 i.e., $|Q_J| \lesssim 2^{-3J}$, as desired. Thus, in this case the refinement process terminates after at most $2J$ splits and achieves \eref{Qgamma}.

 We now turn to case (II) and start  with the refinement for $Q_0=Q$. First we apply (R2) so that $\Gamma$ is now contained in
the triangle $Q_1$ spanned by the vertices $x_{ll} , x_{ul} , x_{ur}$, \dww{see Figure \ref{Fig:II}, (IIa)}.
More generally, whenever $Q_k$ is a triangle, up to rotation or swapping coordinates, there are three possible
configurations as indicated in Figure \ref{Fig:II}, namely, denoting again by $x_{ul}= E^k_p\cap E^k_q$ the
common vertex of the two split-edges:\\

(IIa): $x_{ul}$ is the far-vertex for both $p$ and $q$, see case (IIa) in Figure \ref{Fig:II};

(IIbc): the far-vertex of at least one of the $p,q$ is not $x_{ul}$, see the cases (IIb), (IIc) in  Figure \ref{Fig:II}.\\

Again the idea is to choose the split that best reduces the far-distances of $p$ and $q$.

\begin{figure}[ht]
\centering
\includegraphics[width=0.32\textwidth]{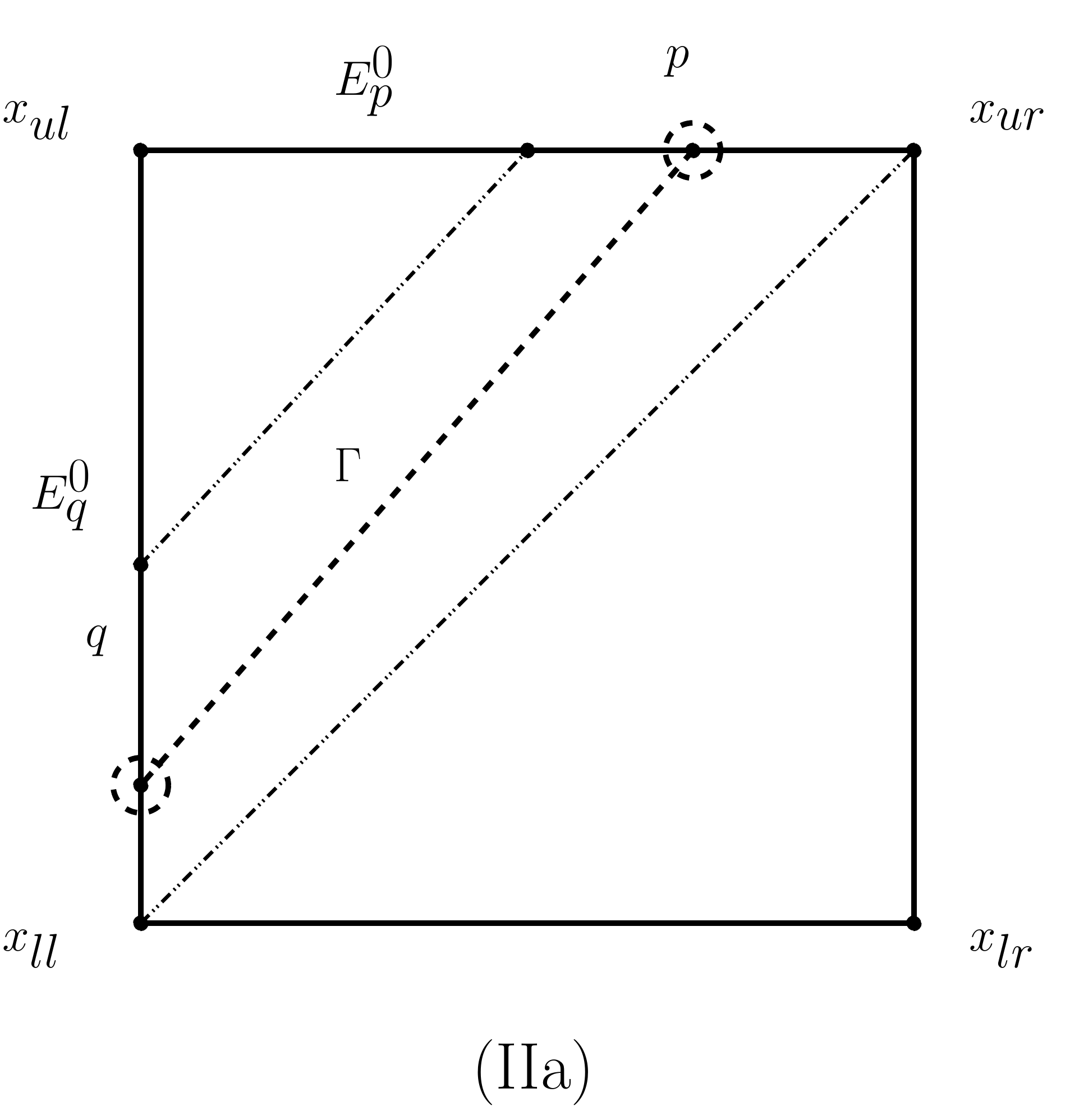} ~
\includegraphics[width=0.32\textwidth]{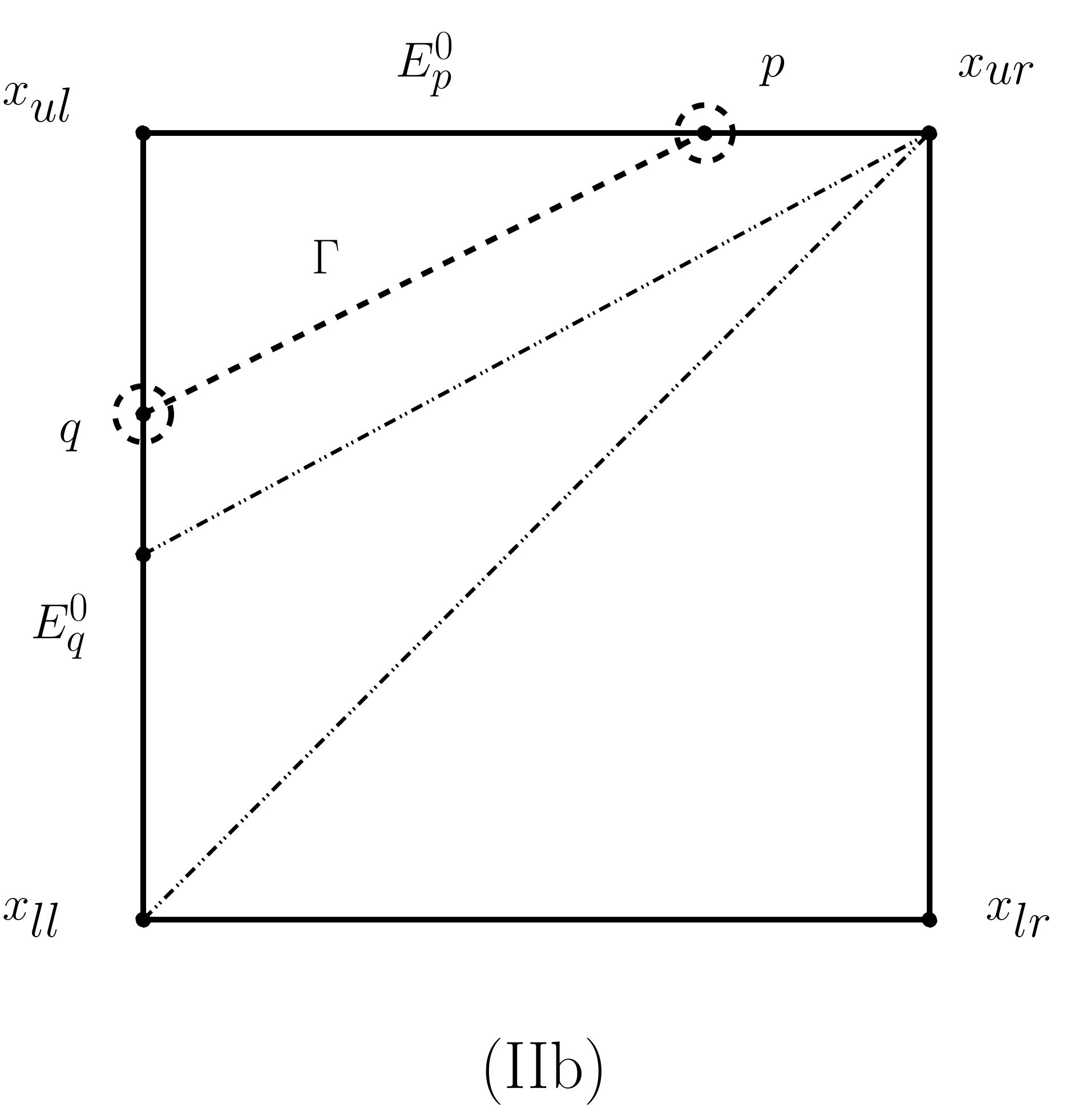} ~
\includegraphics[width=0.32\textwidth]{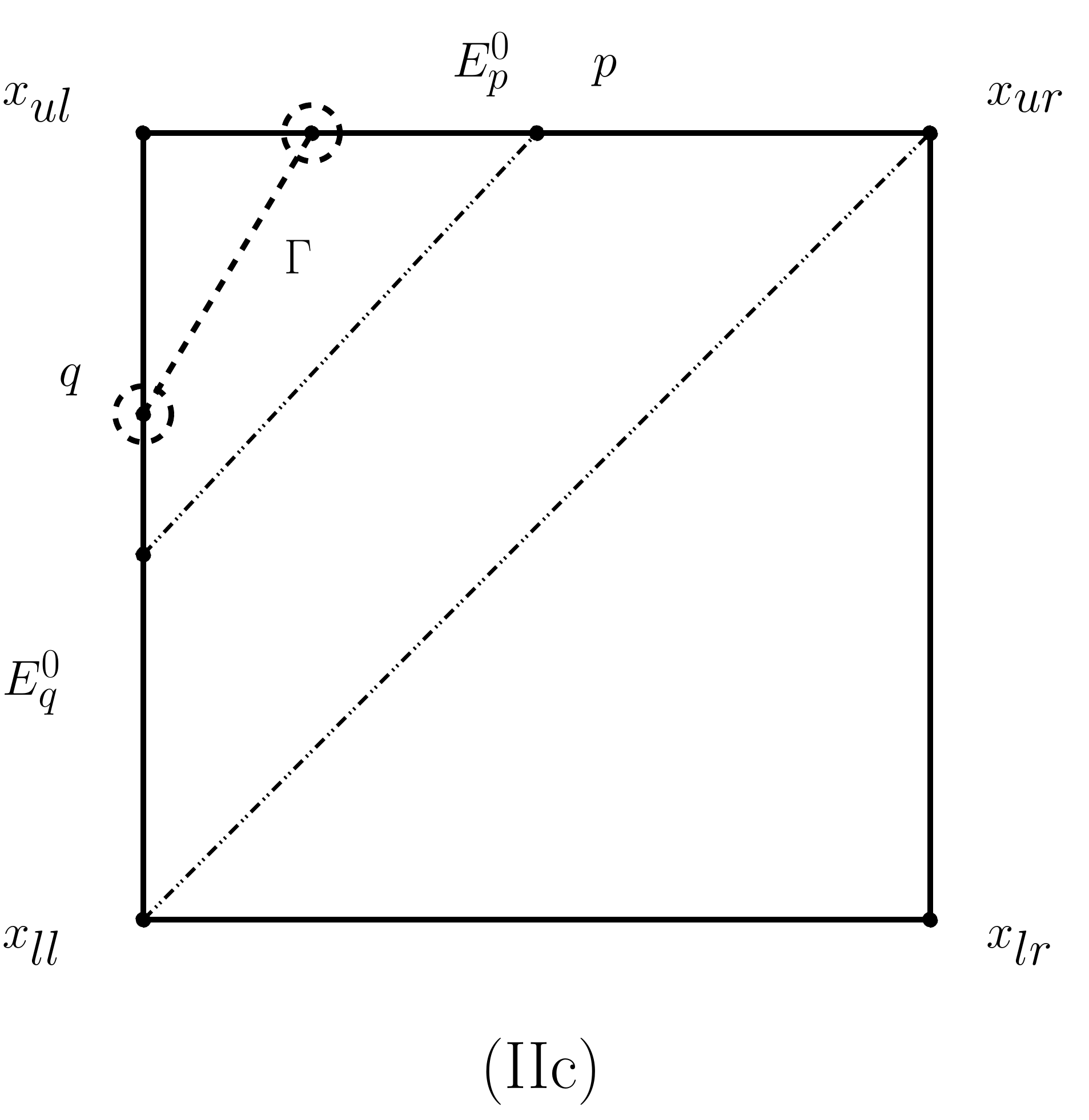}
\caption{Next split generating $Q_2$ in the cases (IIa), (IIb), (IIc)}
\label{Fig:II2}
\end{figure}

In the case (IIa) \dww{of Figure \ref{Fig:II}}   we can apply the split (R3), cutting off the far-vertex $x_{ul}$ and  leaving
$\Gamma \cap Q$ in a quadrilateral $Q_2$, where the two far-distances have at least been halved, see Figure \ref{Fig:II2}, (IIa).
$Q_2$ is now of type (I), so at most $2(J-1)$ further   refinements, as described for case (I), are needed
to produce $Q(\Gamma)=Q_{2J}$ satisfying \eref{Qgamma}. Thus, in this case the refinement procedure
terminates as well with \eref{Qgamma}  after at most  $2J$  splits.

In the case (IIbc) (see cases (IIb) and (IIc) in Figure \ref{Fig:II2}) we can either apply the split (R1) if only one of the far-vertices is different from $x_{ul}$,
splitting the midpoint of the respective split-edge, or one applies (R3) if $x_{ul}=E^k_p\cap E^k_q$ is not a far-vertex
for either one $p$ or $q$, see Figure \ref{Fig:II2}, (IIc). In both cases $Q_{k+1}$ contains $\Gamma \cap Q$ and is again a triangle and hence is of type
(IIa) or (IIbc).
Moreover, since at least one split-edge is halved one has $|Q_{k+1}|\le |Q_k|/2$. Thus, one can have at most $J$ steps
reproducing (IIbc) because such a $Q_J$  would have area at most $2^{-J}2^{-2J}=2^{-3J}$.  In this case $Q_J$ would satisfy \eref{Qgamma}
and the process terminates.
So suppose after $J_0<J$ refinements  $Q_{J_0}$ is of type (IIa).
As shown above, one additional (R3) split produces a quadrilateral containing $\Gamma \cap Q$ and hence
returns to case (I). Moreover,  by the above properties, the product of the length of the split-edges of $Q_{J_0}$ is at most $2^{-J-J_0}$.
By the properties {\bf IPr} of the process for case (I) at most $2J-J_0$ additional refinements are needed
to reduce both split-edges to length at most $2^{-2J}$. The resulting $Q_{2J}$ again satisfies \eref{Qgamma}
which completes the proof.
 \hfill $\Box$\\

 In \eref{Rsplits} the error decay is measured in terms of the number of subdivisions which indicates the work needed to generate
the partition. On the other hand, the number of splits is proportional to the number of generated cells and thus gives
an estimate which is comparable with the one given in Theorem \ref{ShearCart}.
In particular, the refinement scheme of this section offers more possibilities with regard to the shape of the cells, but removes the merging in favor of nestedness. Although the results are eventually not sharp, for the given arguments the additional merging allows a slightly better log factor $(\log N)^{1/2}/N$ in Theorem \ref{ShearCart} as opposed to $(\log N)/N$ in Theorem \ref{thm:split}.
 {The additional log factor originates from the anisotropic refinements necessary for the direction adaption in Step 3. However,} after generating the partition one could merge all  those cells in $Q\in \cP_J(\Gamma)$
that are generated in the  intermediate directional refinements and which do not intersect
$\Gamma$, retaining only  $O(2^J)$ cells.
Revisiting the different cases in Step 3, it is easy to see that these merged cells are
either triangles or quadrilaterals, still forming a nested hierarchy.
Thus an error of order $N^{-1}$ is obtained with  $O(N)$ degrees of freedom \dww{any numerical algorithm based on the split rules (R1) -- (R3) would
require a number of operations of at least the order $N \log N$.}

%\css{[This phrase is repeated; remove here]} \dww{[Comment WD: I agree this has been said before.}
%On the other hand, the above simple construction relies on a perfect oracle
%identifying the intersection of $\Gamma$ with the boundaries
%of the cells in $\cP_J(\Gamma)$ which is not realistic.  Therefore,
%the main purpose of the above considerations
%is to have  {further benchmark} rates  for the nested case, in addition to Theorem \ref{ShearCart},
%that can be compared with the rates obtained by the adaptive method described in Algorithm \ref{alg:uzawa}
%for transport problems.
%

%
\subsection{Adaption for transport problems}
\label{sec:shearlet-implementation}
%
%************************************************************************************************************
In this section we combine the general adaptive scheme Algorithm \ref{alg:uzawa} from Section \ref{sect2.3}
with hierarchies of trial spaces based on anisotropic refinements of the type studied in Section \ref{sect:wang-Q}.
To this end, we have to find realizations of the functions $\stab$ and $\adapt$ which we describe next.
The trial spaces $X_j$ are spanned by piecewise linear functions on a mesh $\cM_j$
composed of cells from collections $\tilde{\cC}_j$, i.e.,
\beqn
\label{Xj}
X_j = \PP_1(\tilde \cC_j),\quad j\geq 0.
\eeqn
The collections $\tilde\cC_j$ consist of anisotropically and
isotropically refined cells of the type introduced in the previous subsection.
Given $X_j$ of the form \eqref{Xj}, we define the {procedure} $\stab$ by
\begin{align*}
  \stab(X_j) & = \PP_2(\tilde\cM_j )\cap C(D) & & \text{where}
& & \tilde\cM_j := \{R^{iso}(P): P \in \tilde \cC_j \}
\;.
\end{align*}
It will be seen that $\stab(X_j)$ is large enough to ensure
$\delta$-proximality with constant $\delta$ significantly smaller than one.

We next explain the anisotropic adaptive refinement
leading to an implementation of the adaption routine $\adapt$ in Algorithm \ref{alg:uzawa}.
The idea is to use a greedy strategy based on the largest ``fluctuation coefficients''.
To describe this, we first initialize
$\tilde\cC_j$ by $\tilde\cC_0 = \cC_0$ defined as in \eqref{eq:initial_cell}.
{The initial isotropic refinement level $J_0$ is typically a small integer but could
depend on a priori knowledge about the convection field $\vb$.}
%\dw{[WD: I describe below also the bulk chasing version specified earlier in the paper. But in the algorithm we actually
%do not specify this and could use the ``maximum criterion''. You can choose and discard one of the two versions.
%Also I didn't understand your notation and modified the representation somewhat.]}

Since we have to rely on a refinement criterion that does not have the information used in
the proof of Theorem \ref{ShearCart} we cannot guarantee that any singularity curve
intersects only parallelograms.
Therefore, given $\tilde\cC_{j-1}$, $j$ odd, we associate with every $Q\in \tilde\cC_{j-1}$
the {\em collection} $\cR^{an}_j(Q)$ of anisotropic refinements $R_j(Q)$ where $R_j(Q)$
is defined as follows:
$R_j(Q) = R^{an}(Q)$ runs through the anisotropic refinements
$R^{an}_{\iota}(Q)$, $\iota = 0,\pm 1$,
defined in \eref{Ran}, when $Q$ is a parallelogram.
When $Q$ is a triangle, \wq{we apply $R^{iso}$ instead.} %$R^{an}$ runs over the three splits (R1) so that the resulting cells are
%always parallelograms or triangles.
When $j$ is even $R_j$ agrees with $R^{iso}$, which
is defined for triangles and parallelograms, producing a single isotropic refinement of $Q$
involving either only triangles or parallelograms.
\wq{We should point out that anisotropic refinement $R^{an}$ is applied only for parallelograms and we always apply isotropic refinement $R^{iso}$ for triangles. }

Given $Q\in \tilde\cC_{j-1}$, we denote for every $R_j(Q)$ by $\Psi_{R_j(Q)}$ an $L_2$-orthonormal basis for
$\PP_1(R_j(Q))$. Recall that $A^*r^K_j \in L_2(D)$.
For $j$ odd, we have to select from several possible anisotropic refinements by
\beqn
\label{Rstar}
R^*_j(Q) := \argmax\,\{\|P_{R_j(Q)}A^* r^K_j\|_{L_2(Q)}: R_j(Q)\in \cR^{an}_j(Q)\},
\eeqn
where $P_{R_j(Q)}$ is the $L_2$-orthogonal projector to $\PP_1(R_j(Q))$. When $j$ is even we have $R_j^*(Q)=R_j(Q)$.
Then, fixing some $\theta \in (0,1)$ we define the collection $\cM$ of {\em marked} cells in $\tilde\cC_{j-1}$ by
%\beqn
%\label{marked}
%\cM(\tilde\cC_{j-1}):= \argmin\,\big\{\sharp(\cS): \cS \subset \tilde\cC_{j-1},\, \sum_{Q\in \cS}\|P_{R^*_j(Q)}A^* r^K_j\|^2_{L_2(Q)}
%\ge \theta \sum_{Q\in \tilde\cC_{j-1}}\|P_{R^*_j(Q)}A^* r^K_j\|^2_{L_2(Q)}\big\}\gk{.}
%\eeqn
%\gk{An alternative choice for $\cM$, which is more convenient for an implementation -- and which we will also
%employ in our numerical experiments -- is}
%
%\dw{[WD: alternatively for the maximum criterion]
\beqn
\label{marked}
\cM(\tilde\cC_{j-1}):=  \{ Q\in  \tilde\cC_{j-1} : \|P_{R^*_j(Q)}A^* r^K_j\|_{L_2(Q)}
\ge \theta \max_{Q\in \tilde\cC_{j-1}} \|P_{R^*_j(Q)}A^* r^K_j\|_{L_2(Q)}
\}.
\eeqn
Then set
\beqn
\label{refinej}
\tilde\cC_j^0 := \mbox{MERGE}(\{Q'\in R^*_j(Q) : Q \in \cM(\tilde\cC_{j-1})\}),\quad
 \tilde\cC^1_j := \{Q \in \tilde\cC_{j-1} : Q \notin \cM(\tilde\cC_{j-1})\},
 \eeqn
providing the refined partition
\[
 \tilde\cC_{j} := \tilde\cC^0_j \cup \tilde\cC^1_j.
\]
\wq{Here, we can not directly apply the operator $\mbox{MERGE}$ as defined in Lemma \ref{lem:merge} since the discontinuity curve $\Gamma$ is not known in this case. For our numerical scheme, we newly define $\mbox{MERGE}(\cC)$ by a set of cells obtained by merging two triangles created from $R^{an}$ to form a parallelogram if there are any such triangles $Q_1,Q_2 \in \cC$ while keeping other cells $Q \in \cC$ as they are.}

Finally, we have to implement $\dataspace$. For simplicity, we use the trivial choice
\begin{equation*}
\dataspace(f, \epsilon, F_h, X_h) = (f, X_h),
\end{equation*}
which essentially means that we ignore data  approximation errors.
\wq{Our numerical adaptive scheme $\solve$ is summarized in {\bf Algorithm 2}.}
\wq{
\begin{algorithm}[htb]
  \caption{2DSOLVER}
  \label{alg:2d_solver}
  \begin{algorithmic}[1]
    \State Initialization: \css{for given target $L_2(D)$-accuracy
           $\epsilon>0$ choose} an initial trial space
           $X_0= \PP_1(\tilde \cC_0)$.
           Set the initial guess, initial error bound and test space
    \begin{align*}
      u_a & = 0, & errbound & =\|f\|_{Y'}, & Z_0 = \stab (X_0),
    \end{align*}
    respectively. Choose parameters $K \in \N$ and $j = 0$.
    \While{$ errbound > \epsilon$}
      \State Compute $\hat u_a$, $\hat r_a$ as the result of $K$ Uzawa iterations with initial value $u_a$, right hand side $f$ and test and trial spaces $X_j$ and $Z_j$.
      \State Compute $\tilde \cC_{j+1}$ from $\tilde \cC_j$ by \eqref{refinej}.
      \State $X_{j+1} := \PP_1(\tilde \cC_{j+1})$
      \State $Z_{j+1} := \stab(X_{j+1})$
      \State Set $j \rightarrow j+1$ and update $u_a$ by
			\begin{align*}
			u_a := \hat u_a + \sum_{Q \in \tilde \cC_{j}}P_{Q}A^*\hat r_a.
			\end{align*}
			\State Update $errbound$.
			\begin{align*}
			errbound := \Big(\sum_{Q \in \tilde \cC_{j}}\|P_{Q}A^*\hat r_a\|^2_{L_2(Q)}\Big)^{1/2}.
			\end{align*}
     \EndWhile
  \end{algorithmic}
\end{algorithm}

%\begin{algorithm}[htb]
%  \caption{2DSOLVER}
%  \label{alg:2d_solver}
%  \begin{algorithmic}[1]
%    \State Initialization: Choose a target accuracy $\epsilon$ and initial trial and test spaces $X_0= \PP_1(\tilde \cC_0)$ and $Z_0 = \stab(X_0)$. Set the initial guess and initial error bound $u_a = 0$, $ \text{errbound} = 2\epsilon$.
%    \While{$\text{errbound} > \epsilon$}
%      \State Compute $\hat u_a$, $\hat r_a$ as the result of $K$ Uzawa iterations with initial value $u_a$, right hand side $f$ and test and trial spaces $X_j$ and $Z_j$.
%      \State Compute $\tilde \cC_{j+1}$ from $\tilde \cC_j$ by \eqref{refinej}.
%      \State $X_{j+1} := \PP_1(\tilde \cC_{j+1})$
%      \State $Z_{j+1} := \stab(X_{j+1})$
%      \State Update $\text{errbound}$ by
%			\begin{align*}
%			\text{errbound} = \sum_{Q \in \tilde \cC_{j}}\|P_{Q}A^*\hat r_a\|^2_{L_2(Q)}.
%			\end{align*}
%			\State Update $u_a: = \hat u_a$ and set $j \rightarrow j+1.$
%     \EndWhile
%  \end{algorithmic}
%\end{algorithm}
}

%%%%%%%%%%%%%%%%%%%%%%%%%%%%%%%%%
\section{Numerical Experiments}\label{sect:numerical}
%%%%%%%%%%%%%%%%%%%%%%%%%%%%%%%%%
We now provide some numerical experiments to illustrate
the performance of the previously introduced anisotropic adaptive
scheme for first order linear transport equations.
For the numerical tests, we consider two simple
transport equations{: One with a jump induced by the right hand side and one with a jump induced by the boundary conditions.} In both cases the
discontinuity curves will be seen to be accurately resolved by
relatively few degrees of freedom, showing the effectiveness of the scheme.
Moreover, we display numerically
estimated values for the stability constant $\delta$ given by
\[
\frac{\inf_{\phi \in Z_j}\norm{u_j -u_j^K - A^*\phi}_{L_2([0,1]^2)}}{\norm{u_j-u_j^K}_{L_2([0,1]^2)}},
\]
where $u_j = \text{argmin}_{v_j \in X_j}\norm{u-v_j}_2$.
Using that $\|\cdot\|_Y = \|A^*\cdot\|_X$ this is exactly the constant
$\delta$ in the definition \eqref{delta} of $\delta$-proximality,
however with $y$ (in \eqref{delta}) replaced by $u_j - u_j^K$.
This is of course only a lower bound for the true
$\delta$-proximality constant, but the given choice for $y$
is exactly the one used in the %major
proofs in \cite{DHSW}.

In the subsequent experiments, the number $K$ of iterations of the adaptive scheme
described in Section \ref{sect:ad} is for simplicity set to $K=10$.
One could as well employ an early termination of the inner iteration based on
a posteriori control of the lifted residuals $r^k_j$.
%*****************************************************************************************
%****************************************************************************************
\subsection{Linear Transport Equation with  homogeneous inflow data}
\begin{figure}[htb]
\begin{center}
\hspace*{\fill}
\subfigure[\label{fig:shear_approx1-a}]{\includegraphics[width=.4\textwidth]{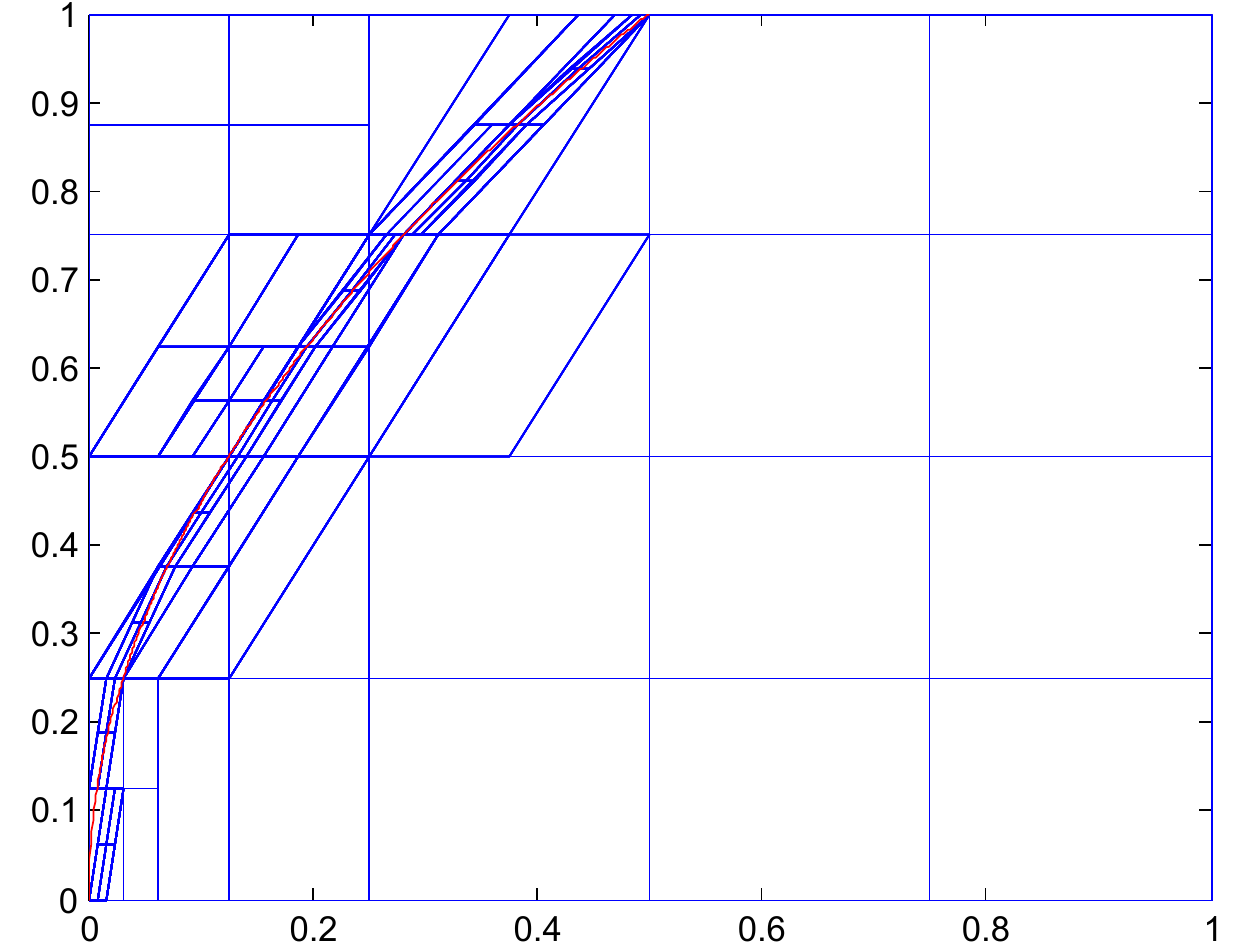}}
\hfill
\subfigure[\label{fig:shear_approx1-b}]{\includegraphics[width=.4\textwidth]{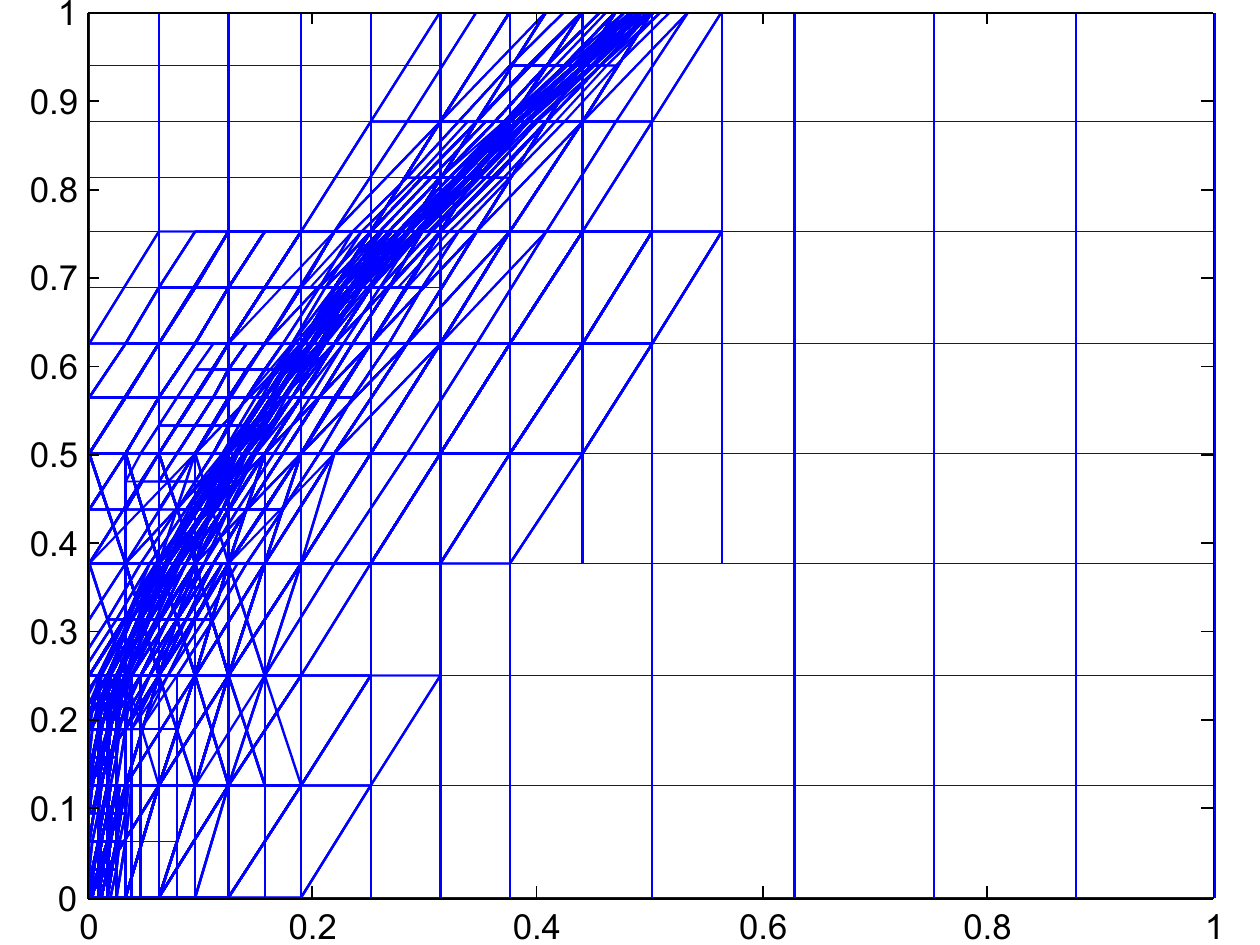}}
\hspace*{\fill} \\
\hspace*{\fill}
\subfigure[\label{fig:shear_approx1-c}]{\includegraphics[width=.4\textwidth]{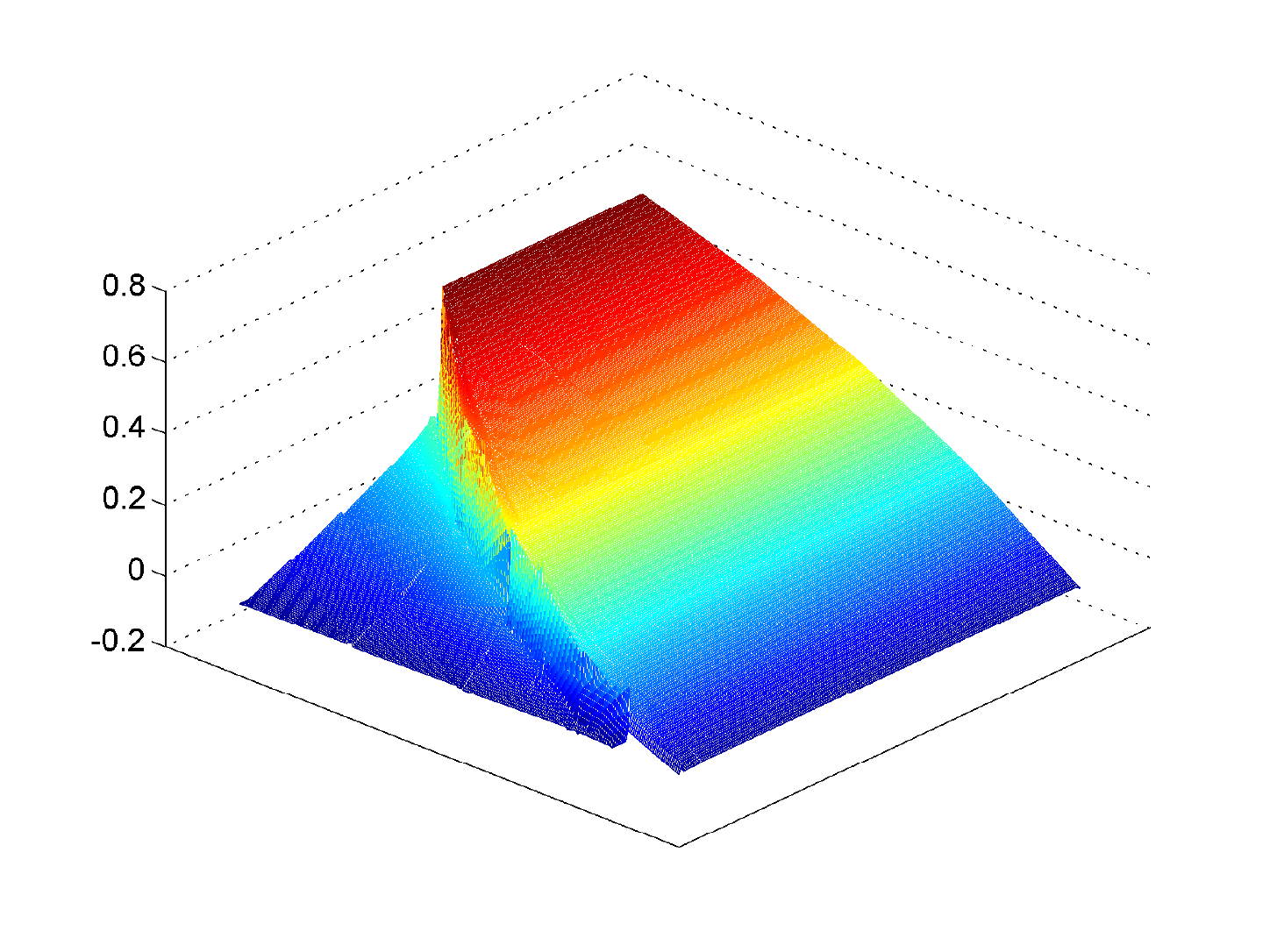}}
\hfill
\subfigure[\label{fig:shear_approx1-d}]{\includegraphics[width=.4\textwidth]{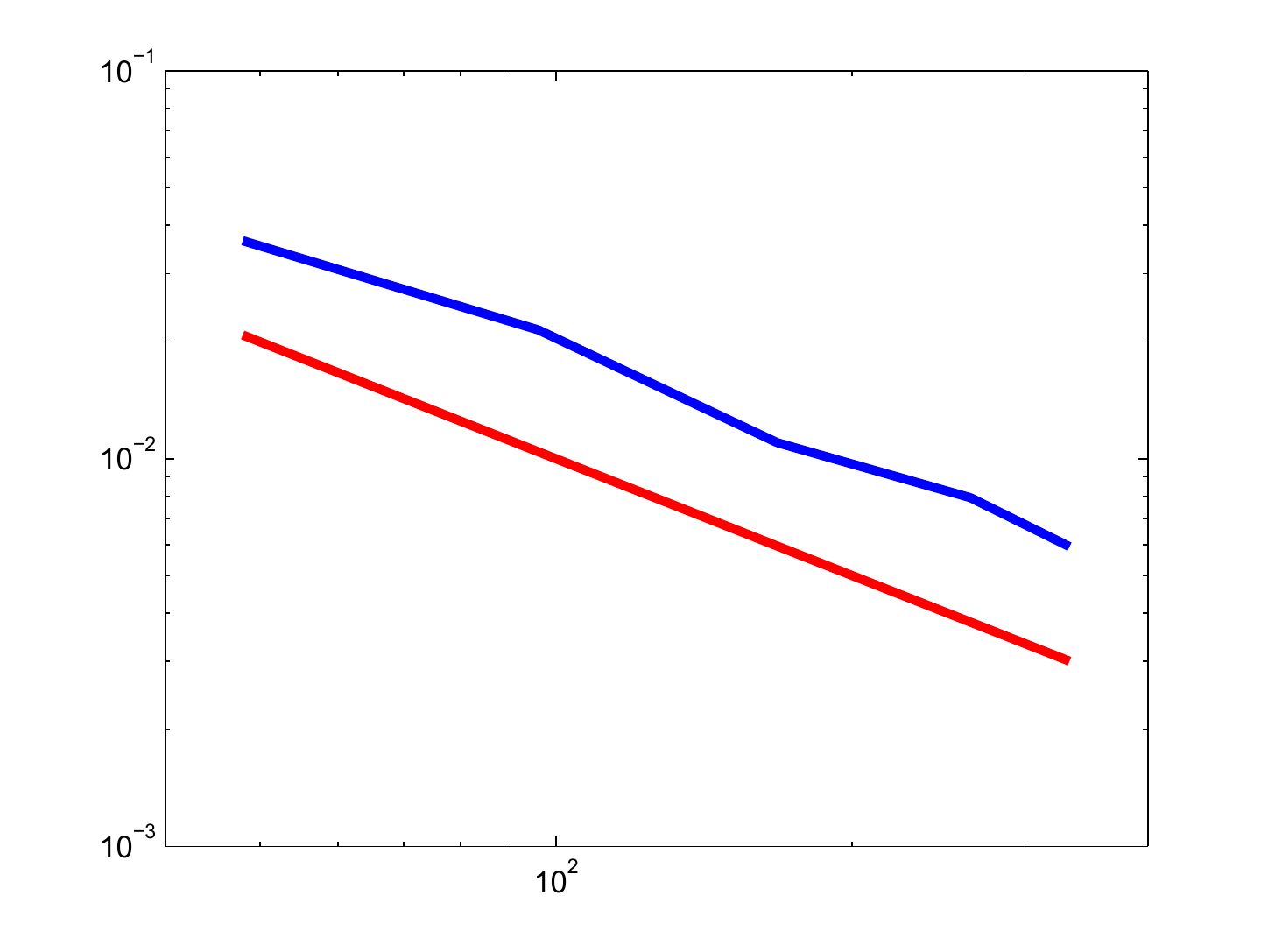}}
\hspace*{\fill}
\end{center}
\caption{\subref{fig:shear_approx1-a} Adaptive \cs{partition} for the trial space $X_5$.
\subref{fig:shear_approx1-b} Adaptive \cs{partition} for the test space $Z_5$.
\subref{fig:shear_approx1-c} Approximate solution (306 basis elements).
\subref{fig:shear_approx1-d} Approximation \cs{error in $L_2(D)$ versus
            the number of degrees of freedom} (blue) and the \cs{theoretical rate} (red).
}
\label{fig:shear_approx1}
\end{figure}

\begin{figure}[htb]
  \hfill
  \subfigure{\includegraphics[height=4cm]{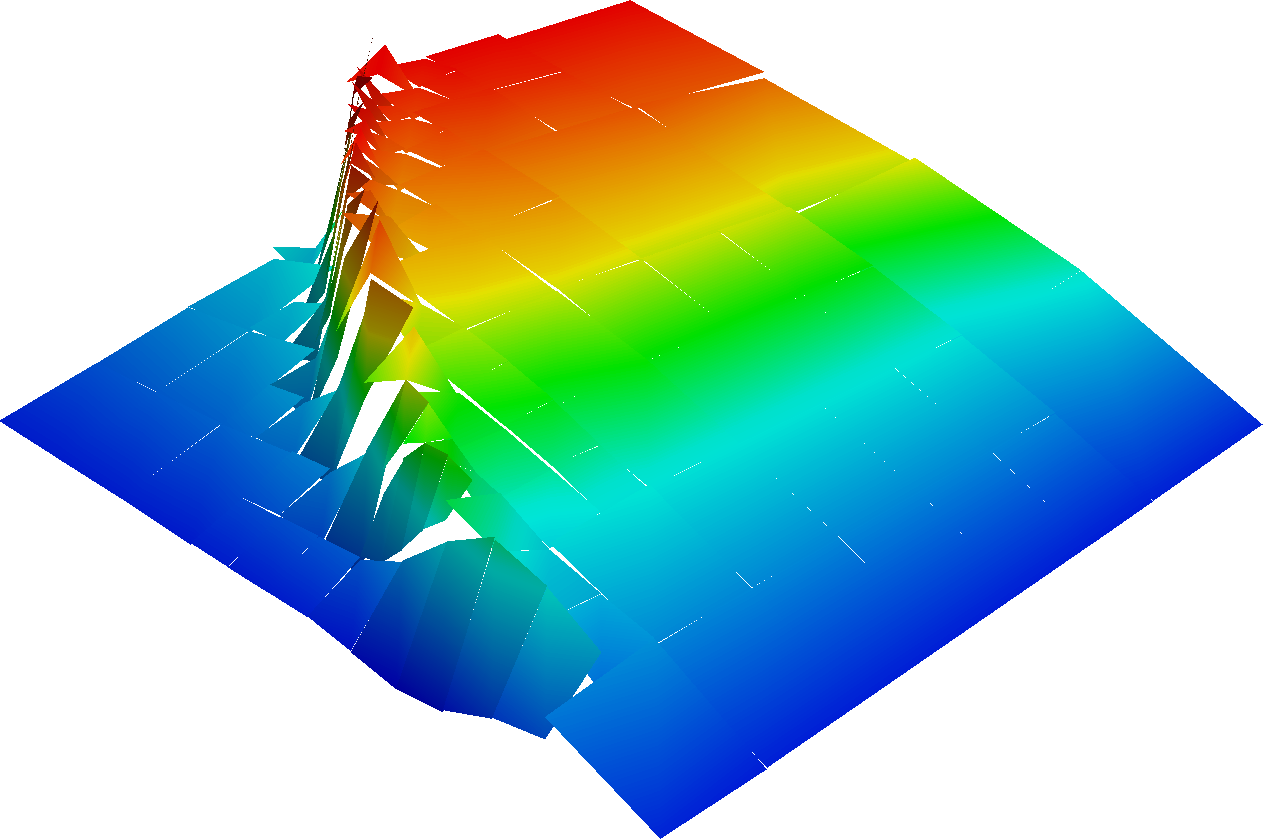}}
  \hfill
  \subfigure{\includegraphics[height=4cm]{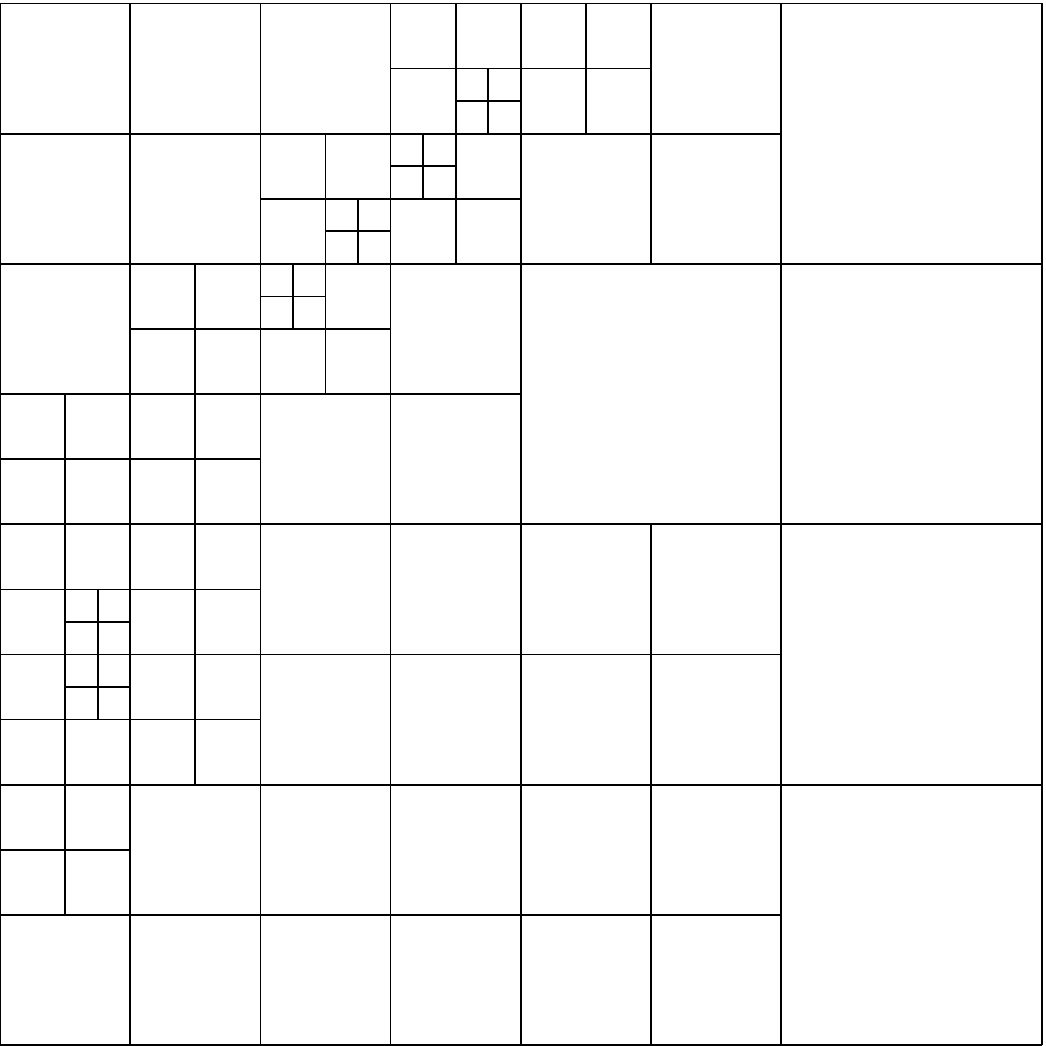}}
  \hspace*{\fill}
  \caption{Solution and \cs{partition} of the $8$th isotropically adaptive cycle for the test problem \eqref{eq:test-curve}. The trial space and test spaces consist of linear and bilinear finite elements, respectively. The \cs{partition for the} test space is obtaind from \cs{the corresponding
parition for} the trial space by refining each cell once.}
  \label{fig:curve-isotropic}
\end{figure}

We first consider a transport equation with zero inflow boundary condition,
whose solution exhibits a discontinuity along the
curvilinear layer given by $x_1 = \frac12 x_2^2$.
More precisely, we consider the equation given by
\begin{equation}
Au =
\left(
  \begin{array}{c}
    x_2 \\
    1 \\
  \end{array}
\right)
 \cdot u + u = f, \quad f = \chi_{\{x_1 > x_2^2/2\}} + 1/2 \cdot \chi_{\{x_1 \leq x_2^2/2\}}.
 \label{eq:test-curve}
\end{equation}
Figure \ref{fig:shear_approx1}
\subref{fig:shear_approx1-a},
\subref{fig:shear_approx1-b}
show the adaptive \cs{partitions} associated with the trial space $X_5$ and
the test space $Z_5$ clearly demonstrating their highly anisotropic
structure reflected by   the refinements in the neighborhood of the
discontinuity curve.
Figure \ref{fig:shear_approx1-c} illustrates the  {anisotropic} approximation given by
306 basis elements.
For a comparison, Figure \ref{fig:curve-isotropic} shows the corresponding results with isotropic refinements only.
We emphasize  that in our scheme, Gibbs like phenomena and spurious oscillations across the jump are almost completely absent and in fact much less pronounced than observed for isotropic discretizations.
%\cs{Agree to this; Gerritt Welper has the files (?)}
Figure \ref{fig:shear_approx1-d} shows the optimality of {the anisotropic} scheme.
The numerical stability is confirmed by Table \ref{fig:shear_delta1}.

\begin{table}[htb]

\begin{filecontents*}{figs/table1.csv}
n (anisotropic), delta (anisotropic), error (anisotropic), rate (anisotropic), n (isotropic), delta (isotropic), error (isotropic), rate (isotropic)
   48 , 0.298138 , 0.036472   ,  ,52   , 0.246551 , 0.0893218 , 0.371732  
   99 , 0.442948 , 0.021484   ,  0.7311,100  , 0.260341 , 0.0753951 , 0.259209
  138 , 0.352767 ,  0.013948  ,  1.3132,160  , 0.286508 , 0.0652673 , 0.306914
  177 , 0.322156 , 0.010937   ,  0.9768,256  , 0.332578 , 0.0549427 , 0.366383 
  237 , 0.316545 , 0.008348   ,  0.9130,388  , 0.376594 , 0.0487029 , 0.289906
  306 , 0.307965 ,  0.006152  ,  1.2100,556  , 0.439148 , 0.0440333 , 0.280168
\end{filecontents*}
\pgfplotstabletypeset[
    col sep=comma,
    every head row/.style={
        before row={%
            \hline
            \multicolumn{4}{|c|}{Anisotropic} &
            \multicolumn{4}{c|}{Isotropic} \\
            \hline
        },
        after row={\hline \hline}
    },
    columns/{n (anisotropic)}/.style={column name=$n$, column type/.add={|}{}},
    columns/{delta (anisotropic)}/.style={column name=Estimated $\delta$},
    columns/{error (anisotropic)}/.style={column name=$L_2$ error},
    columns/{rate (anisotropic)}/.style={column name=Rate, column type/.add={}{|}},
    %
    %column type/.add={|}{},
    %
    columns/{n (isotropic)}/.style={column name=$n$},
    columns/{delta (isotropic)}/.style={column name=Estimated $\delta$},
    columns/{error (isotropic)}/.style={column name=$L_2$ error},
    columns/{rate (isotropic)}/.style={column name=Rate, column type/.add={}{|}},
    every last row/.style={after row=\hline}
]{figs/table1.csv}

\caption{Numerical estimates for the stability constant $\delta$,
         $L_2$ approximation error and \cs{estimated} convergence rates.
%\gw{(The rates for the anisotropic case are missing).}
}
\label{fig:shear_delta1}
\end{table}

%\begin{table}[htb]
%\pgfplotstabletypeset[
%    col sep=comma,
%    %
%    every head row/.style={
%        before row={%
%            \hline
%            \multicolumn{3}{|c|}{Anisotropic} &
%            \multicolumn{4}{c|}{Isotropic} \\
%            \hline
%        },
%        after row={\hline \hline}
%    },
%    %
%    columns/{n (anisotropic)}/.style={column name=$n$, column type/.add={|}{}},
%    columns/{delta (anisotropic)}/.style={column name=Estimated $\delta$},
%    columns/{error (anisotropic)}/.style={column name=$L_2$ error, column type/.add={}{|}},
%    %
%    %column type/.add={|}{},
%    %
%    columns/{n (isotropic)}/.style={column name=$n$},
%    columns/{delta (isotropic)}/.style={column name=Estimated $\delta$},
%    columns/{error (isotropic)}/.style={column name=$L_2$ error},
%    columns/{rate (isotropic)}/.style={column name=Rate, column type/.add={}{|}},
%    %
%    every last row/.style={after row=\hline}
%]{figs/table1.csv}
%\caption{Numerical estimates for the stability constant $\delta$, $L_2$ approximation error and convergence rates.
%%\gw{(The rates for the anisotropic case are missing).}
%}
%\label{fig:shear_delta1}
%\end{table}

% Style to select only points from #1 to #2 (inclusive)
\pgfplotsset{select coords between index/.style 2 args={
    x filter/.code={
        \ifnum\coordindex<#1\def\pgfmathresult{}\fi
        \ifnum\coordindex>#2\def\pgfmathresult{}\fi
    }
}}

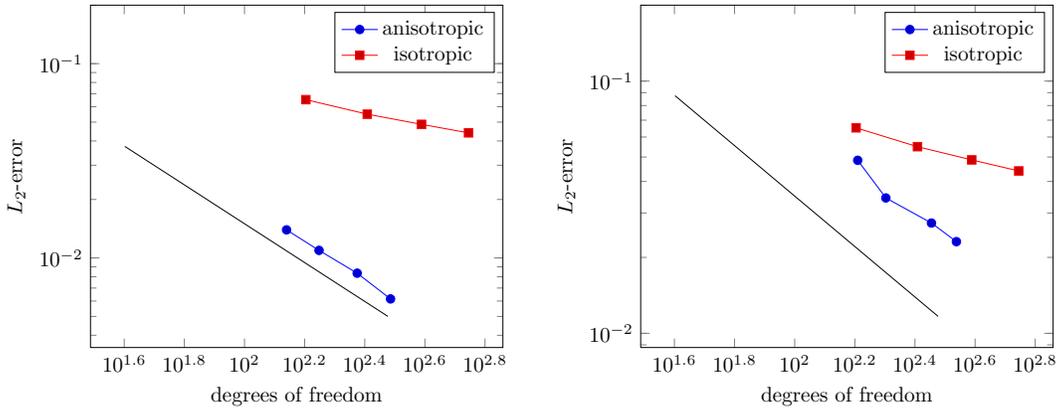
\begin{figure}

\hfill
\begin{tikzpicture}[scale=0.8]
    \begin{loglogaxis}[
        xlabel=degrees of freedom,
        ylabel=$L_2$-error,
        ymax = 0.2]
    \addplot table [x=n (anisotropic), y=error (anisotropic), col sep=comma, select coords between index={2}{7}] {figs/table1.csv};
    \addlegendentry{anisotropic}

    \addplot table [x=n (isotropic), y=error (isotropic), col sep=comma, select coords between index={2}{7}] {figs/table1.csv};
    \addlegendentry{isotropic}

    \addplot [domain=40:300] {1.5*1/x} coordinate[pos=0.4] (A) coordinate[pos=0.7] (B);
    %\addlegendentry{rate $n^{-1}$}

    \draw (A) |- (B);

    \end{loglogaxis}
\end{tikzpicture}
\hfill
\begin{tikzpicture}[scale=0.8]
    \begin{loglogaxis}[
        xlabel=degrees of freedom,
        ylabel=$L_2$-error,
        ymax = 0.2]
    \addplot table [x=n (anisotropic), y=error (anisotropic), col sep=comma, select coords between index={2}{7}] {figs/table2.csv};
    \addlegendentry{anisotropic}

    \addplot table [x=n (isotropic), y=error (isotropic), col sep=comma, select coords between index={2}{7}] {figs/table2.csv};
    \addlegendentry{isotropic}

    \addplot [domain=40:300] {3.5*1/x} coordinate[pos=0.4] (A) coordinate[pos=0.7] (B);
    %\addlegendentry{rate $n^{-1}$}

    \draw (A) |- (B);

    \end{loglogaxis}

\end{tikzpicture}
\hspace*{\fill}

\caption{$L_2$ errors of Examples \eqref{eq:test-curve} and \eqref{eq:test-boundary} and a line indicating the rate $n^{-1}$, respectively.}

\end{figure}

%
%************************************************************************************************************
%************************************************************************************************************
\subsection{Linear Transport Equation with  inhomogeneous inflow data}
\label{sec:LinTrInhDat}
\begin{figure}[h]
\begin{center}
\hspace*{\fill}
\subfigure[\label{fig:shear_approx3-a}]{\includegraphics[width=.4\textwidth]{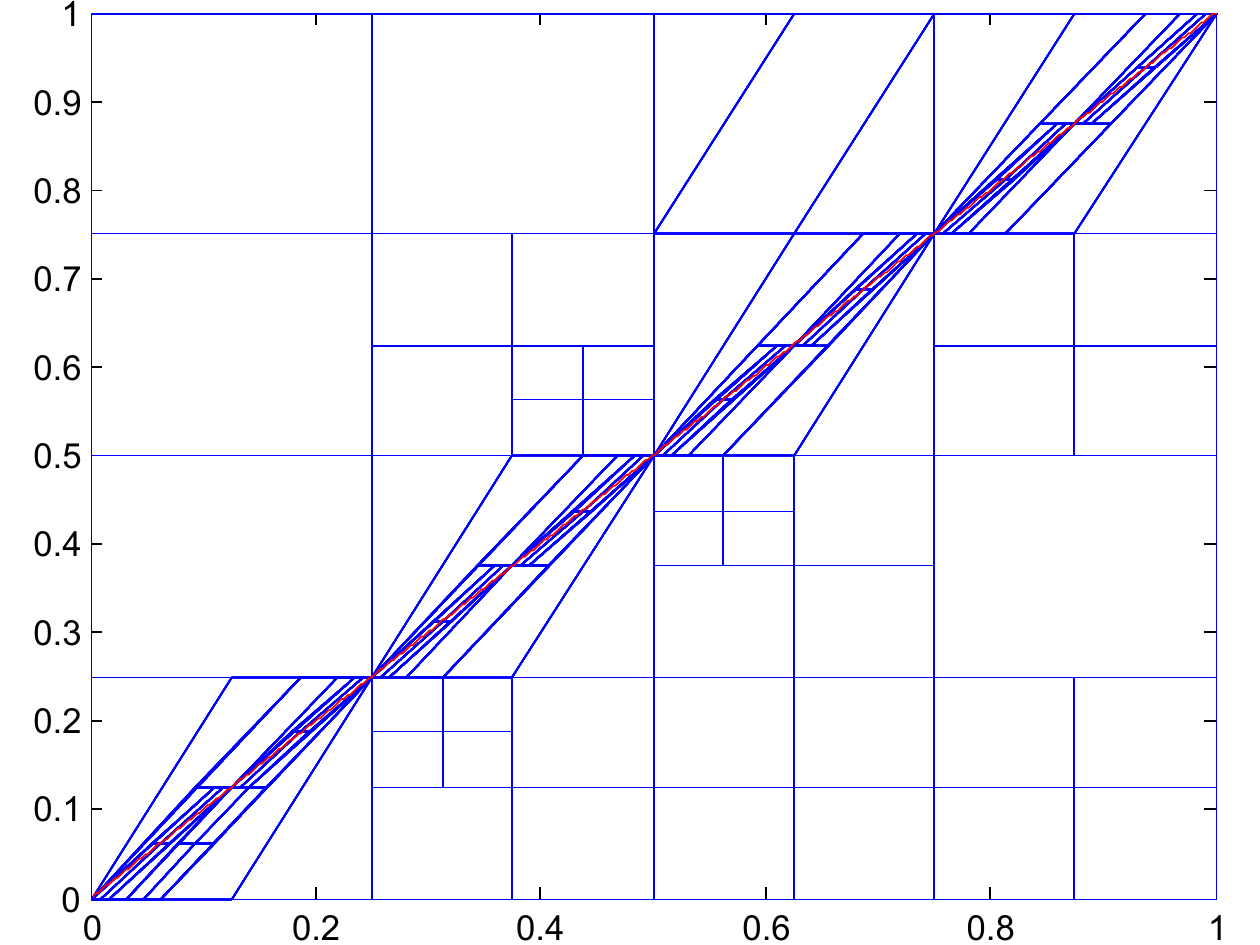}}
\hfill
\subfigure[\label{fig:shear_approx3-b}]{\includegraphics[width=.4\textwidth]{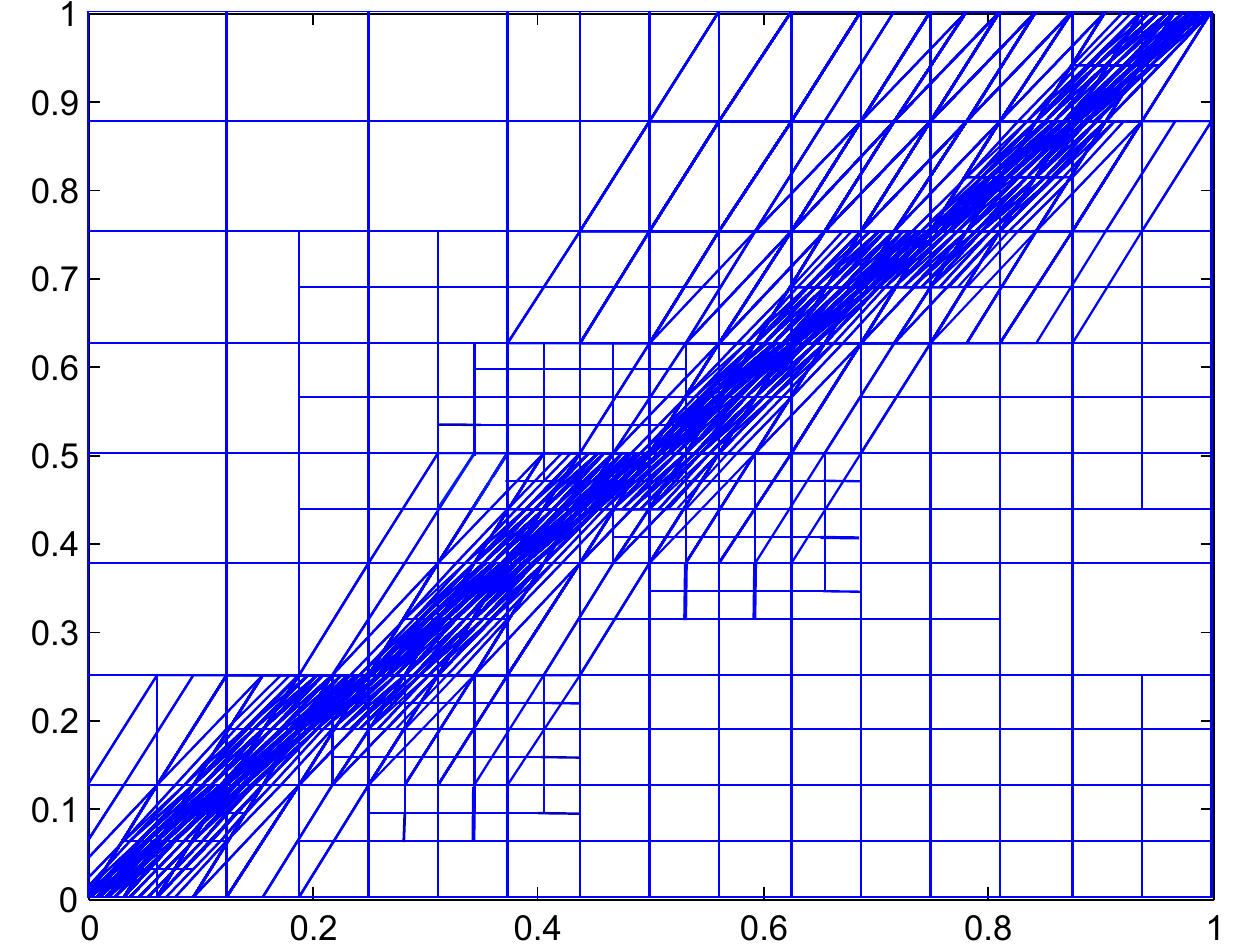}}
\hspace*{\fill} \\
\hspace*{\fill}
\subfigure[\label{fig:shear_approx3-c}]{\includegraphics[width=.4\textwidth]{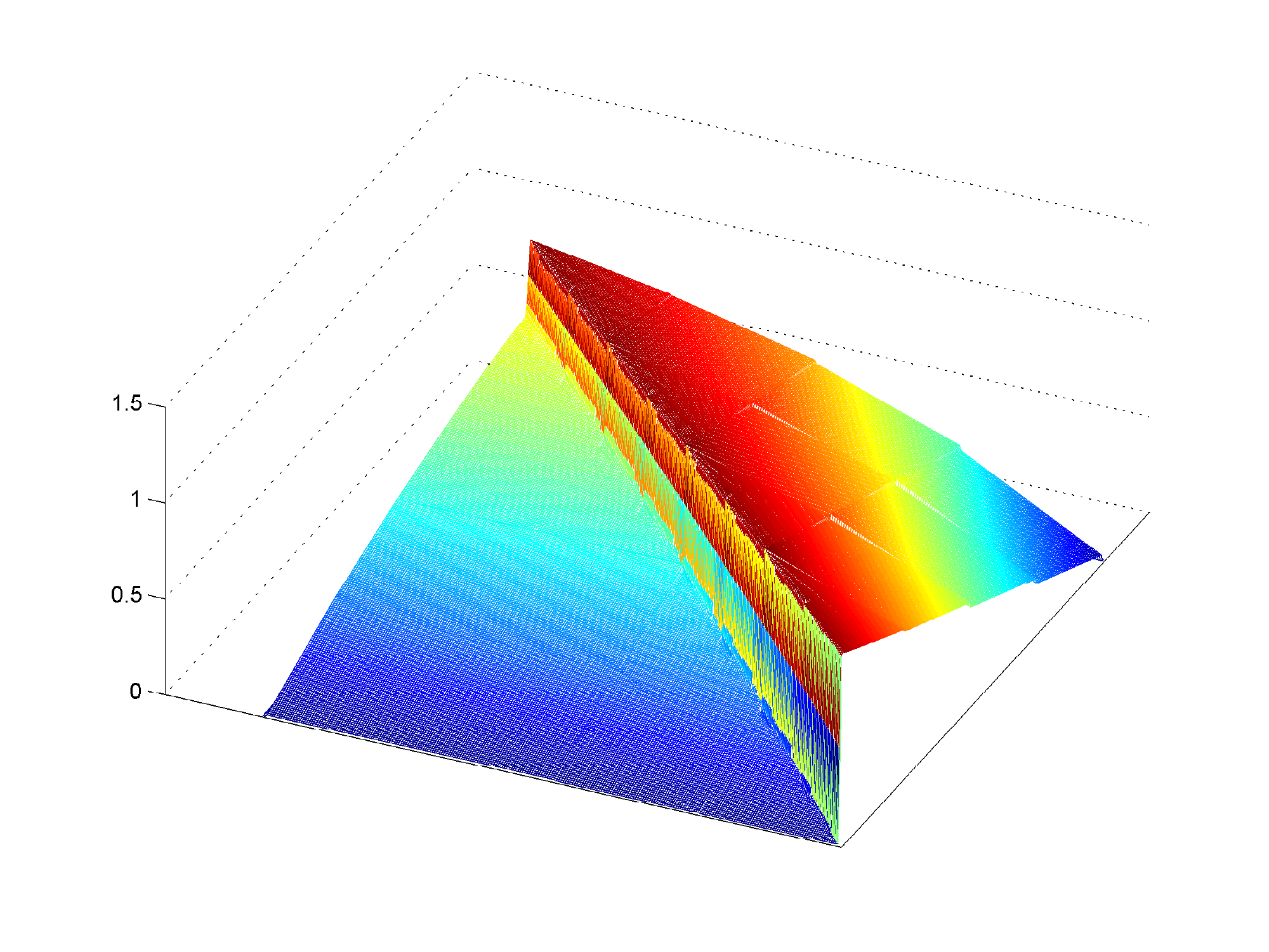}}
\hfill
\subfigure[\label{fig:shear_approx3-d}]{\includegraphics[width=.4\textwidth]{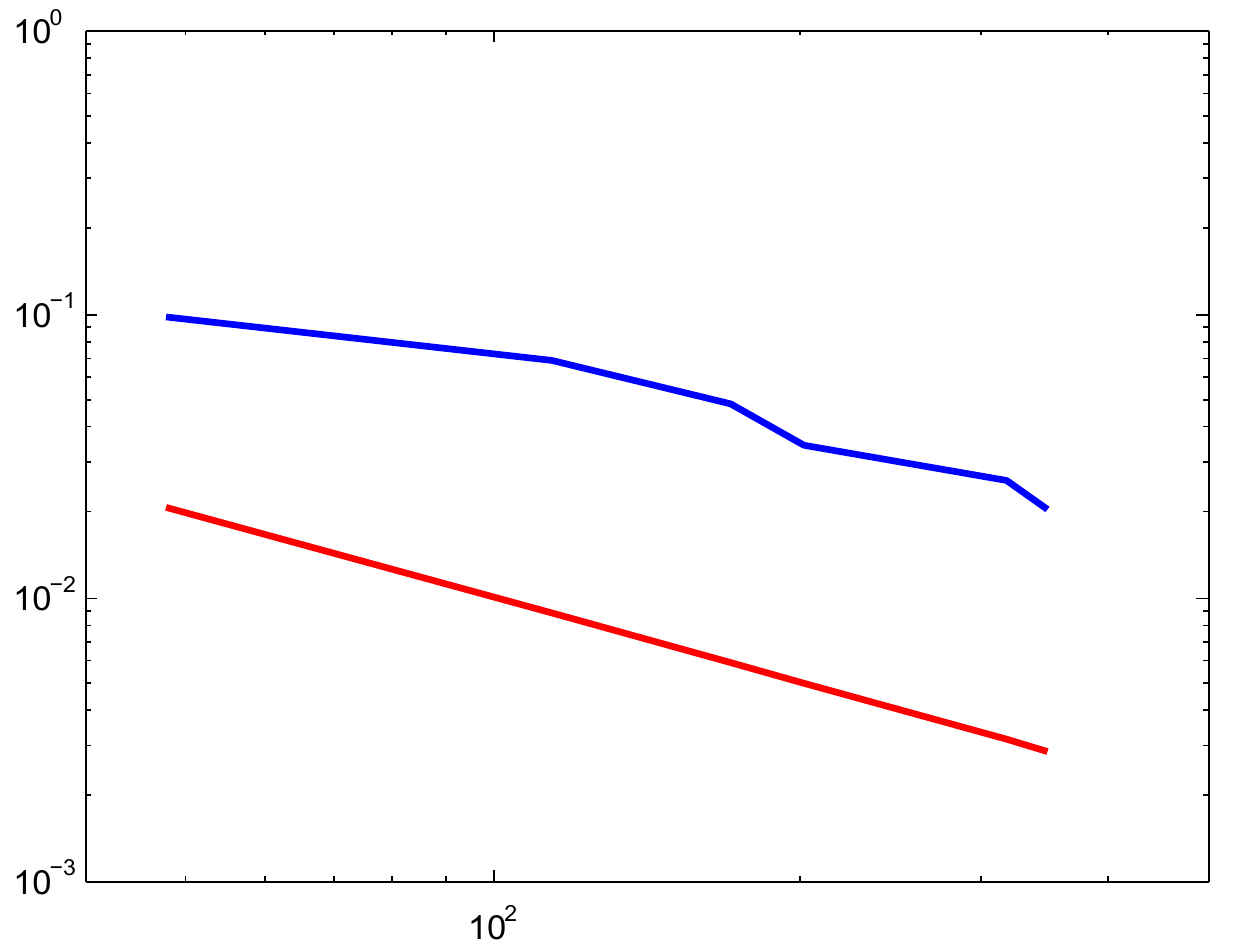}}
\hspace*{\fill}
\end{center}
\caption{\subref{fig:shear_approx3-a} Adaptive \cs{partition}
for the trial space $X_5$.
\subref{fig:shear_approx3-b} Adaptive \cs{partition} for the test space $Z_5$.
\subref{fig:shear_approx3-c} Approximate solution (345 basis elements).
\subref{fig:shear_approx3-d}
 Approximation \cs{error in $L_2(D)$ versus
            the number of degrees of freedom} (blue) and the \cs{theoretical rate} (red).
}
\label{fig:shear_approx3}
\end{figure}

\begin{figure}[htb]
  \hfill
  \subfigure{\includegraphics[height=4cm]{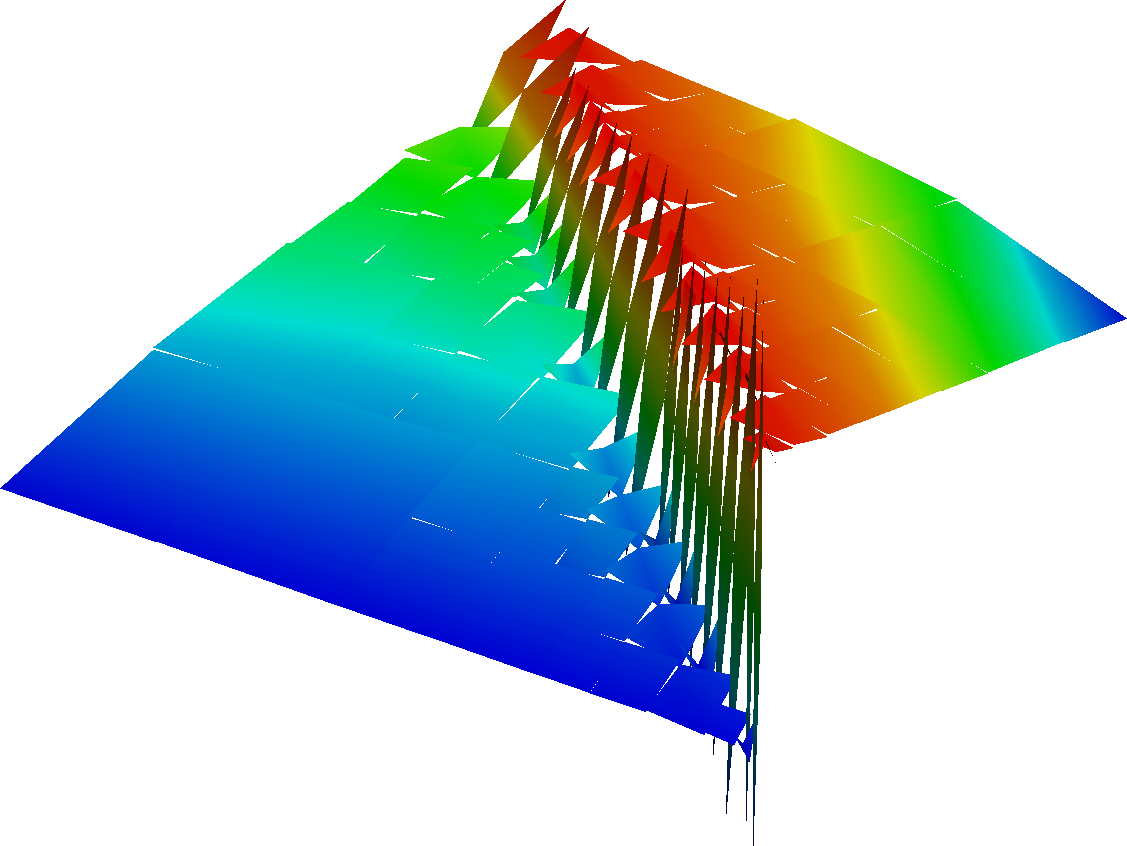}}
  \hfill
  \subfigure{\includegraphics[height=4cm]{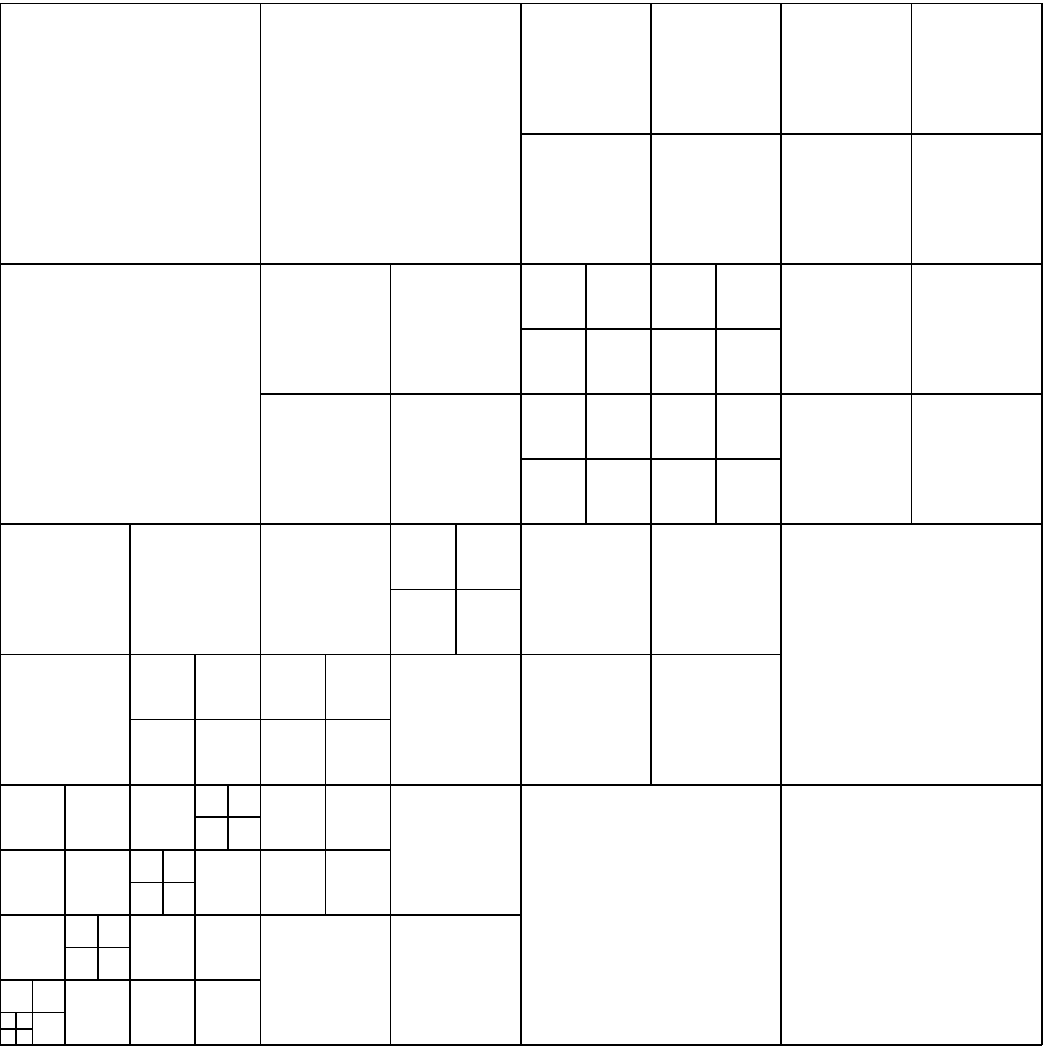}}
  \hspace*{\fill}
  \caption{Solution and \cs{partition} of the $6$th isotropically adaptive cycle for the test problem \eqref{eq:test-boundary}. The trial space and test spaces consist of linear and bilinear finite elements, respectively. The \cs{partition of the} test space is obtaind from the
\cs{the corresponding partition of the} trial space by refining each cell once.}
  \label{fig:boundary-isotropic}
\end{figure}

Second, as a classical benchmark,
we consider a transport equation with non-zero boundary condition, whose solution exhibits a discontinuity
along the shear layer given by $x_1 = x_2$. More precisely, we consider the equation given by
\begin{equation}
Au =
\left(
  \begin{array}{c}
    1 \\
    1 \\
  \end{array}
\right)
 \cdot u + u = f, \quad f = 1/2
\label{eq:test-boundary}
\end{equation}
with boundary conditions
\[
g(x_1,x_2) = 1-x_1 \,\, \text{on} \,\, \{ (x_1,0) \in \Gamma_{-} : 0 < x_1 < 1 \}
\]
and
\[
g(x_1,x_2) = 0 \,\, \text{on} \,\, \{ (0,x_2) \in \Gamma_{-} : 0 < x_2 < 1 \}.
\]
 Figures \ref{fig:shear_approx3} and \ref{fig:boundary-isotropic} show the numerical results for the anisotropic and isotropic cases, respectively. As in the previous case the anisotropic solutions show negligible spurious oscillations and are more stable than  the isotropic ones.
 {The $L_2$ approximation error and $\delta$ estimates can be found in Table} \ref{fig:shear_delta3}.

%\begin{figure}
%\begin{center}
%\begin{tabular}{|c|c|c||c|c|c|c|}
%  \hline
%  \multicolumn{3}{|c||}{Anisotropic} &
%  \multicolumn{4}{|c|}{Isotropic} \\
%  \hline
%  $n$ & Estimated $\delta$ & $\norm{u_j^K - u}_{L_2([0,1]^2)}$
%  & $n$ & Estimated $\delta$ & $\norm{u_j^K - u}_{L_2([0,1]^2)}$ & Estimated rate \\
%  \hline
%  \hline
%   48 & 0.141487 & 0.097105  &  16   & 0.138998 &  0.143281  &          \\
%  108 & 0.188977 & 0.068675  &  28   & 0.150498 &  0.112434  & 0.433229 \\
%  162 & 0.217835 & 0.048527  &  52   & 0.246551 &  0.0893218 & 0.371732 \\
%  201 & 0.184869 & 0.034419  &  100  & 0.260341 &  0.0753951 & 0.259209 \\
%  285 & 0.190835 & 0.027366  &  160  & 0.286508 &  0.0652673 & 0.306914 \\
%  345 & 0.183383 & 0.023095  &  256  & 0.332578 &  0.0549427 & 0.366383 \\
%      &          &           &  388  & 0.376594 &  0.0487029 & 0.289906 \\
%      &          &           &  556  & 0.439148 &  0.0440333 & 0.280168 \\
%      &          &           &  808  & 0.440738 &  0.0414834 & 0.159582 \\
%      &          &           &  1060 & 0.454505 &  0.0385443 & 0.270704 \\
%      &          &           &  1360 & 0.478352 &  0.0354129 & 0.339994 \\
%  \hline
%\end{tabular}
%\end{center}
%\caption{Numerical estimates for the stability constant $\delta$ and the $L_2$ approximation error. }
%\label{fig:shear_delta3}
%\end{figure}

\begin{table}[htb]
\begin{filecontents*}{figs/table2.csv}
n (anisotropic), delta (anisotropic), error (anisotropic), rate (anisotropic), n (isotropic), delta (isotropic), error (isotropic), rate (isotropic)
   48 , 0.141487 , 0.097105 , ,  52   , 0.246551 ,  0.0893218 , 0.371732     
  108 , 0.188977 , 0.068675 , 0.4267,  100  , 0.260341 ,  0.0753951 , 0.259209    
  162 , 0.217835 , 0.048527 , 0.8587,  160  , 0.286508 ,  0.0652673 , 0.306914    
  201 , 0.184869 , 0.034419  , 1.5925,  256  , 0.332578 ,  0.0549427 , 0.366383    
  285 , 0.190835 , 0.027366  , 0.6516,  388  , 0.376594 ,  0.0487029 , 0.289906 
  345 , 0.183383 , 0.023095  , 0.8935,  556  , 0.439148 ,  0.0440333 , 0.280168 
  
  
\end{filecontents*}
\pgfplotstabletypeset[
    col sep=comma,
    every head row/.style={
        before row={%
            \hline
            \multicolumn{4}{|c|}{Anisotropic} &
            \multicolumn{4}{c|}{Isotropic} \\
            \hline
        },
        after row={\hline \hline}
    },
    columns/{n (anisotropic)}/.style={column name=$n$, column type/.add={|}{}},
    columns/{delta (anisotropic)}/.style={column name=Estimated $\delta$},
    columns/{error (anisotropic)}/.style={column name=$L_2$ error},
    columns/{rate (anisotropic)}/.style={column name=Rate, column type/.add={}{|}},
    %
    %column type/.add={|}{},
    %
    columns/{n (isotropic)}/.style={column name=$n$},
    columns/{delta (isotropic)}/.style={column name=Estimated $\delta$},
    columns/{error (isotropic)}/.style={column name=$L_2$ error},
    columns/{rate (isotropic)}/.style={column name=Rate, column type/.add={}{|}},
    every last row/.style={after row=\hline}
]{figs/table2.csv}
\caption{Numerical estimates for the stability constant $\delta$, $L_2$ approximation error and convergence rates.
}
\label{fig:shear_delta3}
\end{table}

%\begin{table}[htb]
%\pgfplotstabletypeset[
%    col sep=comma,
%    %
%    every head row/.style={
%        before row={%
%            \hline
%            \multicolumn{3}{|c|}{Anisotropic} &
%            \multicolumn{4}{c|}{Isotropic} \\
%            \hline
%        },
%        after row={\hline \hline}
%    },
%    %
%    columns/{n (anisotropic)}/.style={column name=$n$, column type/.add={|}{}},
%    columns/{delta (anisotropic)}/.style={column name=Estimated $\delta$},
%    columns/{error (anisotropic)}/.style={column name=$L_2$ error, column type/.add={}{|}},
%    %
%    %column type/.add={|}{},
%    %
%    columns/{n (isotropic)}/.style={column name=$n$},
%    columns/{delta (isotropic)}/.style={column name=Estimated $\delta$},
%    columns/{error (isotropic)}/.style={column name=$L_2$ error},
%    columns/{rate (isotropic)}/.style={column name=Rate, column type/.add={}{|}},
%    %
%    every last row/.style={after row=\hline}
%]{figs/table2.csv}
%\caption{Numerical estimates for the stability constant $\delta$, $L_2$ approximation error and convergence rates.
%%\gw{(The rates for the anisotropic case are missing).}
%}
%\label{fig:shear_delta3}
%\end{table}

In summary, the experiments confirm that the employed anisotropic adaptive refinements essentially recover the
optimal rates
relating the achieved accuracy to the used number of degrees of freedom.
Overall this results in a better rate than the one
achieved by isotropic refinements. Perhaps more importantly, the anisotropy appears to greatly promote the stability of the scheme, i.e.,
the projections form the spaces $Z_j$ (although not much larger in dimension than $X_j$) come sufficiently close to the optimal
(computationally infeasible) test spaces to give rise to very well conditioned variational problems.
%\todo{[Omit next sentence, unless we show isotropic refinements]}
Moreover,
\cs{pointwise over- and undershoots of the approximations analogous
    to Gibbs' phenomena near jump discontinuities when approximating in $L_2(D)$
in these examples}
is much less pronounced than with isotropic refinements.

%\gw{(At the beginning of Section \ref{sect:aniso} a list of example problems is given.
%     In this section we only test two of these. Maybe we should keep that consistent.)}
%{In principle, I agree. It would be interesting to have a right hand side with diagonal jump while a jump
%in the boundary conditions is convected along the above curve, or vice versa. This should give an intersection of discontinuity curves.
%Wang-Q: would that be possible?}

%%%%%%%%%%%%%%%%%%%%%%
\section{The Parametric Case}\label{sect:parametric}
%%%%%%%%%%%%%%%%%%%%%%%%
\newcommand{\Lip}{{\rm Lip}}
\newcommand{\dd}{{\mathbf d}}
%
%{\bf Revised version for integral functionals; Preliminary (to be discussed)}

We now consider parametric transport problems of the type
\eqref{eq:rads} and adhere to the notation and terminology introduced
in Section \ref{sect2.3}. We recall \eqref{eq:Gammapms}, i.e.
\css{the in- and outflow boundaries}
$$
{\partial D}_\pm(\bs)
:=
\{ x\in \partial D : \pm \vb(\bs) \cdot \vec{n}(x) > 0 \},\;\; \bs\in \cS,
$$
now depend on the transport direction $\bs\in \cS$.
Similar to \cite{KGCS2010a,KGCS2010b}, we will now investigate the possibility
of sparse tensor product approximations of the parametric transport problem
\eqref{eq:rads}--\eqref{eq:Gammapms}, using the directionally adaptive multi-level
approximation rates in $L_2(D)$ which were obtained in Theorem \ref{ShearCart}.

We outline in this section some circumstances under which the preceding anisotropic
approximations lead in connection with {\em sparse tensor product} concepts
to the efficient computation of simple quantities of interest
of the parametric fields $u(\cdot \bs)$.
\subsection{Sparse Tensor Approximation Rate Estimates}
\label{sec:ApprRateEst}
We assume given a continuous, linear
functional $G(\cdot) \in (L^2(D))^*$,
independent of the transport direction $\bs$.
The goal of computation is now to
approximate an integral quantity over all directions $\bs$,
such as the total emission:
\beq \label{eq:IGu}
I(G(u)) := \int_\cS G(u(\cdot,\bs)) d\bs \;.
\eeq
In order to obtain superior rates of sparse tensor product approximations,
we require that the solution $u(\bs,x)$ of \eqref{eq:rads}--\eqref{eq:Gammapms}
belongs to the cartoon class $\cC$ uniformly over all directions $\bs$,
with a certain amount of regularity with respect to $\bs$.
We work under the following assumptions.
\begin{assumption}\label{As:Cartoon}

\begin{enumerate}
\item
Condition \eqref{coerccond} holds uniformly with respect to the direction $\bs\in \cS$,
\item
%For each parameter instance $\bs\in \cS$, there exist
%tuples $\Omega(\bs),\kappa(\bs), L(\bs),M(\bs)$ as in \eqref{cartoon}
%such that problem \eqref{eq:rads} admits a unique
%solution $u(\bs,\cdot)\in {\cC}(\Omega(\bs),\kappa(\bs), L(\bs),M(\bs))$,
%[WD: I don't understand this condition. The cartoon class can be viewed as a {\em collection}
%of $\Omega's$ which therefore should not appear as characterizing parameters. We don't index e.g. an
%$L_2$-ball of functions by the functions $f$ themselves but by the radius, say. The cartoon class is characterized
%by the bounding parameters $\kappa,L,M,\omega$, limiting admissible curvatures, curve length, second oder
%dervatives, and curve separation. For this setting Theorem \ref{ShearCart} has been proved.
%So, shouldn't this assumption read:\\
There exist parameters $\kappa,L,M,\omega$ (as in Theorem \ref{ShearCart}) such that
problem \eqref{eq:rads} admits for each $\bs\in \cS$  a unique solution $u(\bs,\cdot)\in \cC(\kappa,L,M,\omega)$
and
%
%\item
%the constants $\kappa(\bs), L(\bs),M(\bs)$ are uniformly bounded
%w.r. to $\bs\in \cS$ by $\kappa,L,M$, and the cartoon sets $\Omega(\bs)$
%admit the representation \eqref{eq:boundary},
%for all $\bs \in \cS$, with \eqref{Eprop} valid uniformly w.r. to $\bs$,
%\dw{[WD: why is this condition needed? The approximation properties were established for
%the general cartoon class not only for the horizon model. The solutions we consider happen to
%fall into the horizon model but not in the form  \eqref{Eprop} but with $x_2,x_1$ interchanged.
%Again, I don't think we need this and should skip it.}
%
\item
$G\circ u \in C^{\cs{0,\alpha}}(\cS)$ for some $0<\alpha \leq 1$.
We remark that this requirement
necessitates regularity with respect to $\bs\in \cS$ \emph{and}
in the class $\cC$ of admissible cartoons.
\end{enumerate}
\end{assumption}
Note that for (uniformly in $\bs$) $C^1$-convection fields, piecewise smooth boundary conditions
and data $f$ in a cartoon class with $C^1$-pieces separated by a $C^2$-curve,
the solutions belong to a cartoon class with parameters depending on the data
(right hand side, boundary data, convection field)
and thus satisfies Assumption \ref{As:Cartoon}(2).
We assume we have at hand a family $\{ I_\ell \}_{\ell \geq 0}$
of quadrature approximations to the direction-integral
$I$ in \eqref{eq:IGu} which,
for integrand functions $\varphi \in C^{0,\alpha}(\cS)$,
for some $0<\alpha \leq 1$, converges with $O(h_\ell^\alpha)$
\beq\label{eq:I_ellErr}
\left| I(\varphi) - I_\ell(\varphi) \right|
\lesssim
h^\alpha_\ell
\cs{| \varphi |_{C_{h_\ell}^{0,\alpha}(\cS)}}, \quad \ell \to \infty,
\;.
\eeq
where
$| \varphi |_{C_{h}^{0,\alpha}(\cS)}
:=
\sup_{\bs, \bs':|\bs -\bs'|\ge h}|\bs-\bs'|^{-\alpha}|\varphi(\bs)-\varphi(\bs')|$
is a ``conditional'' Lipschitz \cs{semi-}norm which is for $h>0$
weaker than $\| \varphi \|_{C^{0,\alpha}(\cS)}$.
The composed Midpoint- and Trapezoidal rules satisfy \eqref{eq:I_ellErr}.

For $j\geq 0$, and for $u\in \cC\subset L^2(D)$,
we denote by \cs{$B^x_ju$} the anisotropic
(and hence nonlinear) approximations constructed in Theorem \ref{ShearCart}
%\dw{[WD: we don't construct a shearlet-approximation.]}.
%We  recall from the proof of Theorem \ref{ShearCart} that, for each given
%$\bs_i \in \cS$,  that under the above assumptions
%$u(\bs_i,\cdot) \in \cC(\kappa,L,M,\omega)$, there exists a sequence
%of %shearlet-based adaptive,
%anisotropic approximation operators
%%\dw{[WD: I don't think we have constructed a well defined operator
%as a concise rule that applies
%to any element in the cartoon class.]}
%$B_J^{x}$ (depending among others on the direction $\bs_i$)
%and a constant $C>0$
%(depending only on the parameters of the class $\cC$)
such that \cs{there exists a constant $C>0$ such that for all $J\geq 1$}
\beq\label{eq:DsShear}
\|u- B_J u\|_{L_\infty(\cS;L_2(D))}\le C J^{1/2} 2^{-J} \;.
%\| u - B_J^{x}u \|_{C^{0,\alpha}(\cS;L_2(D))} \leq C J^{1/2} 2^{-J}\;,
\eeq
\cs{We recall that} the constant $C$ depends only on the cartoon parameters.
%
%\dw{[WD: I don't see why this is true. I think we have proved under the above assumptions only
%$$
%\|u- B_J^x u\|_{L_\infty(\cS;L_2(D))}\le C J^{1/2} 2^{-J}.
%$$}
%with constant $C>0$ independent of $J$
%and such that the number of degrees of freedom $N_J^{x}$
%which are activated in directional adaptive approximation
%is not larger than $2^J$.

For a favorable sparse approximation error bound,
we require the \emph{stronger condition} that there exists $C>0$
independent of $j$ such that
\beq
\label{eq:DsShearC}
%\quad
| G\circ(u - (B_j u)(\bs,\cdot))|
\leq
C |\bs - \bs'|^\alpha j^{1/2} 2^{-j}, \quad \forall \bs,\bs'\in \cS\,\mbox{ s.t. }  |\bs - \bs'|\ge 2^{-j}
\;.
\eeq
This assumption is deliberately optimistic because its consequences will later just serve as a benchmark for
the subsequent numerical realizations.
%
%\dw{[WD: I don't understand this condition because there is only an $\bs$ on the left hand side of the
%inequality while $\bs, \bs'$ appear on the right. Note that the right hand side vanishes for $\bs = \bs'$ which
%makes no sense.]}
%In view of $G(\cdot)\in L^2(D)' = L_2(D)$, %\dw{[WD: we usually denote the dual space by a prime.]}
%this may be rewritten as a stronger form of the error bound \eqref{eq:DsShearC}:
%there exists a constant $C>0$ such that, for $\alpha \in ]0,1]$ as
%in Assumption \ref{As:Cartoon} item (4), and for every $J \geq 1$
%holds
%
%\beq \label{eq:DsShearCLip}
%\| u - B_J^{x} u\|_{C^{0,\alpha}(\cS;L^2(D))} \leq C J^{1/2} 2^{-J}
%\eeq
%
%\dw{[WD: again I have a problem with this statement because it is identical to \eref{eq:DsShear}?]}
%Condition \eqref{eq:DsShearC} is, for the present parametric
%transport problem, analogous to the requirement of
%``mix\dw{ed}''-regularity in hyperbolic cross and sparse grid approximations.
%It corresponds to the convergence bound in Theorem \ref{ShearCart}
%in a stronger (w.r. to $\bs$) norm.

%The ensuing error bound requires
%only a slightly weaker condition than \eqref{eq:DsShearC},
%namely the validity of \eqref{eq:DsShearC} on a finite subset
%of all transport directions $\bs\in \cS$ which are required by the
%quadratures $I_\ell$ for $\ell = 0,1,...,J$ in \eqref{eq:I_ellErr}.

The quadrature approximation $I_L$ of the integral over $\cS$
is based on $N^s_L = O(2^L)$
(usually equi-spaced in $\cS$) directions $\bs$,
denoted by $\bs^\ell_i\in \cS$, $i=1,...,N^s_\ell = O(2^\ell)$.
We define
$$
\delta I_\ell := I_\ell - I_{\ell-1}\;,\quad \ell = 0,1,...
$$
with the convention $I_{-1} :=0$.
Then, for every $L\geq 1$,
$$
I_L = \sum_{\ell=0}^L \delta I_\ell \;,\quad \mbox{where}\quad
\delta I_\ell := (I_\ell - I_{\ell-1})
\;.
$$
We approximate $I(G(u))$ by
\beq \label{eq:SpApprox}
\cs{
I_L(G(u);\{j(\ell)\}_{\ell=0}^L)
:=
}
\sum_{\ell=0}^L ( \delta I_\ell(G(u_{j(\ell)})) )
\eeq
i.e., by quadrature-differences $\delta I_\ell$ applied
to (functionals of) anisotropically approximated solutions
in the spatial domain $D$, with spacial resolution level $j(\ell)$
to be determined. Note that each $\delta I_\ell$ requires
evaluating $G(u)$ at $O(2^\ell)$ distinct directions
$\{\bs^\ell_i\}_{i=1}^{N^s_\ell}$ with $N^s_\ell = O(2^\ell)$.
Each of these evaluations $G(u(\cdot,\bs^\ell_i))$ can be done
by approximating $u(\cdot,\bs^\ell_i)$ in $L^2(D)$. Then
evaluating $G( )$ on this approximation
(which we assume to be possible in $O(2^{j(\ell)})$ operations
%\dw{[WD: where have we assumed this?]}),
%
we subsequently form $\delta I_\ell$.
\cs{Since we are only interested in an integral
functional} $G( )$,
there is no need to merge partitions which are adapted to different directions.
We estimate the error
$$
E_L
:=
\left| I(G(u)) - \sum_{\ell=0}^L ( \delta I_\ell(G(u_{j(\ell)})) ) \right|\;,
$$
and in doing so, we obtain a choice for $j(\ell)$.
To this end, we write
$$
\begin{array}{rcl}
E_L & \leq & \displaystyle
\left| I(G(u)) - I_L(G(u)) \right|
+
\left|  I_L(G(u)) - \sum_{\ell=0}^L (\delta I_\ell)(G(u_{j(\ell)})) \right|
\;.
\end{array}
$$
The first term, denoted by $E^1_L$, is bounded by $\lesssim 2^{-\alpha L}$
provided that $\bs \mapsto G(u(\cdot,\bs)) \in C^{0,\alpha}(\cS)$
which is Assumption \ref{As:Cartoon}(3).
We estimate next the second term, denoted by $E^2_L$.

Using the linearity of $\delta I_\ell$ and of $G( )$,
$$
\begin{array}{rcl}
E^2_L & %\leq & \displaystyle
%\left| I_L(G(u)) - \sum_{\ell=0}^L (\delta I_\ell)(G(u_{j(\ell)})) \right|
%\\
%&
= & \displaystyle
\left| \sum_{\ell=0}^L (\delta I_\ell)(G(u)) -  (\delta I_\ell)(G(u_{j(\ell)})) \right|
\leq
\sum_{\ell=0}^L \left| (\delta I_\ell)(G(u) - G(u_{j(\ell)})) \right|
\\
& = &  \displaystyle
\sum_{\ell=0}^L \left| (\delta I_\ell)(G(u-u_{j(\ell)})) \right|
\;.
\end{array}
$$
We use the error bounds for the quadrature operator,
\eqref{eq:I_ellErr}, and the triangle inequality
to infer
$$
\begin{array}{rcl}
\left| (\delta I_\ell)(G(u-u_{j(\ell)})) \right|
&\leq& \displaystyle
\left| (I-I_\ell)(G(u-u_{j(\ell)})) \right|
+
\left|  (I-I_{\ell-1})(G(u-u_{j(\ell-1)})) \right|
\\
& \lesssim & \displaystyle
2^{-\alpha \ell} \| G(u-u_{j(\ell)}) \|_{C_{2^{-j(\ell)}}^{0,\alpha}(\cS)}
\\
%& \lesssim & \displaystyle
%2^{-\ell} \| G( u-u_{j(\ell)} ) \|_{\Lip(\cS)}
%\\
& \lesssim & \displaystyle
2^{-\alpha \ell} j(\ell)^{1/2} 2^{-j(\ell)},
\end{array}
$$
where we have invoked \eref{eq:DsShearC} in the last step.
%
%where we used the direction independence and %the linearity of $G(\cdot)$,
%\eqref{eq:DsShearC}.
%\dw{[WD: it seems to me that, instead of discussing/assuming the various inequaluties above,
%the essential property is
%\beqn
%\label{wish}
%\| G(u-u_{j(\ell)}) \|_{C^{0,\alpha}(\cS)}\lesssim   j(\ell)^{1/2} 2^{-j(\ell)}
%\eeqn
%which we should perhaps require/assume as that and exhibit a simple example where it can be verified. I am afraid, however,
%that this property is not true in general, see a short diescussion at the end of this section.]}
We now choose
\beq \label{eq:jl}
j(\ell) \simeq \alpha( L - \ell )\;, \quad L\geq 1\;,\;\; \ell = 0,1,...,L
\eeq
and  find
$$
E^2_L \leq
\sum_{\ell=0}^L \left| (\delta I_\ell)(G(u-u_{j(\ell)})) \right|
=
2^{-\alpha L} \sum_{\ell=0}^L ( L - \ell )^{1/2}
\lesssim
L^{3/2} 2^{-\alpha L}
\;.
$$
Combining with the bound $E^1_L \lesssim 2^{-\alpha L}$,
we arrive at the error bound
\beq\label{eq:SpTBd}
E_L \lesssim L^{3/2} 2^{-\alpha L}
\;.
\eeq
The total number of degrees of freedom/ total work
to construct the approximation \eqref{eq:SpApprox} is
bounded by
(assuming that $G(u_j)$ can be evaluated
\cs{to accuracy $O(2^{-\alpha L})$}
at a computational cost that stays proportional to
the number of degrees of freedom $N^x_j = O(2^j)$ involved in
building the \cs{direction-adaptive}
approximation \cs{$B^x_j u$ of $u$} in $L^2(D)$),
\beq\label{eq:WorkLBd}
W_L
\lesssim
\sum_{\ell = 0}^L 2^{\ell} 2^{j(\ell)}
\simeq
\sum_{\ell = 0}^L 2^{\ell} 2^{\alpha(L-\ell)}
\lesssim L^{\theta(\alpha)}
2^L
\;,
\eeq
with $\theta(\alpha) = 0$ if $0 < \alpha < 1$ and $\theta(1) = 1$.

This complexity bound, obtained under the assumption \eref{eq:DsShearC}, will serve
as a benchmark for the actual performance of such a sparse-tensor product scheme
when combined with the adaptive solver \solve
described in Section \ref{sec:shearlet-implementation}
and tested for a single convection field in Section \ref{sect:numerical}.

Combining \eqref{eq:SpTBd} and \eqref{eq:WorkLBd}, we obtain
essentially (up to \cs{logarithmic terms})
the same error versus work bound as for the approximation
of the solution in $L^2(D)$ for a single transport direction.
Then the approximations to $u(\cdot,\bs)$ will, of course, not
be those constructed in Theorem \ref{ShearCart},
but will be generated adaptively without prior knowledge on $u(\cdot,\bs)$,
\cs{and without certified (quasi-)optimality of the adaptive algorithm.}

Notice that this argument did not use the particular construction
in Theorem \ref{ShearCart}, and the same reasoning applies to any other directionally adaptive
approximation provided that it yields the above rates \eref{eq:DsShear} and \eref{eq:DsShearC}.
%\dw{[WD: with the infinity in $\bs$ bound]}.

Also remark that in \cite{KGCS2010a,KGCS2010b}, we obtained
approximation rates in $L^2(\cS\times D)$ which scaled
(up to log terms) analogous to approximation rates of
\emph{adaptively refined, isotropic} approximations in $D$;
these rates are, for piecwise linear polynomial approximations
in $D$, for parametric transport problems, inferior to the
rates %which are afforded by the directionally adapted approximations which are
furnished by the approximation result, Theorem \ref{ShearCart}.

The rates  in \cite{KGCS2010a,KGCS2010b} were possible due to
the adaptive refinements in $D$ being isotropic.
Applying the same reasoning as in \cite{KGCS2010a,KGCS2010b}
in the context of \cs{directionally adaptive approximations}
$B^x_j u$ of the solution $u(\cdot, \bs)$ to refinement level $j$
in the domain $D$ which is furnished
by Theorem \ref{ShearCart} to achieve an $L^2(D\times \cS)$-rate
could increase the asymptotic complexity of the approximation,
due to merging the \emph{non-nested, directionally adapted}
partitions in $D$ with directions sweeping all of $\cS$.
The evaluation of $G(B^x_ju(\cdot,\bs^\ell_i))$
for each ``quadrature node'' $\bs^\ell_i \in \cS$
which arises in the quadrature approximations obviates this merging
step, and affords the work bound \eqref{eq:WorkLBd}.
%\dw{[WD: I am not quite sure I understand the last paragraph. Isn't it true that even when using nested partitions, merging locally refined
%partions for different directions could double the complexity?]}\\
%
\subsection{Numerical Experiment}
\label{sec:NumEx}
We consider a model parametric transport
problem with non-zero inflow boundary data,
whose solution exhibits a discontinuity
along the \cs{$C^2$ singular support curve}
given by
$x_1 = \frac{\tan \theta}{2} x^2_2$ with $\theta \in [0,\frac{\pi}{4}]$.
More precisely, we consider the equation given by
\begin{equation}\label{eq:para}
Au =
\left(
  \begin{array}{c}
    \tan\theta x_2 \\
    1 \\
  \end{array}
\right)
 \cdot \nabla u + u = 1
\end{equation}
with boundary conditions
\[
g(x_1,x_2) = 1-x_1 \,\, \text{on} \,\, \{ (x_1,0) \in \Gamma_{-} : 0 < x_1 < 1 \}
\]
and
\[
g(x_1,x_2) = 0 \,\, \text{on} \,\, \{ (0,x_2) \in \Gamma_{-} : 0 < x_2 < 1 \}
\]
for  {$\bs = (\cos \theta, \sin \theta) \in \cS = {\mathbb{S}}^1$} with $\theta \in [0, \tfrac{\pi}{4}]$.
 {Writing $u(\theta)= u(\bs,\cdot)$,}  we approximate
\[I(G(u)) := \int_{[0,\frac{\pi}{4}]}G(u(\cdot,\theta))d\theta\]
by
\begin{equation}\label{eq:SpApprox2}
\hat{I}_{L}(G(u)) = \sum_{\ell = 0}^{L}I_{\ell}(G(u_{j(\ell)}))-I_{\ell-1}(G(u_{j(\ell)})),
\end{equation}
where $I_{\ell}$ is quadrature approximation of the integral over $[0,\frac{\pi}{4}]$ with $2^{\ell}$ equi-spaced samples in $[0,\frac{\pi}{4}]$ for $\ell \ge 0$ and $I_{-1} = 0$. That is
\begin{equation}\label{eq:quad}
I_{\ell}(G(u_{j(\ell)})) = \frac{\pi}{2^{\ell+2}}\sum_{i = 0}^{2^{\ell}-1}G(u_{j(\ell)}(\cdot,\theta^{\ell}_i)) \quad \text{and} \quad
I_{\ell-1}(G(u_{j(\ell)})) = \frac{\pi}{2^{\ell+1}}\sum_{i = 0}^{2^{\ell-1}-1}G(u_{j(\ell)}(\cdot,\theta^{\ell-1}_i))
\end{equation}
where $j(\ell) = L-\ell$, $\theta^{\ell}_i = \frac{i\pi}{2^{\ell+2}}$ and $u_{j(\ell)}(\cdot,\theta^{\ell}_i)$ is an approximate solution to \eqref{eq:para} for a fixed $\theta^{\ell}_i \in [0,\pi/4]$ obtained by our adaptive refinement solver $\solve$ with $O(2^{j(\ell)})$ adpatively chosen basis elements. Here, the linear functional $G \in (L_2(D))^*$ is given as
\beqn
\label{average}
G(u) = \int_{D_0} u(x) dx \quad \text{where} \,\,
D_0 = [\tfrac{1}{8},\tfrac{5}{8}]\times [\tfrac{1}{4}, \tfrac{3}{4}]
\;.
\eeqn
Note that we choose $D_0 \subset D$ so that it \cs{contains}
a region where the jump discontinuity of a solution $u$ for \eqref{eq:para}
occurs for each $\theta \in [0,\tfrac{\pi}{4}]$.

\cs{
Assuming available a
directionally adaptive Petrov-Galerkin solver $\solve$ in the physical domain,
\eqref{eq:SpApprox} is implemented as follows: we assume that
subroutine $\solve$ requires $N_{\ell}$,
an upper bound on the number of basis elements to be activated
by $\solve$ and the directional angle $\theta \in [0,\tfrac{\pi}{4}]$
as input arguments, to compute $\hat{I}_L$ from \eqref{eq:SpApprox} and
\eqref{eq:jl}.
}

\begin{algorithm}[htb]
  \caption{}
  \label{alg:numeric_para}
  \begin{algorithmic}[1]
  \State \State Initialization: \css{Choose a maximal refinement level $L>0$.}
  \State  Compute $u_{j(0)}(\cdot,0) := \solve(N_{0},0)$ with $N_{0} = O(2^L)$.
  \State Compute $G(u_{j(0)}(\cdot,0))$.
  \State Apply quadrature approximation \eqref{eq:quad} to compute $I_{0}(G(u_{j(0)}))$.
  \State $\hat I_L \leftarrow  I_{0}(G(u_{j(0)}))$.
    \For{$\ell = L,\dots,1$}
      \State  Compute $u_{j(\ell)}(\cdot,\theta^{\ell}_i) := \solve(N_{\ell},\theta^{\ell}_i)$ for $i = 0,\dots,2^{\ell}-1$ and $N_{\ell} = O(2^{L-\ell})$.
      \State Compute $G(u_{j(\ell)}(\cdot,\theta^{\ell}_i))$ for $i = 0,\dots,2^{\ell}-1$.
      \State Apply quadrature approximation \eqref{eq:quad} to compute $I_{\ell}(G(u_{j(\ell)}))-I_{\ell-1}(G(u_{j(\ell)}))$.
      \State $\hat{I}_L \leftarrow \hat{I}_L + I_{\ell}(G(u_{j(\ell)}))-I_{\ell-1}(G(u_{j(\ell)}))$.
      \EndFor
  \end{algorithmic}
\end{algorithm}

\begin{figure}
\begin{center}
\includegraphics[width=.4\textwidth]{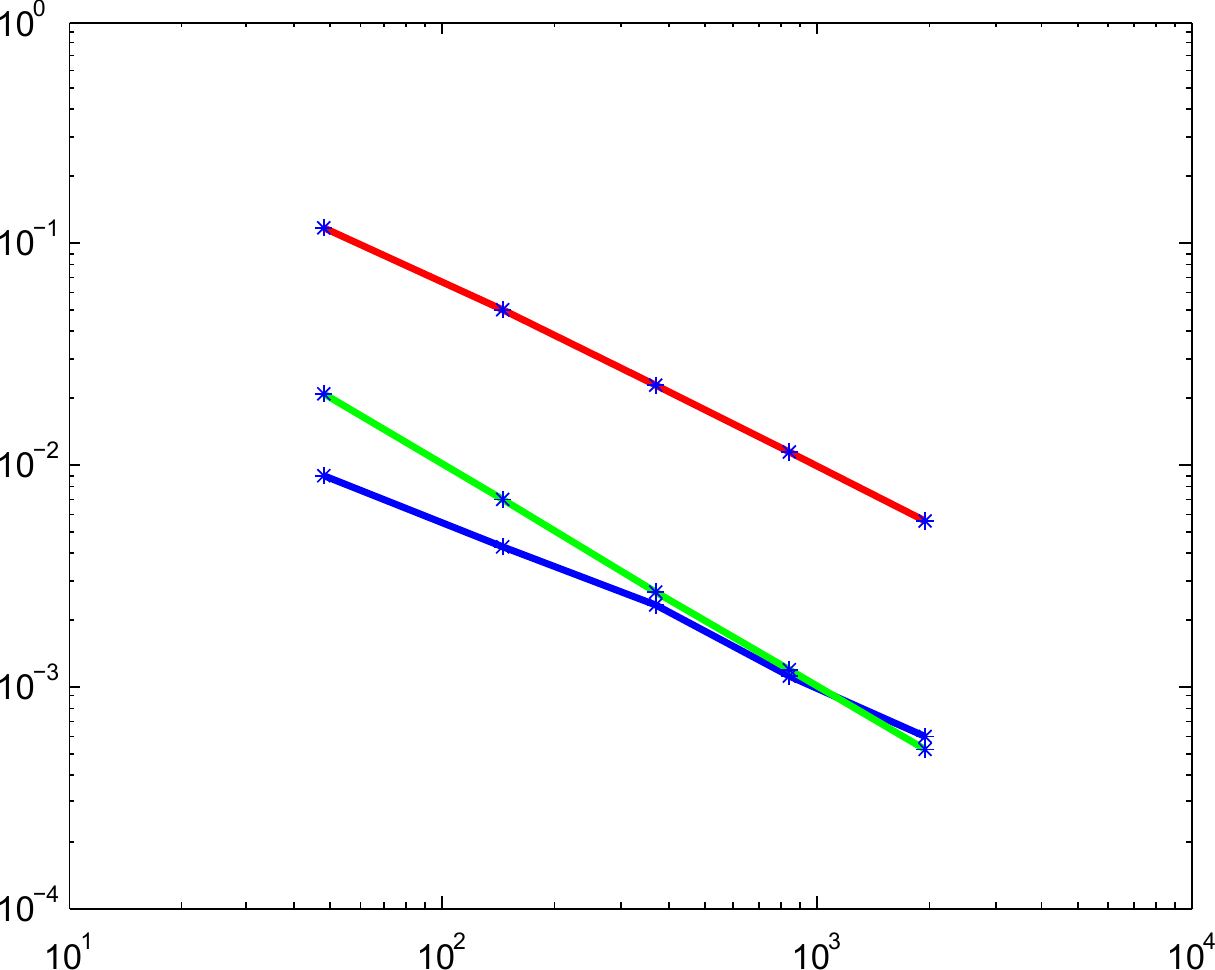}
\end{center}
%\vspace{-50pt}
\caption{\cs{Estimated approximation error in $L^2(D)$ versus total number of
             degrees of freedom (blue), and}
             the benchmark rates $n^{-1}$ (green) and $\log_2(n)n^{-1}$ (red).}
\label{fig:SPTtest}
\end{figure}

Figure \ref{fig:SPTtest} and Table \ref{fig:table2}
show approximation errors and their comparison
with the benchmark rates $n^{-1}$ and $\log_2(n)n^{-1}$ for
\[E_L = |I(G(u))-\hat{I}_L(G(u))|.\]
\begin{table}
\begin{center}
\begin{tabular}{|c|c|c|}
  \hline
  \text{Dof:} $n$ & $n^{-1}$ ($\log_2(n)n^{-1}$) & $E_L$ \\
  \hline
  48 &  0.020833 (0.116353)& 0.012404 \\
  144 & 0.006944 (0.049791)& 0.005242 \\
  372 & 0.002688 (0.022954)& 0.002407 \\
  849 & 0.001177 (0.003613)& 0.001266 \\
  1944& 0.000514 (0.001842)& 0.000740 \\
  %345 & 0.183383 & 0.023095 \\
  \hline
\end{tabular}
\end{center}
\caption{Numerical estimates for the benchmark rates $n^{-1}$ and $\log_2(n)n^{-1}$ and the approximation error.}
\label{fig:table2}
\end{table}

%\vspace{-50pt}

\bibliographystyle{alpha}
\newcommand{\etalchar}[1]{$^{#1}$}
%\bibliography{ref}
%%%%%%%%%%%%%%%%%%%%%%%%%%%%%%%%%

%

\end{document}